\documentclass[12pt,a4paper]{article}

\def\?{
?\vadjust{\vbox to 0pt{\vss\hbox{\kern\hsize\kern1em\large\bf ?!}}}}

\usepackage{pifont}
\usepackage{amsfonts}
\usepackage{mathrsfs}
\usepackage{amssymb}
\usepackage{amsmath}

\usepackage{enumerate}

\newtheorem{theorem}{Theorem}[section]
\newtheorem{definition}[theorem]{Definition}
\newtheorem{lemma}[theorem]{Lemma}
\newtheorem{corollary}[theorem]{Corollary}

\newtheorem{example}[theorem]{Example}

\makeatletter

\newcommand{\Rmnum}[1]{\expandafter\@slowromancap\romannumeral #1@}
\makeatother

\newcommand\cm{\operatorname{cm}}
\newcommand\lcm{\operatorname{lcm}}
\newcommand\supp{\operatorname{supp}}
\newcommand\Mod{\operatorname{Mod}}
\newcommand\sgp{sgp}
\newcommand\fir{\operatorname{fir}}
\newcommand\path{\operatorname{path}}
\newcommand\wt{\operatorname{wt}}
\newcommand\inwt{\operatorname{inwt}}

\begin{document}

\title{Gr\"obner--Shirshov bases and their
calculation\footnote{Supported by
the NNSF of China (11171118),
the Research Fund for the Doctoral Program of Higher Education of China
(20114407110007),
the NSF of Guangdong Province (S2011010003374)
and the Program on International Cooperation and Innovation,
Department of Education,
Guangdong Province (2012gjhz0007).}}
\author{
L. A. Bokut\footnote {Supported by
RFBR 12-01-00329,
LSS--3669.2010.1,
SB RAS Integration grant No.2009.97 (Russia)
and Federal Target Grant
``Scientific and educational personnel of innovation Russia''
for 2009-2013
(government contract No.02.740.11.5191).} \\
{\small Sobolev Institute of Mathematics
and Novosibirsk State University, Novosibirsk 630090, Russia}\\
{\small \ School of Mathematical Sciences,
South China Normal University}\\
{\small Guangzhou 510631, P. R. China}\\
{\small  bokut@math.nsc.ru}\\
\\
Yuqun Chen\\
{\small \ School of Mathematical Sciences,
South China Normal University}\\
{\small Guangzhou 510631, P. R. China}\\
{\small yqchen@scnu.edu.cn}}

\date{}

\maketitle
\noindent\textbf{Abstract:}
In this survey we give an~exposition of
the theory of Gr\"{o}bner--Shirshov bases
for associative algebras and Lie algebras.
We mention some new Composition-Diamond lemmas and applications.



\ \

\ \

\newpage
{\Large Notation}

\ \

CD-lemma:
Composition-Diamond lemma.

GS basis:
Gr\"{o}bner--Shirshov basis.

LS word (basis):
Lyndon--Shirshov word (basis).

ALSW(X):
the set of all associative Lyndon--Shirshov words in~%
$X$.

NLSW(X):
the set of all non-associative Lyndon--Shirshov words in~%
$X$.

PBW theorem:
the Poincare-Birkhoff-Witt theorem.

$X^*$:
the free monoid generated by~%
$X$.

$[X]$:
the free commutative monoid generated by~%
$X$.

$X^{\ast \ast }$:
the set of all non-associative words
$(u)$
in~%
$X$.

$gp\langle X|S\rangle$:
the group generated by~%
$X$
with defining relations~%
$S$.

$\sgp\langle X|S\rangle$:
the semigroup generated by~%
$X$
with defining relations~%
$S$.

$k$:
a~field.

$K$:
a~commutative algebra over~%
$k$
with unity.

$k\langle X\rangle$:
the free associative algebra over~%
$k$
generated by~%
$X$.

$k\langle X|S\rangle$:
the associative algebra over~%
$k$
with generators~%
$X$
and defining relations~%
$S$.

$S^c$:
a Gr\"obner--Shirshov completion  of~%
$S$.

$Id(S)$:
the ideal generated by a~set~%
$S$.

$\bar s$:
the maximal word of a~polynomial~%
$s$
with respect to some ordering~%
$<$.

$Irr(S)$:
the set of all monomials avoiding the subword~%
$\bar s$
for all
$s\in S$.

$k[ X]$:
the polynomial algebra over~%
$k$
generated by~%
$X$.

$Lie(X)$:
the free Lie algebra over~%
$k$
generated by~%
$X$.

$Lie_K(X)$:
the free Lie algebra generated by~%
$X$
over a~commutative algebra~%
$K$.

\section{Introduction}

In this survey we review the method of Gr\"obner--Shirshov%
\footnote{Though Shirshov \cite{Sh62b} 1962 was the first
to come up with the idea of a~`Gr\"obner--Shirshov basis'
for Lie and non-commutative polynomial algebras,
his paper became practically unknown outside Russia.
In the meantime,
Buchberger's `Gr\"obner basis'
(Thesis  1965 \cite{bu65}, paper 1970 \cite{bu70})
for (commutative) polynomials
became very popular in science.
As a~result,
the first author suggested the name `Gr\"obner--Shirshov basis'
for non-commutative and non-associative polynomials.
For (commutative) differential polynomials
an~analogous,
or better to say,
closely related `basis' is called a~Ritt--Kolchin characteristic set,
due to Ritt \cite{Ritt1950} 1950 and Kolchin \cite{Kolchin1973} 1973,
and rediscovered by Wen-Ts\"un Wu \cite{Wu1978} 1978.}
(GS for short)
bases for different classes of linear universal algebras,
together with an~overview of 
calculation of these bases
in a~variety of specific cases.

A.~I.~Shirshov
(also spelled  A.~I.~\v{S}ir\v{s}ov)
in his
pioneering work (\cite{Sh62b}, 1962)
posed  the following fundamental question:

How to find a~linear basis of a~Lie algebra
defined by generators and relations?

He gave an~infinite algorithm to solve this problem
using a~new notion of the composition
(later the `$s$-polynomial'  in Buchberger's terminology  \cite{bu65,bu70})
of two Lie polynomials
and a~new notion of completion of a~set of Lie polynomials
(adding nontrivial compositions;
the critical pair/completion
(cpc-) algorithm in the later terminology of
Knuth and Bendix \cite{Knuth-Bendix}
and Buchberger \cite{Bu87,BuCL}).

Shirshov's algorithm goes as follows.
Consider a~set~%
$S \subset Lie(X)$
of Lie polynomials in the free algebra
$k\langle X\rangle$
on~%
$X$
over a~field~%
$k$
(the algebra of non-commutative polynomials on~%
$X$
over~%
$k$).
Denote by
$S'$
the superset of~%
$S$
obtained by adding all non-trivial Lie compositions
(`Lie~%
$s$-polynomials')
of the elements of~%
$S$.
The problem of triviality of a~Lie polynomial
modulo a~finite
(or recursive)
set~%
$S$
can be solved algorithmically
using Shirshov's Lie reduction algorithm
from his previous paper \cite{Sh58}, 1958.
In general,
an~infinite sequence
$$
S\subseteq S'\subseteq S''\subseteq\dots \subseteq S^{(n)}\subseteq\dots
$$
of Lie multi-compositions arises.
The union
$S^c$
of this sequence
has the property that
every Lie composition of elements of
$S^c$
is trivial modulo
$S^c$.
This is what is now called a~Lie GS basis.

Then a~new `Composition-Diamond lemma%
\footnote{The name
`Composition-Diamond lemma'
combines the Neuman Diamond Lemma \cite{Newman},
the Shirshov Composition Lemma \cite{Sh62b}
and the Bergman Diamond Lemma \cite{Be78}.}
for Lie algebras'
(Lemma 3 in \cite{Sh62b})
implies that
the set
$Irr(S^c)$
of all
$S^c$-irreducible
(or
$S^c$-reduced)
basic Lie monomials
$[u]$
in~%
$X$
is a~linear basis of the Lie algebra
$Lie(X|S)$
generated by~%
$X$
with defining relations~%
$S$.
Here a~basic Lie monomial means
a~Lie monomial in a~special linear basis of the free Lie algebra
$Lie(X)\subset k\langle X\rangle$,
known as the Lyndon--Shirshov
(LS for short)
basis
(Shirshov \cite{Sh62b} and Chen-Fox-Lyndon \cite{CFL},
see below).
An~LS monomial
$[u]$
is called
$S^c$-irreducible
(or
$S^c$-reduced)
whenever~%
$u$,
the associative support of
$[u]$,
avoids the word~%
$\bar s$
for all
$s\in S$,
where~%
$\bar s$
is the maximal word of~%
$s$
as an~associative polynomial
(in the deg-lex ordering).
To be more precise,
Shirshov used his reduction algorithm at each step
$S$,
$S'$,
$S'',\dots$.
Then we have a~direct system
$S\rightarrow S'\rightarrow S''\rightarrow\dots$
and
$S^c=\underrightarrow{lim}S^{(n)}$
is what is now called a~\textit{minimal GS basis}
(a minimal GS basis is not unique,
but a~reduced GS basis is,
see below).
As a~result,
Shirshov's algorithm
gives a~solution to the above problem for Lie algebras.

Shirshov's algorithm,
dealing with the word problem,
is an~infinite algorithm
like the Knuth--Bendix algorithm \cite{Knuth-Bendix}, 1970
dealing with the identity problem
for every variety of universal algebras%
\footnote{We use the standard algebraic terminology
`the word problem',
`the identity problem',
see O.~G.~Kharlampovich, M.~V.~Sapir \cite{Sapir94}
for instance.}.
The initial data for the Knuth--Bendix algorithm
is the defining identities of a~variety.
The output of the algorithm,
if any,
is a~`Knuth--Bendix basis' of identities of the variety
in the class of all universal algebras of a~given signature
(not a~GS basis of defining relations,
say,
of a~Lie algebra).

Shirshov's algorithm gives linear bases
and algorithmic decidability of the word problem
for one-relation Lie algebras \cite{Sh62b},
(recursive)
linear bases for Lie algebras with
(finite)
homogeneous defining relations \cite{Sh62b},
and linear bases for free products of Lie algebras
with known linear bases \cite{Sh62d}.
He also proved the Freiheitssatz
(freeness theorem)
for Lie algebras \cite{Sh62b}
(for every one-relation Lie algebra
$Lie(X|f)$,
the subalgebra
$\langle X\backslash\{x_{i_0}\}\rangle$,
where
$x_{i_0}$
appears in
$f$,
is a~free Lie algebra).
The Shirshov problem \cite{Sh62b}
of the decidability of the word problem for Lie algebras
was solved negatively in \cite{Bo72}.
More generally,
it was proved \cite{Bo72} that
some recursively presented Lie algebras with undecidable word problem
can be embedded into finitely presented Lie algebras
(with undecidable word problem).
It is a~weak analogue
of the Higman embedding theorem for groups \cite{Higman}.
The problem \cite{Bo72}
whether an~analogue of the Higman embedding theorem
is valid for Lie algebras
is still open.
For associative algebras
a~similar problem \cite{Bo72}
was solved positively by V.~Y.~Belyaev \cite{Belyaev}.
A~simple example of a~Lie algebra
with undecidable word problem was given by G.~P.~Kukin \cite{Kukin}.

Actually,
a~similar algorithm for associative algebras
is implicit in Shirshov's paper \cite{Sh62b}.
The reason is that he treats
$Lie(X)$
as the subspace of Lie polynomials
in the free associative algebra
$k\langle X\rangle$.
Then to define a~Lie composition
$\langle f,g\rangle_w$
of two Lie polynomials
relative to an~associative word
$w=\lcm(\bar f,\bar g)$,
he defines firstly the associative composition
(non-commutative `$s$-polynomial')
$(f,g)_w=fb-ag$,
with
$a,b\in X^*$.
Then he inserts some brackets
$\langle f,g\rangle_w=[fb]_{\bar f}-[ag]_{\bar g}$
by applying his special bracketing lemma of \cite{Sh58}.
We can obtain
$S^c$
for every
$S\subset k\langle X\rangle$
in the same way as for Lie polynomials
and in the same way as for Lie algebras
(`CD-lemma for associative algebras')
to infer that
$Irr(S^c)$
is a~linear basis of the associative algebra
$k\langle X|S\rangle$
generated by~%
$X$
with defining relations~%
$S$.
All proofs are similar to those in \cite{Sh62b} but much easier.

Moreover,
the cases of semigroups and groups
presented by generators and defining relations
are just special cases of associative algebras
via semigroup and group algebras.
To summarize,
Shirshov's algorithm gives linear bases and normal forms of elements
of every Lie algebra,
associative algebra,
semigroup or group presented by generators and defining relations!
The algorithm works in many cases
(see below).

The theory of Gr\"obner bases and Buchberger's algorithm
were initiated
by B.~Buchberger (Thesis \cite{bu65} 1965, paper \cite{bu70} 1970)
for commutative associative algebras.
Buchberger's algorithm is a~finite algorithm
for finitely generated commutative algebras.
It is one of the most useful and famous algorithms
in modern computer science.

Shirshov's paper \cite{Sh62b}
was in the spirit of the program of A.~G.~Kurosh (1908-1972)
to study non-associative
(relatively)
free algebras and free products of algebras,
initiated in Kurosh's paper \cite{Kurosh}, 1947.
In that paper he proved non-associative analogs
of the Nielsen-Schreier and Kurosh theorems for groups.
It took quite a~few years to clarify the situation
for Lie algebras in Shirshov's papers \cite{Sh53}, 1953
and \cite{Sh62b}, 1962
closely related to his theory of GS bases.
It is important to note that
Kurosh's program quite unexpectedly led to
Shirshov's theory of GS bases
for Lie and associative algebras \cite{Sh62b}.

A~step in Kurosh's program was made by
his student A.~I.~Zhukov in his Ph.D.\ Thesis \cite{Zhukov}, 1950.
He algorithmically solved the word problem for non-associative algebras.
In a~sense,
it was the beginning of the theory of GS bases
for non-associative algebras.
The main difference with the future approach of Shirshov
is that
Zhukov did not use a~linear ordering of non-associative monomials.
Instead he chose an~arbirary monomial of maximal degree
as the `leading' monomial of a~polynomial.
Also,
for non-associative algebras
there is no `composition of intersection'
(`$s$-polynomial').
In this sense it cannot be a~model for Lie and associative algebras%
\footnote{After his Ph.D.\ Thesis of 1950,
A.~I.~Zhukov moved to the present
Keldysh Institute of Applied Mathematics (Moscow)
to do computational mathematics.
S.~K.~Godunov in `Reminiscence about numerical schemes',
arxiv.org/pdf/0810.0649, 2008,
mentioned his name in relation to the creation of
the famous Godunov numerical method.
So,
A.~I.~Zhukov was a~forerunner of two important computational
methods!}.

A.~I.~Shirshov,
also a~student of Kurosh's,
defended his Candidate of Sciences Thesis \cite{S53a}
at Moscow State University in 1953.
His thesis together with the paper
that followed \cite{Sh58}, 1958
may be viewed as a~background for
his later method of GS bases.
In the thesis,
he proved the free subalgebra theorem for free Lie algebras
(now known as Shirshov-Witt theorem,
see also Witt \cite{Witt56}, 1956)
using the elimination process
rediscovered by Lazard \cite{Lazard60}, 1960.
He used the elimination process later \cite{Sh58}, 1958
as a~general method to prove the properties of regular (LS) words,
including an~algorithm of (special) bracketing of an~LS word
(with a~fixed LS subword).
The latter algorithm is of some importance
in his theory of GS bases for Lie algebras
(particularly in the definition of the composition of two Lie polynomials).
Shirshov also proved the free subalgebra theorem
for (anti-) commutative non-associative algebras \cite{Sh54}, 1954.
He used that later in \cite{Sh62a}, 1962
for the theory of GS bases of
(commutative,
anti-commutative)
non-associative algebras.
Shirshov (Thesis \cite{S53a}, 1953)
found the
(`Hall--Shirshov')
series of bases of a~free Lie algebra
(see also \cite{Sh62c} 1962,
the first issue of Malcev's Algebra and Logic).%
\footnote{It must be pointed out that A.~I.~Malcev (1909-1967)
inspired Shirshov's works very much.
Malcev was an~official opponent
(referee)
of his (second) Doctor of Sciences Dissertation at MSU in 1958.
The first author,
L.~A.~Bokut,
remembers this event at the Science Council Meeting,
chaired by A.~N.~Kolmogorov,
and Malcev's words ``Shirshov's dissertation is a~brilliant one!''.
Malcev and Shirshov worked together
at the present Sobolev Institute of Mathematics in Novosibirsk
since 1959
until Malcev's sudden death at 1967,
and have been friends despite the age difference.
Malcev headed the Algebra and Logic Department
(by the way,
the first author is a~member of the departement since 1960)
and Shirshov was the first deputy director of the institute
(whose director was S.~L.~Sobolev).
In those years,
Malcev was interested in
the theory of algorithms of mathematical logic and
algorithmic problems of model theory.
Thus,
Shirshov had an~additional motivation to work
on algorithmic problems for Lie algebras.
Both Maltsev and Kurosh were delighted with
Shirshov's results of \cite{Sh62b}.
Malcev successfully nominated the paper
for an~award of the Presidium of
the Siberian Branch of the Academy of Sciences
(Sobolev and Malcev were the only Presidium members
from the Institute of Mathematics at the time).}

The LS basis is a~particular case of
the Shirshov or Hall--Shirshov series of bases
(cf.\ Reutenauer \cite{Reutenauer},
where this series is called the `Hall series').
In the definition of his series,
Shirshov used Hall's inductive procedure
(see Ph.~Hall \cite{PH}, 1933, M.~Hall \cite{MH}, 1950):
a~non-associative monomial
$w=((u)(v))$
is a~basic monomial whenever
\begin{enumerate}[(1)]
\item
$(u),\ (v)$
are basic;
\item
$(u)>(v)$;
\item
if
$(u)=((u_1)(u_2))$
then
$(u_2)\leq(v)$.
\end{enumerate}
However,
instead of ordering by the degree function
(Hall words),
he used an~arbirary linear ordering of non-associative monomials
satisfying
$$
((u)(v))>(v).
$$
For example,
in his Thesis \cite{S53a}, 1953
he used the ordering by the content of monomials
(the \textit{content} of,
say,
the monomial
$(u) = ((x_2x_1)((x_2x_1)x_1))$
is the vector
$(x_2,x_2,x_1,x_1,x_1)$).
Actually,
the content
$\widehat{u}$
of
$(u)$
may be viewed as a~commutative associative word that
equals
$u$
in the free commutative semigroup.
Two contents are compared lexicographically
(a proper prefix of a~content is greater than the content).

If we use the lexicographic ordering,
$(u)\succ(v)$
if
$u\succ v$
lexicographically
(with the condition
$u\succ uv,\ v\neq1$),
then we obtain the LS basis.%
\footnote{The Lyndon--Shirshov basis for the alphabet
$x_1,x_2$
is different from the above Shirshov content basis
starting with monomials of degree~7.}
For example,
for the alphabet
$x_1$,
$x_2$
with
$x_2\succ x_1$
we obtain basic Lyndon--Shirshov monomials by induction:
\begin{gather*}
x_2,\ x_1,\
[x_2x_1],\
[x_2[x_2x_1]]=[x_2x_2x_1],\
[[x_2x_1]x_1]=[x_2x_1x_1],\\
[x_2[x_2x_2x_1]]=[x_2x_2x_2x_1],\
[x_2[x_2x_1x_1]]=[x_2x_2x_1x_1],\\
[[x_2x_1x_1]x_1]=[x_2x_1x_1x_1],\
[[x_2x_1][x_2x_1x_1]]=[x_2x_1x_2x_1x_1],
\end{gather*}
and so on.
They are exactly all Shirshov regular (LS) Lie monomials
and their associative supports are exactly all Shirshov regular words
with a~one-to-one correspondence between two sets
given by the Shirshov elimination (bracketing) algorithm
for (associative) words.

Let us recall that
an~elementary step of Shirshov's elimination algorithm
is to join the minimal letter of a~word to previous ones by bracketing
and to continue this process
with the lexicographic ordering of the new alphabet.
For example,
suppose that
$x_2\succ x_1$.
Then we have the succession of bracketings
\begin{gather*}
x_2x_1x_2x_1x_1x_2x_1x_1x_1x_1x_2x_1x_1,\\
[x_2x_1][x_2x_1x_1][x_2x_1x_1x_1][x_2x_1x_1],\\
[x_2x_1][[x_2x_1x_1][x_2x_1x_1x_1]][x_2x_1x_1],\\
[[x_2x_1][[x_2x_1x_1][x_2x_1x_1x_1]]][x_2x_1x_1],\\
[[[x_2x_1][[x_2x_1x_1][x_2x_1x_1x_1]]][x_2x_1x_1]];\\
x_2x_1x_1x_1x_2x_1x_1x_2x_1x_2x_2x_1,\\
[x_2x_1x_1x_1][x_2x_1x_1][x_2x_1]x_2[x_2x_1],\\
[x_2x_1x_1x_1][x_2x_1x_1][x_2x_1][x_2[x_2x_1]];\\
x_2x_1x_1x_1\prec x_2x_1x_1\prec x_2x_1\prec x_2x_2x_1.
\end{gather*}
By the way, the second series of partial bracketings
illustrates Shirshov's factorization theorem \cite{Sh58} of 1958 that
every word is a~non-decreasing product of LS words
(it is often mistakenly called Lyndon's theorem,
see \cite{Berstel-J-P}).

The Shirshov special bracketing \cite{Sh58} goes as follows.
Let us give as an~example
the special bracketing of the LS word
$w=x_2x_2x_1x_1x_2x_1x_1x_1$
with the LS subword
$u=x_2x_2x_1$.
The Shirshov standard bracketing is
$$
[w]=[x_2[[[x_2x_1]x_1][x_2x_1x_1x_1]]].
$$
The Shirshov special bracketing is
$$
[w]_u=[[[u]x_1][x_2x_1x_1x_1]].
$$
In general,
if
$w=aub$
then the Shirshov standard bracketing gives
$[w]=[a[uc]d]$,
where
$b=cd$.
Now,
$c=c_1\cdots c_t$,
each
$c_i$
is an~LS-word,
and
$c_1\preceq\cdots\preceq c_t$
in the lex ordering
(Shirshov's factorization theorem).
Then we must change the bracketing of
$[uc]$:
$$
[w]_u=[a[\dots[[u][c_1]]\dots[c_t]]d]
$$
The main property of
$[w]_u$
is that
$[w]_u$
is a~monic associative
polynomial with the maximal monomial
$w$;
hence,
$\overline{[w]_u}=w$.

Actually,
Shirshov \cite{Sh62b}, 1962 needed 
a~`double' relative bracketing of a~regular word
with two disjoint LS subwords.
Then he implicitly used the following property:
every LS subword of
$c=c_1\cdots c_t$
as above is a~subword of some
$c_i$
for
$1\leq i\leq t$.

Shirshov defined regular (LS) monomials \cite{Sh58}, 1958,
as follows:
$(w)=((u)(v))$
is a~regular monomial
iff:
\begin{enumerate}[(1)]
\item
$w$
is a
regular word;
\item
$(u)$
and
$(v)$
are regular monomials
(then automatically
$u\succ v$
in the lex ordering);
\item
if
$(u)=((u_1)(u_2))$
then
$u_2\preceq v$.
\end{enumerate}

Once again,
if we formally omit all Lie brackets in Shirshov's paper \cite{Sh62b}
then essentially the same algorithm and essentially the same CD-lemma
(with the same but much simpler proof)
yield a~linear basis for associative algebra
presented by generators and defining relations.
The differences are the following:
\begin{itemize}
\item
no need to use LS monomials and LS words,
since the set
$X^*$
is a~linear basis of the free associative algebra
$k\langle X\rangle$;
\item
the definition of associative composition for monic polynomials
$f$~%
and~%
$g$,
$$
(f,g)_w = fb-ag, \quad
w=\bar{f}b=a\bar{g}, \quad
deg(w)<deg(\bar f)+deg(\bar g),
$$
or
$$
(f,g)_w=f-agb, \quad
w=\bar{f}=a\bar{g}b, \quad
w,a,b \in X^*,
$$
are much simpler than
the definition of Lie composition for monic Lie polynomials
$f$~%
and~%
$g$,
$$
\langle f,g\rangle_w = [fb]_{\bar f}-[ag]_{\bar g}, \quad
w=\bar{f}b=a\bar{g}, \quad
deg(w)<deg(\bar f)+deg(\bar g),
$$
or
$$
\langle f,g\rangle_w = f-[agb]_{\bar g}, \quad
w=\bar{f}=a\bar{g}b, \quad
w,a,b, \bar{f}, \bar{g} \in X^*,
$$
where
$[fb]_{\bar f}$,
$[ag]_{\bar g}$,
and
$[agb]_{\bar g}$
are the Shirshov special bracketings of the LS words~%
$w$
with fixed LS subwords
$\bar{f}$~%
and~%
$\bar{g}$
respectively.
\item
The definition of elimination of the leading word~%
$\bar{s}$
of an~associative monic polynomial~%
$s$
is straightforward:
$a\bar{s}b\rightarrow a(r_s)b$
whenever
$s=\bar{s} - r_s$
and
$a,b\in X^*$.
However,
for Lie polynomials,
it is much more involved
and uses the Shirshov special bracketing:
$f\rightarrow f-[agb]_{\bar g}$
whenever
$\bar{f}=a\bar{g}b$.
\end{itemize}

We can formulate the main idea of Shirshov's proof as follows.
Consider a~complete set
$S$
of monic Lie polynomials
(all compositions are trivial).
If
$w=a_1\bar {s_1}b_1= a_2\bar {s_2}b_2$,
where~%
$w$,
$a_i$,
$b_i \in X^*$
and~%
$w$
is an~LS word,
while
$s_1,s_2\in S$,
then the Lie monomials
$[a_1s_1b_1]_{\overline{s_1}}$
and
$[a_2s_2b_2]_{\overline{s_2}}$
are equal modulo the smaller Lie monomials in
$Id(S)$:
$$
[a_1s_1b_1]_{\overline{s_1}} = [a_2s_2b_2]_{\overline{s_2}} +
\sum_{i>2} \alpha_i[a_is_ib_i]_{\overline{s_i}},
$$
where
$\alpha_i\in k,\ s_i\in S$
and
$\overline {[a_is_ib_i]_{\overline{s_i}}}=a_i\bar {s_i}b_i < w$.
Actually,
Shirshov proved a~more general result:
if
$\overline{(a_1s_1b_1)} = a_1\overline{s_1}b_1$
and
$\overline{(a_2s_2b_2)}= a_2\overline{s_2}b_2$
with
$w=a_1\overline{s_1}b_1=a_2\overline{s_2}b_2$
then
$$
(a_1s_1b_1) = (a_2s_2b_2) + \sum_{i>2} \alpha_i(a_is_ib_i),
$$
where
$\alpha_i\in k,\ s_i\in S$
and
$\overline {(a_is_ib_i)}=a_i\bar {s_i}b_i < w$.
Below we call a~Lie polynomial
$(asb)$
a~Lie normal~%
$S$-word
provided that
$\overline {(asb)}=a\bar {s}b$.

This is precisely where he used
the notion of composition and other notions and properties mentioned above.

It is much easier to prove an~analogue of this property
for associative algebras
(as well as commutative associative algebras):
given a~complete monic set~%
$S$
in
$k\langle X\rangle$
($k[X]$),
for
$w=a_1\bar{s_1}b_1 = a_2\bar {s_2}b_2$
with
$a_i,b_i\in X^*$
and
$s_1,s_2\in S$
we have
$$
a_1s_1b_1 = a_2s_2b_2 + \sum_{i>2} \alpha_i a_is_ib_i,
$$
where
$\alpha_i\in k,\ s_i\in S$
and
$a_i\bar {s_i}b_i < w$.

Summarizing,
we can say with confidence that the work (Shirshov
\cite{Sh62b}) implicitly contains the CD-lemma for associative
algebras as a~simple exercise that requires no new ideas.
The first author,
L.~A.~Bokut,
can confirm that
Shirshov clearly understood this and told him that
``the case of associative algebras is the same''.
The lemma was formulated explicitly in Bokut \cite{Bo76}, 1976
(with a~reference to Shirshov's paper \cite{Sh62b}),
Bergman \cite{Be78}, 1978,
and Mora \cite{Mora}, 1986.

Let us emphasize once again that
the CD-Lemma for associative algebras
applies to every semigroup
$P=\sgp\langle X|S\rangle$,
and in particular to every group,
by way of the semigroup algebra
$kP$
over a~field~%
$k$.
The latter algebra has the same generators and defining relations as~%
$P$,
or
$kP=k\langle X|S\rangle$.
Every composition of the binomials
$u_1-v_1$
and
$u_2-v_2$
is a~binomial
$u-v$.
As a~result,
applying Shirshov's algorithm to a~set of semigroup relations~%
$S$
gives rise to a~complete set of semigroup relations
$S^c$.
The
$S^c$-irreducible words in~%
$X$
constitute the set of normal forms of the elements of~%
$P$.

Before we go any further,
let us give some well-known examples of algebra,
group,
and semigroup presentations 
by generators and defining relations
together with linear bases,
normal forms,
and GS bases
for them
(if~known).
Consider a~field~%
$k$
and a~commutative ring or commutative~%
$k$-algebra~%
$K$.
\begin{itemize}
\item
The Grassman algebra over
$K$
is
$$
K\langle X|x_i^2=0, x_ix_j+x_jx_i=0,\ i>j\rangle.
$$
The set of defining relations is a~GS basis
with respect to the deg-lex ordering.
A~%
$K$-basis is
$$
\{x_{i_1}\cdots x_{i_n}|x_{i_j}\in X,\ j=1,\dots,n,\  i_1<\dots< i_n,\ n\geq0\}.
$$

\item
The Clifford algebra over
$K$
is
$$
K\langle X| x_ix_j+x_jx_i= a_{ij},\ 1\leq i,j\leq n\rangle,
$$
where
$(a_{ij})$
is an
$n\times n$
symmetric matrix over~%
$K$.
The set of defining relations is a~GS basis
with respect to the deg-lex ordering.
A~%
$K$-basis is
$$
\{x_{i_1}\cdots x_{i_n}|x_{i_j}\in X,\ j=1,\dots,n,\  n\geq0,\ i_1<\dots< i_n\}.
$$

\item
The universal enveloping algebra of a~Lie algebra~%
$L$
is
$$
U_K(L)=K\langle X| x_ix_j-x_jx_i=\sum \alpha^k_{ij}x_k,\ i>j\rangle.
$$
If $L$ is a~free~%
$K$-module
with a~well-ordered~%
$K$-basis
$$
X=\{x_i|i\in I\},\ [x_ix_j]=\sum \alpha^k_{ij}x_k,\ i>j,\ i,j\in I,
$$
then the set of defining relations is a~GS basis of
$U_K(L)$.
The PBW theorem follows:
$U_K(L)$
is a~free~%
$K$-module
with a~%
$K$-basis,
as for polynomials,
$$
\{ x_{i_1}\cdots x_{i_n} \ |\ i_1\leq \cdots\leq i_n,\ i_t\in I,\ t=1,\dots,
n,\ n\geq 0\}.
$$

\item
A.~Kandri-Rody and V.~Weispfenning \cite{Kandri-Rody}
invented an~important class of
(noncommutative polynomial)
`algebras of solvable type',
which includes  universal enveloping algebras.
An~algebra of solvable type is
$$
R=k\langle X| s_{ij}=x_ix_j-x_jx_i-p_{ij},\ i>j,\ p_{ij}<x_ix_j\rangle,
$$
and the compositions
$(s_{ij},s_{jk})_w=0$
modulo
$(S, w)$,
where
$w=x_ix_jx_k$
with
$i>j>k$.
Here
$p_{ij}$
is a~noncommutative polynomial with all terms less than
$x_ix_j$.
They created a~theory of GS bases for every algebra of this class;
thus,
they found a~linear basis of every quotient of~%
$R$.

\item
A~general presentation
$U_k(L)=k\langle X|S^{(-)}\rangle$
of a~universal enveloping algebra over a~field~%
$k$,
where
$L=Lie(X|S)$
with
$S\subset Lie(X)\subset k\langle X\rangle$
and
$S^{(-)}$
is~%
$S$
as a~set of associative polynomials.
The PBW theorem in a~form of Shirshov's theorem.
The following conditions are equivalent:
\begin{enumerate}[(i)]
\item
the set~%
$S$
is a~Lie GS basis;
\item
the set
$S^{(-)}$
is a~GS basis for
$k\langle X\rangle$;
\item
a~linear
basis for
$U_k(L)$
consists of words
$u_1u_2\cdots u_n$,
where
$u_i$
are
$S$-irreducible LS words with
$u_1\preceq u_2\preceq\dots\preceq u_n$
(in the lex-ordering)%
\footnote{Given a~Lie polynomial
$s\in Lie(X)\subset k\langle X\rangle$,
the maximal associative word~%
$\bar s$
of~%
$s$
is an~LS word and the regular LS monomial
$[\bar s]$
is the maximal LS monomial of~%
$s$
(both in deg-lex ordering).},
see \cite{BoMa97,BoMa99};
\item
a~linear basis for~%
$L$
consists of the~%
$S$-irreducible
LS Lie monomials
$[u]$
in~%
$X$;
\item
a~linear basis for
$U_k(L)$
consists of the polynomials
$u=[u_1]\cdots [u_n]$,
where
$u_1\preceq\dots \preceq u_n$
in the lex ordering,
$n\geq0$,
and each
$[u_i]$
is an~%
$S$-irreducible
non-associative LS word in~%
$X$.
\end{enumerate}

\item
Free Lie algebras
$Lie_K(X)$
over~%
$K$.
M.~Hall,
A.~I.~Shirshov,
and R.~Lyndon provided different linear~%
$K$-bases
for a~free Lie algebra
(the Hall--Shirshov series of bases,
in particular,
the Hall basis,
the Lyndon--Shirshov basis,
the basis compatible
with the free solvable (polynilpotent) Lie algebra) \cite{ro99},
see also \cite{bokut63}.
Two anticommutative GS bases of
$Lie_K(X)$
were found in \cite{BoChLi09,BoChLi10},
which yields the Hall and Lyndon--Shirshov linear bases respectively.

\item
The Lie~%
$k$-algebras
presented by Chevalley generators and defining relations of types
$A_n$,
$B_n$,
$C_n$,
$D_n$,
$G_2$,
$F_4$,
$E_6$,
$E_7$,
and
$E_8$.
Serre's theorem provides linear bases and multiplication tables
for these algebras
(they are finite dimensional simple Lie algebras over~%
$k$).
Lie GS bases for these algebras are found in
\cite{BokutKlein96,BokutKlein98,BokutKlein99}.

\item
The Coxeter group
$$
W=\sgp\langle S| s_i^2=1,\ m_{ij}(s_i,s_j)=m_{ji}(s_j,s_i)\rangle
$$
for a~given Coxeter matrix
$M=(m_{ij})$.
J.~Tits \cite{Tits}
(see also \cite{Bjoner})
algorithmically solved the word problem for Coxeter groups.
Finite Coxeter groups are presented by `finite' Coxeter matrices
$A_n$,
$B_n$,
$D_n$,
$G_2$,
$F_4$,
$E_6$,
$E_7$,
$E_8$,
$H_3$,
and
$H_4$.
Coxeter's theorem provides normal forms and Cayley tables
(these are finite groups generated by reflections).
GS bases for finite Coxeter groups are found in \cite{bs}.

\item
The Iwahory--Hecke (Hecke) algebras~%
$H$
over~%
$K$
differ from the group algebras
$K(W)$
of Coxeter groups
in that
instead of
$s_i^2=1$
there are relations
$(s_i-q_i^{1/2})(s_i+q_i^{1/2})=0$
or
$(s_i-q_i)(s_i+1)=0$,
where
$q_i$
are units of~%
$K$.
Two
$K$-bases for
$H$
are known;
one is natural,
and the other is the Kazhdan--Lusztig canonical basis \cite{Lusztig2003}.
The GS bases for the Iwahory--Hecke algebras
are known for the finite Coxeter matrices.
A~deep connection of the Iwahory--Hecke algebras of type
$A_n$
and braid groups
(as well as link invariants)
was found by V.~F.~R.~Jones \cite{Jones}.

\item
Affine Kac--Moody algebras \cite{Kac}.
The Kac--Gabber theorem provides linear bases for
these algebras under the symmetrizability condition on the Cartan matrix.
Using this result,
E.~N.~Poroshenko found the GS bases of these algebras
\cite{poroshenkoB,poroshenkoC-D,poroshenkoA}.

\item
Borcherds--Kac--Moody algebras
\cite{Borcherds86,Borcherds88,Borcherds90,Kac}.
GS bases are not known.

\item
Quantum enveloping algebras
(V.~G.~Drinfeld, M.~Jimbo).
Lusztig's theorem \cite{Lusztig}
provides  linear canonical bases of these algebras.
Different approaches were developed by
C.~M.~Ringel \cite{Ringel90,Ringel96},
J.~A.~Green \cite{GreenJA},
and V.~K.~Kharchenko
\cite{Kharchenko99,Kharchenko02,Kharchenko05,Kharchenko08,Kharchenko10}.
GS bases of quantum enveloping algebras
are unknown except for the case
$A_n$,
see  \cite{BoMa,ChShaoShum,Rosso89,Ya}.

\item
Koszul algebras.
The quadratic algebras with a~basis of standard monomials,
called PBW-algebras,
are always Koszul
(S.~Priddy \cite{Priddy}),
but not conversely.
In different terminology,
PBW-algebras are algebras with quadratic GS bases.
See~\cite{Polishchuk}.

\item
Elliptic algebras (B.~Feigin, A.~Odesskii)
These are associative algebras presented by~%
$n$
generators and
$n(n-1)/2$
homogeneous quadratic relations
for which the dimensions of the graded components
are the same as for the polynomial algebra in~%
$n$
variables.
The first example of this type was Sklyanin algebra (1982)
generated by
$x_1$,
$x_2$,
and
$x_3$
with the defining relations
$[x_3,x_2]=x_1^2$,
$[x_2,x_1]=x_3^2$,
and
$[x_1,x_3]=x_2^2$.
See~\cite{Odesskii}.
GS bases are not known.

\item
Leavitt path algebras.
GS bases for these algebras are found in A.~Alahmedi et al \cite{AAJEfim}
and applied by the same authors to determine the structure
of the Leavitt path algebras of polynomial growth in \cite{AAJEfim2013}.

\item
Artin braid group
$Br_n$.
The Markov--Artin theorem provides
the normal form and semi-direct structure of the group
in the Burau generators.
Other normal forms of
$Br_n$
were obtained by Garside, Birman--Ko--Lee, and Adjan--Thurston.
GS bases for
$Br_n$
in the Artin--Burau,
Artin--Garside,
Birman--Ko--Lee,
and Adjan--Thurston generators were found in
\cite{Bo08,Bo09,Bo-Ch-Shum,ChZhong-Braid}
respectively.

\item
Artin--Tits groups.
They differ from Coxeter groups in the absence of the relations
$s_i^2=1$.
Normal forms are known in the spherical case,
see E.~Brieskorn, K.~Saito \cite{BrieskornSaito}.
GS bases are not known
except for braid groups
(the Artin--Tits groups of type~%
$A_n$).

\item
The groups of Novikov--Boon type
(Novikov \cite{Novikov},
Boon \cite{Boon},
Collins \cite{Collins},
Kalorkoti \cite{Kalorkoti82,Kalorkoti06,Kalorkoti09,Kalorkoti11})
with unsolvable word or conjugacy problem.
They are groups with standard bases
(standard normal forms or standard GS bases),
see \cite{bokut66,bokut67,Bok1968,ChChLuo}.

\item
Adjan's \cite{Adjan} and Rabin's \cite{Rabin} constructions of groups
with unsolvable isomorphism problem and Markov properties.
A~GS basis is known for Adjan's construction \cite{BokutChainikov08}.

\item
Markov's \cite{Markov} and Post's \cite{Post47} semigroups
with unsolvable word problem.
The GS basis of Post's semigroup is found in  \cite{Xia-Ming}.

\item
Markov's construction of semigroups
with unsolvable isomorphism problem and Markov properties.
The GS basis for the construction is not known.

\item
Plactic monoids.
A~theorem due to Richardson, Schensted, and Knuth
provides a~normal form of the elements of these monoids
(see M.~Lothaire \cite{Lothaire-etc}).
New approaches to plactic monoids
via GS bases in the alphabets of row and column generators
are found in \cite{ChLiJing-plactic}.

\item
The groups of quotients of the multiplicative semigroups
of power series rings with topological quadratic relations
of the type
$k\langle\langle x,y,z,t|xy=zt\rangle\rangle$
embeddable
(without the zero element)
into groups but in general not embeddable into division algebras
(settling a~problem of Malcev).
The relative standard normal forms of these groups
found in \cite{bokut69i-iii,bokut69iv}
are the reduced words for
what was later called a~\textit{relative} GS basis \cite{BokutShum}.
\end{itemize}

To date,
the method of GS bases has been adapted,
in particular,
to the following classes of linear universal algebras,
as well as for operads,
categories,
and semirings.
Unless stated otherwise,
we consider all linear algebras over a~field~%
$k$.
Following the terminology of Higgins and Kurosh,
we mean by a~((differential) associative)~%
$\Omega$-algebra
a~linear space
((differential) associative algebra)
with a~set of multi-linear operations~%
$\Omega$:
\begin{itemize}
\item
Associative algebras,
Shirshov \cite{Sh62b},
Bokut \cite{Bo76},
Bergman \cite{Be78};
\item
Associative algebras over a~commutative algebra,
A.~A.~Mikhalev, Zolotykh \cite{MZ};
\item
Associative
$\Gamma$-algebras,
where~%
$\Gamma$
is a~group,
Bokut, Shum \cite{BokutShum};
\item
Lie algebras,
Shirshov \cite{Sh62b};
\item
Lie algebras over a~commutative algebra,
Bokut, Chen, and Chen
\cite{BoChChen-Liecomm};
\item
Lie p-algebras over
$k$
with
$\operatorname{char} k=p$,
A.~A.~Mikhalev
\cite{Mikhalev92};
\item
Lie superalgebras,
A.~A.~Mikhalev \cite{Mikhalev89,Mikhalev96};
\item
Metabelian Lie algebras,
Chen, Chen \cite{CC12};
\item
Quiver (path) algebras,
Farkas, Feustel, and Green \cite{FFG};
\item
Tensor products of associative algebras,
Bokut, Chen, and Chen \cite{BCC08};
\item
Associative differential algebras,
Chen, Chen, and Li \cite{ChChLi-GSB-diff};
\item
Associative
$(n-)$conformal
algebras over
$k$
with
$\operatorname{char} k=0$,
Bokut, Fong, and Ke \cite{BFK04},
Bokut, Chen, and Zhang \cite{BoChZhang-n-conf};
\item
Dialgebras, Bokut, Chen, and Liu \cite{BCL08};
\item
Pre-Lie (Vinberg--Koszul--Gerstenhaber, right (left) symmetric) algebras,
Bokut, Chen, and Li \cite{BoChLi-GSB-rightsym},
\item
Associative Rota--Baxter algebras over
$k$
with
$\mathop{\mathrm{char}}k=0$,
Bokut, Chen, and Deng \cite{BCD08};
\item
$L$-algebras,
Bokut, Chen, and Huang \cite{BCH13};
\item
Associative~%
$\Omega$-algebras,
Bokut, Chen, and Qiu \cite{BoChQiu-CD-Omega};
\item
Associative differential~%
$\Omega$-algebras,
Qiu and Chen \cite{ChQiu-Cd-diff};
\item
$\Omega$-algebras,
Bokut, Chen, and Huang \cite{BCH13};
\item
Differential Rota--Baxter commutative associative algebras,
Guo, Sit, and Zhang \cite{Guoli2012};
\item
Semirings,
Bokut, Chen, and Mo \cite{BCM13};
\item
Modules over an~associative algebra,
Green \cite{Green},
Kang, Lee \cite{kl1,KL},
Chibrikov \cite{Chibrikov04};
\item
Small categories,
Bokut, Chen, and Li \cite{BoChLi-CD-category};
\item
Non-associative algebras,
Shirshov \cite{Sh62a};
\item
Non-associative algebras over a~commutative algebra,
Chen,  Li, and Zeng \cite{CLZ};
\item
Commutative non-associative algebras,
Shirshov \cite{Sh62a};
\item
Anti-commutative non-associative algebras,
Shirshov \cite{Sh62a};
\item
Symmetric operads,
Dotsenko and Khoroshkin \cite{DK10}.
\end{itemize}

\bigskip

At the heart of the GS method for a~class of linear algebras
lies a~CD-lemma for a~free object of the class.
For the cases above,
the free objects are
the free associative algebra
$k\langle X\rangle$,
the doubly free associative
$k[Y]$-algebra
$k[Y]\langle X\rangle$,
the free Lie algebra
$Lie(X)$,
and the doubly free Lie
$k[Y]$-algebra
$Lie_{k[Y]}(X)$.
For the tensor product of two associative algebras
we need to use the tensor product of two free algebras,
$k\langle X\rangle\otimes k\langle Y\rangle$.
We can view every semiring
as a~double semigroup with two associative products
$\cdot$~%
and~%
$\circ$.
So,
the CD-lemma for semirings is
the CD-lemma for the semiring algebra
of the free semiring
$Rig(X)$.
The CD-lemma for modules is
the CD-lemma for the doubly free module
$\Mod_{k\langle Y\rangle}(X)$,
a~free module over a~free associative algebra.
The CD-lemma for small categories is
the CD-lemma for the `free partial~%
$k$-algebra'
$kC\langle X\rangle$
generated by an~oriented graph~%
$X$
(a sequence
$z_1z_2\cdots z_n$,
where
$z_i\in X$,
is a~partial word in~%
$X$
iff it is a~path;
a~partial polynomial is
a~linear combination of partial words
with the same source and target).

All CD-lemmas have essentially the same statement.
Consider a~class~%
$\bold {V}$
of linear universal algebras,
a~free algebra
$\bold {V}(X)$
in
$\bold {V}$,
and a~well-ordered $k$-basis of terms
$N(X)$
of
$\bold {V}(X)$.
A~subset
$S\subset \bold {V}(X)$
is called a~GS basis
if every \textit{composition} of the elements of~%
$S$
is \textit{trivial}
(vanishes upon the elimination of the leading terms~%
$\bar {s}$
for
$s\in S$).
Then the following conditions are equivalent:
\begin{enumerate}[(i)]
\item
$S$
is a~GS basis.
\item
If
$f\in Id(S)$
then the leading term~%
$\bar{f}$
contains the subterm~%
$\bar{s}$
for some
$s\in S$.
\item
The set of~%
$S$-irreducible terms
is a~linear basis for the~%
$\bold {V}$-algebra
$\bold {V}\langle X|S\rangle$
generated by~%
$X$
with defining relations~%
$S$.
\end{enumerate}
In some cases
($(n-)$
conformal algebras,
dialgebras),
conditions (i)~and~(ii) are not equivalent.
To be more precise,
in those cases we have
$(i)\Rightarrow (ii)\Leftrightarrow (iii)$.

\bigskip

Typical compositions are compositions of intersection and inclusion.
Shirshov \cite{Sh62a,Sh62b} avoided inclusion composition.
He suggested instead that
a~GS basis must be minimal
(the leading words do not contain each other as subwords).
In some cases,
new compositions must be defined,
for example,
the composition of left (right) multiplication.
Also,
sometimes we need to combine all these compositions.
We present here a~new approach to
the definition of a~composition,
based on the concept of the least common multiple
$\lcm(u,v)$
of two terms
$u$~%
and~%
$v$.

In some cases
(Lie algebras,
($n$-)
conformal algebras)
the `leading' term~%
$\bar{f}$
of a~polynomial
$f\in \bold {V}(X)$
lies outside
$\bold {V}(X)$.
For Lie algebras,
we have
$\bar{f}\in k\langle X\rangle$,
for ($n$-) conformal algebras~%
$\bar{f}$
belongs to an~%
`$\Omega$-semigroup'.

Almost all CD-lemmas require the new notion of a~`normal~%
$S$-term'.
A~term
$(asb)$
in
$\{X, \Omega\}$,
where
$s\in S$,
with only one occurrence of~%
$s$
is called a~normal~%
$S$-term
whenever
$\overline{(asb)} = (a(\bar s)b)$.
Given
$S\subset k\langle X\rangle$,
every
$S$-word
(that is,
an~%
$S$-term)
is a~normal~%
$S$-word.
Given
$S\subset Lie(X)$,
every Lie~%
$S$-monomial
(Lie~%
$S$-term)
is a~linear combination of normal Lie~%
$S$-terms
(Shirshov \cite{Sh62b}).

One of the two key lemmas asserts that
if~%
$S$
is complete under compositions of multiplication
then every element of the ideal generated by~%
$S$
is a~linear combination of normal~%
$S$-terms.
Another key lemma says that
if~%
$S$
is a~GS basis
and the leading words of two normal~%
$S$-terms are the same
then these terms are the same modulo lower normal~%
$S$-terms.
As we mentioned above,
Shirshov proved these results \cite{Sh62b} for
$Lie(X)$
(there are no compositions of multiplication
for Lie and associative algebras).

This survey continues
our surveys with P.~S.~Kolesnikov, Y.~Fong, W.-F.~Ke,  and K.~P.~Shum
\cite{BC07,survey08,BCS,BFKK00,BK03,BK05}.

The paper is organized as follows.
Section~2 is for associative algebras,
Section~3 is for semigroups and groups,
Section~4 is for Lie algebras,
and the short Section~5 is for~%
$\Omega$-algebras and operads%
\footnote{The first definitions of the symmetric operad
were given by A.~G.~Kurosh's student V.~A.~Artamonov
under the name `clone of multilinear operations' in 1969,
see A.~G.~Kurosh \cite{Kurosh69}
and V.~A.~Artamonov \cite{Artamonov},
cf.\ J.~Lambek (1969) \cite{Lambek}
and P.~May (1972) \cite{May}.}.

To conclude this introduction,
we give some information about the work of Shirshov;
for more on this,
see the book \cite{Shirshov-Selected}.
A.~I.~Shirshov (1921-1981)
was a~famous Russian mathematician.
His name is associated with notions and results
on the Gr\"{o}bner--Shirshov bases,
the Composition-Diamond lemma,
the Shirshov--Witt theorem,
the Lazard-Shirshov elimination,
the Shirshov height theorem,
Lyndon--Shirshov words,
Lyndon--Shirshov basis
(in a~free Lie algebra),
the Hall--Shirshov series of bases,
the Cohn--Shirshov theorem for Jordan algebras,
Shirshov's theorem on the Kurosh problem,
and the Shirshov factorization theorem.
Shirshov's ideas were used by
his students Efim Zelmanov to solve the restricted Burnside problem
and Aleksander Kemer to solve the Specht problem.

We thank P.~S.~Kolesnikov, Yongshan Chen and Yu Li
for valuable comments and help in writing some parts of the survey.

\subsection{Digression on the history of Lyndon--Shirshov bases and
Lyndon--Shirshov words}

Lyndon \cite{Lyndon}, 1954,
defined standard words,
which are the same as Shirshov's regular words \cite{Sh58}, 1958.
Unfortunately,
the papers (Lyndon \cite{Lyndon})
and (Chen-Fox-Lyndon \cite{CFL}, 1958)
were practically unknown before 1983.
As a~result,
at that
time almost all  authors
(except four who used the names Shirshov and Chen--Fox--Lyndon,
see below)
refer to the basis and words
as Shirshov regular basis and words,
cf.\ for instance
\cite{Bakhturin11,Bakhturin92,Cohnbook65,Razmyslov,Ufnarovski95,Zelmanov92}.
To the best of our knowledge,
none of the authors mentioned Lyndon's paper \cite{Lyndon}
as a~source of `Lyndon words' before 1983(!).

In the following papers
the authors mentioned
both (Chen-Fox-Lyndon \cite{CFL})
and (Shirshov \cite{Sh58})
as a~source of `Lyndon--Shirshov basis' and
`Lyndon--Shirshov words':
\begin{itemize}
\item
M.~P.~Sch\"utzenberger, S.~Sherman \cite{Sch63}, 1963;
\item
M.~P.~Sch\"utzenberger \cite{Sch65}, 1965;
\item
G.~Viennot \cite{Viennot78}, 1978;
\item
J.~Michel \cite{Michel76}, 1976.
\end{itemize}
The authors of \cite{Sch63} thank P.~M.~Cohn
for pointing out Shirshov's paper \cite{Sh58}.
They also formulate Shirshov's(!) factorization theorem \cite{Sh58}.
They mention \cite{CFL,Sh58}
as a~source of `LS words'.
M.~P.~Sch\"utzenberger also mentions \cite{Sch65}
Shirshov's factorization theorem,
but in this case he attributes it to
both Chen-Fox-Lyndon \cite{CFL}
and Shirshov \cite{Sh58}.
Actually,
he cites \cite{CFL} by mistake,
as that result is absent from the paper,
see Berstel--Perrin \cite{Berstel-J-P}%
\footnote{From \cite{Berstel-J-P}:
``A famous theorem concerning Lyndon words asserts that
any word
$w$
can be factorized in a~unique way as a
non-increasing product of Lyndon words,
i.e. written
$w=x_1x_2\dots x_n$
with
$x_1\geq x_2\geq \dots \geq x_n$.
This theorem has imprecise origin.
It is usually credited to Chen--Fox--Lyndon,
following the paper of Sch\"utzenberger \cite{Sch65}
in which
it appears as an~example of factorization of free monoids.
Actually,
as pointed out to one of us by D.~Knuth in 2004,
the reference
\cite{CFL} does not contain explicitly this statement.'' }.

Starting with the book of M.~Lothaire,
\textsl{Combinatorics on words} (\cite{Lothaire-etc}, 1983),
some authors called the words and basis
`Lyndon words' and `Lyndon basis';
for instance,
see C.~Reutenauer, Free Lie algebras  (\cite{Reutenauer}, 1993).

\section{Gr\"obner--Shirshov bases for associative algebras}

In this section
we give a~proof of Shirshov's CD-lemma for associative algebras
and Buchberger's theorem for commutative algebras.
Also,
we give the Eisenbud--Peeva--Sturmfels lifting theorem,
the CD-lemmas for modules
(following S.-J.~Kang and K.-H.~Lee \cite{KL}
and E.~S.~Chibrikov \cite{Chibrikov04}),
the PBW theorem and the PBW theorem in Shirshov's form,
the CD-lemma for categories,
the CD-lemma for associative algebras over commutative algebras
and the Rosso--Yamane theorem for
$U_q(A_n)$.

\subsection{Composition-Diamond lemma for associative algebras}

Let
$k$
be a~field,
$k\langle X\rangle$
be the free associative algebra over~%
$k$
generated by~%
$X$
and
$X^{*}$
be the free monoid generated by~%
$X$,
where the empty word is the identity,
denoted by~1.
Denote the length (degree) of a~word
$w\in X^*$
by
$|w|$
or
$deg(w)$.
Suppose that
$X^*$
is a~well-ordered set.
Take
$f\in k\langle X\rangle$
with the leading word~%
$\bar{f}$
and
$f=\alpha \bar{f}-r_f$,
where
$0\neq\alpha\in k$
and
$\overline{r_f}<\bar f$.
We call~%
$f$
\textit{monic} if
$\alpha=1$.

A~well-ordering~%
$>$
on
$X^*$
is called a~\textit{monomial ordering}
whenever it is compatible with the multiplication of words,
that is,
for all
$u, v\in X^*$
we have
$$
u > v \Rightarrow w_{1}uw_{2} > w_{1}vw_{2},
\text{ for all }
w_{1}, \ w_{2}\in  X^*.
$$
A~standard example of monomial ordering on
$X^*$
is the deg-lex ordering,
in which two words are compared first by the degree
and then lexicographically,
where~%
$X$
is a~well-ordered set.

Fix a~monomial ordering
$<$
on
$X^*$
and take two monic polynomials
$f$~%
and~%
$g$
in
$k\langle X\rangle$.
There are two kinds of compositions:
\begin{enumerate}[(i)]
\item
If~%
$w$
is a~word such that
$w=\bar{f}b=a\bar{g}$
for some
$a,b\in X^*$
with
$|\bar{f}|+|\bar{g}|>|w|$
then the polynomial
$(f,g)_w=fb-ag$
is called the \textit{intersection composition} of
$f$~%
and~%
$g$
with respect to~%
$w$.
\item
If
$w=\bar{f}=a\bar{g}b$
for some
$a,b\in X^*$
then the polynomial
$(f,g)_w=f - agb$
is called the \textit{inclusion composition} of
$f$~%
and~%
$g$
with respect to~%
$w$.
\end{enumerate}
Then
$\overline{(f,g)_w}<w$
and
$(f,g)_w$
lies in the ideal
$Id\{f,g\}$
of
$k\langle X\rangle$
generated by
$f$~%
and~%
$g$.

In the composition
$(f,g)_w$,
we call~%
$w$
an~ambiguity
(or the least common multiple
$\lcm(\bar{f},\bar{g})$,
see below).

Consider
$S\subset k\langle X\rangle$
such that
very
$s\in S$
is monic.
Take
$h\in k\langle X\rangle$
and
$w\in X^*$.
Then~%
$h$
is called \textit{trivial} modulo
$(S,w)$,
denoted by
$$
h\equiv0 \mod (S,w),
$$
if
$h=\sum\alpha_i a_i s_i b_i$,
where
$\alpha_i\in k$,
$a_i,b_i\in X^{*}$,
and
$s_i\in S$
with
$a_i \overline{s_i} b_i<w$.

The elements
$asb$,
$a,b\in X^{*}$,
and
$s\in S$
are called~%
$S$-words.

A~monic set
$S\subset k\langle X\rangle$
is called a~GS basis in
$k\langle X\rangle$
with respect to the monomial ordering~%
$<$
if every composition of polynomials in~%
$S$
is trivial modulo~%
$S$
and the corresponding~%
$w$.

A~set~%
$S$
is called a~minimal GS basis in
$k\langle X\rangle$
if~%
$S$
is a~GS basis in
$k\langle X\rangle$
avoiding inclusion compositions;
that is,
given
$f,g\in S$
with
$f\neq g$,
we have
$\overline{f}\neq a\overline{g}b$
for all
$a,b\in X^*$.

Put
$$
Irr(S) = \{ u \in X^* |  u \neq a\bar{s}b ,s\in S,a ,b \in X^*\}.
$$
The elements of
$Irr(S)$
are called~%
$S$-irreducible
or
$S$-reduced.

A~GS basis~%
$S$
in
$k\langle X\rangle$
is reduced
provided that
$\supp(s)\subseteq  Irr(S\setminus\{s\})$
for every
$s\in S$,
where
$\supp(s)=\{u_1,u_2,\dots,u_n\}$
whenever
$s=\sum_{i=1}^n\alpha_iu_i$
with
$0\neq\alpha_i\in k$
and
$u_i\in X^*$.
In other words,
each
$u_i$
is an
$S\setminus\{s\}$-irreducible word.

The following lemma is key for proving
the CD-lemma for associative algebras.

\begin{lemma}\label{L4.1}
If~%
$S$
is a~GS basis in
$k \langle X\rangle$
and
$w=a_1\overline{s_1}b_1=a_2\overline{s_2}b_2$,
where
$a_1,b_1,a_2,b_2\in X^*$
and
$s_1,s_2\in S$,
then
$a_1s_1b_1\equiv a_2s_2b_2\ \mod (S,w)$.
\end{lemma}

\noindent
{\bf Proof.}
There are three cases to consider.

Case 1.
Assume that
the subwords
$\bar s_1$
and
$\bar s_2$
of~%
$w$
are disjoint,
say,
$|a_2|\geq |a_1|+|\bar s_1|$.
Then,
$a_2=a_1\bar s_1 c$
and
$b_1=c\bar s_2 b_2$
for some
$c\in X^*$,
and so
$w_1=a_1\bar s_1 c \bar s_2 b_2$.
Now,
\begin{eqnarray*}
a_1 s_1 b_1-a_2 s_2 b_2&=&a_1 s_1 c \bar s_2
b_2-a_1\bar s_1 c  s_2 b_2\\
&=&a_1 s_1 c (\bar s_2 - s_2) b_2+a_1(s_1-\bar s_1) c  s_2 b_2.
\end{eqnarray*}
Since
$\overline{\overline{s_2}-s_2}<\bar s_2$
and
$\overline{s_1-\overline{s_1}}<\bar s_1$,
we conclude that
$$
a_1 s_1 b_1-a_2 s_2 b_2=\sum\limits_i
\alpha_iu_is_1v_i+\sum\limits_j \beta_ju_js_2v_j
$$
with
$\alpha_i,\beta_j\in k$
and
$S$-words
$u_is_1v_i$
and
$u_js_2v_j$
satisfying
$
u_i\bar s_1v_i,u_j\bar s_2v_j<w.
$

Case 2.
Assume that
the subword~%
$\bar s_1$
of~%
$w$
contains
$\bar s_2$
as a~subword.
Then
$\bar s_1=a\bar s_2b$
with
$a_2=a_1a$
and
$b_2=bb_1$,
that is,
$w=a_1a\bar s_2bb_1$
for some
$S$-word
$a s_2 b$.
We have
\begin{eqnarray*}
a_1 s_1 b_1-a_2 s_2 b_2 = a_1 s_1 b_1-a_1 a~s_2 b b_1=
a_1(s_1-as_2b)b_1=a_1(s_1,s_2)_{\overline{s_1}} b_1.
\end{eqnarray*}
The triviality of compositions implies that
$
a_1s_1b_1\equiv
a_2s_2b_2\ \mod (S,w).
$

Case 3.
Assume that
the subwords
$\bar{s}_1$
and
$\bar{s}_2$
of
$w$
have a~nonempty intersection.
We may assume that
$a_2=a_1a$
and
$b_1=bb_2$
with
$w=\bar{s}_1b=a\bar{s}_2$
and
$|w|<|\bar{s}_1|+|\bar{s}_2|$.
Then,
as in Case~2,
we have
$
a_1s_1b_1\equiv a_2s_2b_2\ \mod (S,w).
$
\hfill
$\blacksquare$

\begin{lemma}\label{L4.2}
Consider a~set
$S\subset k \langle X\rangle$
of monic polynomials.
For every
$f\in k \langle X\rangle$
we have
\begin{equation*}
f=\sum\limits_{u_i\leq \bar f }\alpha_iu_i+
\sum\limits_{a_j\overline{s_j}b_j\leq\bar f}\beta_ja_js_jb_j
\end{equation*}
where
$\alpha_i,\beta_j\in k$,
$u_i\in Irr(S)$,
and
$a_js_jb_j$
are
$S$-words.
So,
$Irr(S)$
is a~set of linear generators of the algebra
$f\in k \langle X|S\rangle$.
\end{lemma}

\noindent {\bf Proof.}
Induct on~%
$\bar f$.\hfill
$\blacksquare$

\begin{theorem}
\label{CD asso}(The CD-lemma for associative algebras)
Choose a~monomial ordering~%
$<$
on
$X^*$.
Consider a~monic set
$S\subset k \langle X\rangle$
and  the ideal
$Id(S)$
of
$k \langle X\rangle$
generated by~%
$S$.
The following statements are equivalent:
\begin{enumerate}[\rm (i)]
\item
$S$
is a~Gr\"{o}bner--Shirshov basis in
$k\langle X\rangle$.
\item
$f\in Id(S)\Rightarrow \bar{f}=a\bar{s}b$
for some
$s\in S$
and
$a,b\in  X^*$.
\item
$Irr(S) = \{ u \in X^* |  u \neq a\bar{s}b ,s\in S,a ,b \in X^*\}$
is a~linear basis of the algebra
$k\langle X | S \rangle$.
\end{enumerate}
\end{theorem}

\noindent {\bf Proof.}
(i)$\Rightarrow$(ii).
Assume that
$S$
is a~GS basis
and take
$0\neq f\in Id(S)$.
Then,
we have
$
f=\sum_{i=1}^n\alpha_ia_is_ib_i
$
where
$\alpha_i\in k$,
$a_i,b_i\in X^*$,
and
$s_i\in S$.
Suppose that
$w_i=a_i\overline{s_i}b_i$
satisfy
$$
w_1=w_2=\cdots=w_l>w_{l+1}\geq\cdots.
$$
Induct on
$w_1$~%
and~%
$l$
to show that
$\overline{f}=a\overline{s}b$
for some
$s\in S \ \mbox{and} \ a,b\in X^*$.
To be more precise,
induct on
$(w_1,l)$
with the lex ordering of the pairs.

If
$l=1$
then
$\overline{f}=\overline{a_1s_1b_1}=a_1\overline{s_1}b_1$
and hence the claim holds.
Assume that
$l\geq 2$.
Then
$w_1=a_1\overline{s_1}b_1= a_2\overline{s_2}b_2$.
Lemma \ref{L4.1} implies that
$
a_1s_1b_1\equiv a_2s_2b_2 \mod(S,w_1).
$
If
$\alpha_1+\alpha_2\neq 0$
or
$l>2$
then the claim follows by induction on~%
$l$.
For the case
$\alpha_1+\alpha_2= 0$
and
$l=2$,
induct on
$w_1$.
Thus,
(ii) holds.

(ii)$\Rightarrow$(iii).
By Lemma \ref{L4.2},
$Irr(S)$
generates
$k\langle X|S\rangle$
as a~linear space.
Suppose that
$\sum\limits_{i}\alpha_iu_i=0$
in
$k\langle X|S\rangle$,
where
$0\neq\alpha_i\in k$
and
$u_i\in {Irr(S)}$.
It means that
$\sum\limits_{i}\alpha_iu_i\in{Id(S)}$
in
${k\langle X\rangle}$.
Then
$\overline{\sum\limits_{i}\alpha_iu_i}=u_j\in{Irr(S)}$
for some~%
$j$,
which contradicts (ii).

(iii)$\Rightarrow$(i).
Given
$f,g\in{S}$,
Lemma \ref{L4.2} and (iii) yield
$
(f,g)_{w}\equiv0\ \mod(S,w).
$
Therefore,
$S$
is a~GS basis.
 \hfill
$\blacksquare$

\medskip

A~new exposition of the proof
of Theorem \ref{CD asso} (CD-lemma for associative algebras).

\medskip

Let us start with the concepts of
non-unique \textit{common multiple}
and \textit{least common multiple} of two words
$u,v\in X^*$.
A~common multiple
$\cm(u,v)$
means that
$\cm(u,v)=a_1ub_1 = a_2vb_2$
for some
$a_i,b_i\in X^*$.
Then
$\lcm(u,v)$
means that
some
$\cm(u,v)$
contains some
$\lcm(u,v)$
as a~subword:
$\cm(u,v)=c\cdot \lcm(u,v)\cdot d$
with
$c,d\in X^*$,
where
$u$~%
and~%
$v$
are the same subwords in both sides. 
To be precise,
\begin{multline*}
\lcm(u,v)\in \{ucv, c\in X^*
\text{ (a trivial $\lcm(u,v)$)};
\\
u=avb, \ a,b\in X^*
\text{ (an inclusion $\lcm(u,v)$)};
\\
ub=av, \ a,b\in X^*, \
|ub|<|u|+|v| \text{ (an intersection $\lcm(u,v)$)}\}.
\end{multline*}
Define the general composition
$(f,g)_{\lcm(\bar f, \bar g)}$
of monic polynomials
$f,g\in k\langle X\rangle$
as
$$
(f,g)_{\lcm(\bar f, \bar g)} = \lcm(\bar f, \bar g)|_{\bar f\mapsto f}
- \lcm(\bar f, \bar g)|_{\bar g\mapsto g}.
$$
The only difference with the previous definition of composition is that
we include the case of trivial
$\lcm(\bar f, \bar g)$.
However,
in this case the composition is trivial,
$$
(f,g)_{\bar f c \bar g}\equiv 0 \mod (\{f,g\},\bar f c \bar g).
$$
It is clear that if
$a_1\bar fb_1 = a_2\bar gb_2$
then,
up to the ordering of
$f$~%
and~%
$g$,
$$
a_1fb_1 -a_2gb_2 = c\cdot (f,g)_{\lcm(\bar f, \bar g)}\cdot d.
$$
This implies Lemma \ref{L4.1}.
The main claim
(i)$\Rightarrow$(ii)
of Theorem \ref{CD asso}
follows from Lemma \ref{L4.1}.

\medskip

\noindent{\bf Shirshov algorithm.}
If a~monic subset
$S \subset k \langle X\rangle$
is not a~GS basis
then we can add to~%
$S$
all nontrivial compositions,
making them monic.
Iterating this process,
we eventually obtain a~GS basis
$S^{c}$
that contains~%
$S$
and generates the same ideal,
$Id(S^c)=Id(S)$.
This
$S^{c}$
is called the GS completion of~%
$S$.
Using the reduction algorithm
(elimination of the leading words of polynomials),
we may obtain a~minimal GS basis
$S^{c}$
or a~reduced GS basis.

\medskip

The following theorem gives a~linear basis for the ideal
$Id(S)$
provided that
$S \subset k \langle X\rangle$
is a~GS basis.

\begin{theorem}\label{t4.4}
If
$S \subset k \langle X\rangle$
is a~Gr\"{o}bner--Shirshov basis
then,
given
$u\in X^*\setminus Irr(S)$,
by Lemma \ref{L4.2} there exists
$\widehat{u}\in kIrr(S)$
with
$\overline{\widehat{u}}< u$
(if
$\widehat{u}\neq0$)
such that
$u-\widehat{u}\in Id(S)$
and the set
$\{u-\widehat{u}|u\in X^*\setminus Irr(S)\}$
is a~linear basis for the ideal
$Id(S)$
of
$k \langle X\rangle$.
\end{theorem}

\noindent {\bf Proof.}
Take
$0\neq f\in Id(S)$.
Then by the CD-lemma for associative algebras,
$\bar{f}=a_1\overline{s_1}b_1=u_1$
for
some
$s_1\in S$
and
$a_1,b_1\in  X^*$,
which implies that
$\bar{f}=u_1\in X^*\setminus Irr(S)$.
Put
$f_1=f-\alpha_1(u_1-\widehat{u_1})$,
where
$\alpha_1$
is the coefficient of the leading term of~%
$f$
and
$\overline{\widehat{u_1}}< u_1$
or
$\widehat{u_1}=0$.
Then
$f_1\in
Id(S)$
and
$\overline{f_1}<\bar{f}$.
By induction on
$\bar f$,
the
set
$\{u-\widehat{u}|u\in X^*\setminus Irr(S)\}$
generates
$Id(S)$
as a~linear space.
It is clear that
$\{u-\widehat{u}|u\in X^*\setminus Irr(S)\}$
is a~linearly independent set.
\hfill
$\blacksquare$

\begin{theorem}
Choose a~monomial ordering
$>$
on
$X^*$.
For every ideal
$I$
of
$k\langle X\rangle$
there exists a~unique reduced Gr\"{o}bner--Shirshov basis
$S$
for
$I$.
\end{theorem}

\noindent {\bf Proof.}
Clearly,
a~Gr\"{o}bner--Shirshov basis
$S\subset k\langle X\rangle$
for the ideal
$I=Id(S)$
exists;
for example,
we may take
$S=I$.
By Theorem \ref{CD asso},
we may assume that
the leading terms of the elements of~%
$S$
are distinct.
Given
$g\in S$,
put
$$
\Delta_g=\{f\in S|f\neq g
\mbox{ and } \overline{f}=a\bar{g}b
\mbox{ for some } a,b\in X^*\}
$$
and
$S_1=S\setminus\cup_{g\in S}\Delta_g$.

For every
$f\in Id(S)$
we show that
there exists an
$s_1\in S_1$
such that
$\overline{f}=a\overline{s_1}b \mbox{ for some } a,b\in X^*$.

In fact,
Theorem \ref{CD asso} implies that
$\overline{f}=a'\bar{h}b'$
for some
$a',b'\in X^*$
and
$h\in S$.
Suppose that
$h\in S\setminus S_1$.
Then we have
$h\in \cup_{g\in S}\Delta_g$,
say,
$h\in \Delta_g$.
Therefore,
$h\neq g$
and
$\overline{h}=a\bar{g}b$
for some
$a,b\in X^*$.
We claim that
$\bar{h}>\bar{g}$.
Otherwise,
$\bar{h}<\bar{g}$.
It follows that
$\bar{h}=a\bar{g}b>a\bar{h}b$
 and so we have the infinite descending chain
$$
\bar{h}>a\bar{h}b >a^2\bar{h}b^2 >a^3\bar{h}b^3>\dots,
$$
which contradicts the assumption that
$>$
is a~well ordering.

Suppose that
$g\not\in S_1$.
Then,
by the argument above,
there exists
$g_1\in S$
such that
$g\in \Delta_{g_1}$
and
$\overline{g}>\overline{g_1}$.
Since
$>$
is a~well ordering,
there must exist
$s_1\in S_1$
such that
$\overline{f}=a_1\overline{s_1}b_1 \mbox{ for some } a_1,b_1\in
 X^*$.

Put
$f_1=f-\alpha_1a_1s_1b_1$,
where
$\alpha_1$
is the coefficient of the leading term of~%
$f$.
Then
$f_1\in Id(S)$
and
$\overline{f}>\overline{f_1}$.

By induction on
$\overline{f}$,
we know that
$f\in Id(S_1)$,
and hence
$I=Id(S_1)$.
Moreover,
Theorem \ref{CD asso} implies that
$S_1$
is clearly a~minimal GS basis for the ideal
$Id(S)$.

Assume that~%
$S$
is  a~minimal GS basis for~
$I$.

For every
$s\in S$
we have
$s=s'+s''$,
where
$\supp(s')\subseteq Irr(S\setminus\{s\})$
and
$s''\in Id(S\setminus\{s\})$.
Since~%
$S$
is a~minimal GS basis,
it follows that
$\overline{s}=\overline{s'}$
for every
$s\in S$.

We claim that
$S_2=\{s'|s\in S\}$
is a~reduced GS basis for~
$I$.
In fact,
it is clear that
$S_2\subseteq Id(S)=I$.
By Theorem \ref{CD asso},
for every
$f\in Id(S)$
we have
$\overline{f}=a_1\overline{s_1}b_1 =a_1\overline{s_1'}b_1$
for some
$a_1,b_1\in  X^*$.

Take two reduced GS bases
$S$~%
and~%
$R$
for the ideal~%
$I$.
By Theorem \ref{CD asso},
for every
$s\in S$,
$$
\overline{s}=a\overline{r}b, \quad
\overline{r}=c\overline{s_1}d
$$
for some
$a,b,c,d\in X^*$,
$r\in R$,
and
$s_1\in S$,
and hence
$\overline{s}=ac\overline{s_1}db$.
Since
$\bar s\in \supp(s)\subseteq Irr(S\setminus\{s\})$,
we have
$s=s_1$.
It follows that
$a=b=c=d=1$,
and so
$\overline{s}=\overline{r}$.

If
$s\neq r$
then
$0\neq s-r\in I=Id(S)=Id(R)$.
By Theorem \ref{CD asso},
$\overline{s-r}=a_1\overline{r_1}b_1=c_1\overline{s_2}d_1$
for some
$a_1,b_1,c_1,d_1\in X^*$
with
$\overline{r_1},\overline{s_2}<\overline{s}=\overline{r}$.
This means that
$s_2\in S\setminus\{s\}$
and
$r_1\in R\setminus\{r\}$.
Noting that
$\overline{s-r}\in \supp(s)\cup \supp(r)$,
we have either
$\overline{s-r}\in \supp(s)$
or
$\overline{s-r}\in  \supp(r)$.
If
$\overline{s-r}\in \supp(s)$
then
$\overline{s-r}\in Irr(S\setminus\{s\})$,
which contradicts
$\overline{s-r}=c_1\overline{s_2}d_1$;
if
$\overline{s-r}\in \supp(r)$
then
$\overline{s-r}\in Irr(R\setminus\{r\})$,
which contradicts
$\overline{s-r}=a_1\overline{r_1}b_1$.
This shows that
$s= r$,
and then
$S\subseteq R$.
Similarly,
$R\subseteq S$.
\hfill
$\blacksquare$

\medskip

\noindent{\bf Remark 1.}
In fact,
a~reduced GS basis is unique
(up to the ordering)
in all possible cases below.

\medskip

\noindent{\bf Remark 2.}
Both associative and Lie CD-lemmas are valid
when we replace the base field~%
$k$
by an~arbitrary commutative ring~%
$K$
with identity
because we assume that all GS bases consist of monic polynomials.
For example,
consider a~Lie algebra~%
$L$
over~%
$K$
which is a~free~
$K$-module
with a~well-ordered~%
$K$-basis
$\{a_i|i\in I\}$.
With the deg-lex ordering on
$\{a_i|i\in I\}^*$,
the universal enveloping associative algebra
$U_K(L)$
has a (monic) GS basis
$$
\{a_ia_j-a_ja_i=\sum\alpha_{ij}^ta_t| i>j,\ i,j\in I\},
$$
where
$\alpha_{ij}^t\in K$
and
$[a_i,a_j]=\sum\alpha_{ij}^ta_t$
in~%
$L$,
and the CD-lemma for associative algebras over~%
$K$
implies that
$L\subset U_K(L)$
and
$$
\{a_{i_1}\cdots a_{i_n}|
i_1\leq\dots\leq i_n,\
n\geq0,\
i_1,\dots,i_n\in I\}
$$
is a~%
$K$-basis for
$U_K(L)$.

In fact,
for the same reason,
all CD-lemmas in this survey are valid
if we replace the base field~%
$k$
by an~arbitrary commutative ring~%
$K$
with identity.
If this is the case then
claim (iii) in the CD-lemma should read:
$K(X|S)$
is a~free~%
$K$-module with a~%
$K$-basis
$Irr(S)$.
But in the general case,
Shirshov's algorithm fails:
if~%
$S$
is a~monic set then
$S'$,
the set obtained by adding to~%
$S$
all non-trivial compositions,
is not a~monic set in general,
and the algorithm may stop with no result.

\subsection{Gr\"{o}bner bases for commutative algebras
and their lifting to Gr\"{o}bner--Shirshov bases}

Consider the  free commutative associative algebra
$k[X]$.
Given a~well ordering~%
$<$
on
$X=\{x_i|i\in I\}$,
$$
[X]=\{x_{i_1}\dots
x_{i_t}|i_1\leq\dots\leq i_t,\ i_1,\dots, i_t\in I,\ t\geq0\}
$$
is a~linear basis for
$k[X]$.

Choose a~monomial ordering~%
$<$
on
$[X]$.
Take two monic polynomials
$f$~%
and~%
$g$
in
$k[X]$
such that
$w=\lcm(\bar{f},\bar{g})=\bar{f}a=\bar{g}b$
for some
$a,b\in [X]$
with
$|\bar{f}|+|\bar{g}|>|w|$
(so,
$\bar{f}$~%
and~%
$\bar{g}$
are not coprime in
$[X]$).
Then 
$(f,g)_w=fa-gb$
is called the~%
$s$-polynomial of
$f$~%
and~%
$g$.

A~monic subset
$S\subseteq k[X]$
is called a~Gr\"{o}bner basis with respect to the monomial ordering~%
$<$
whenever all
$s$-polynomials
of two arbitrary polynomials in
$S$
are trivial modulo~%
$S$.

An~argument similar to the proof of the CD-lemma for associative algebras
justifies the following theorem due to B.~Buchberger.

\begin{theorem}\label{CD asso comm}(Buchberger Theorem)\
Choose a~monomial ordering~%
$<$
on
$[X]$.
Consider a~monic set
$S\subset k [X]$
and the ideal
$Id(S)$
of
$k [X]$
generated by~%
$S$.
The following statements are equivalent:
\begin{enumerate}[\rm (i)]
\item
$S$
is a~Gr\"{o}bner basis in
$k[X]$.
\item
$f\in Id(S)\Rightarrow \bar{f}= \bar{s}a$
for some
$s\in S$
and
$a\in  [X]$.
\item
$Irr(S) = \{ u \in [X] |  u \neq \bar{s}a ,s\in S,a  \in [X]\}$
is a~linear basis for the algebra
$k[X|S]=k[X]/Id(S)$.
\end{enumerate}
\end{theorem}

\noindent {\bf Proof.}
Denote by
$\lcm(u,v)$
be the usual (unique) least common multiple of two commutative words
$u,v\in [X]$:
\begin{multline*}
\lcm(u,v)\in \{uv
\text{ (the trivial $\lcm(u,v)$)};
\\
au=bv,\ a,b\in [X],\ |au| < |u| + |v|
\text{ (the nontrivial $\lcm(u,v)$)}\}.
\end{multline*}
If
$\cm(u,v) =a_1u=a_2v$
 is a~common multiple of
$u$~%
and~%
$v$
then
$\cm(u,v)= b\cdot \lcm(u,v)$.

The~%
$s$-polynomial
of two monic polynomials
$f$~%
and~%
$g$
is
$$
(f,g)_{_{\lcm(\bar f, \bar g)}} =
\lcm(\bar f, \bar g)|_{\bar f\mapsto f} - \lcm(\bar f, \bar g)|_{\bar g\mapsto g}.
$$

An~analogue of Lemma \ref{L4.1} is valid for
$k[X]$
because if
$a_1\bar s_1=a_2\bar s_2$
for two monic polynomials
$s_1$
and
$s_2$
then
$$
a_1s_1-a_2s_2 = b\cdot (s_1,s_2)_{\lcm(\bar s_1, \bar s_2)}.
$$
Lemma \ref{L4.1} implies the main claim
(i)$\Rightarrow$(ii)
of Buchberger's theorem.
\hfill
$\blacksquare$

\medskip

\begin{theorem}
Given an~ideal~%
$I$
of
$k[X]$
and a~monomial ordering~%
$<$
on
$[X]$,
there exists a~unique reduced Gr\"{o}bner basis~%
$S$
for~%
$I$.
Moreover,
if~%
$X$
is finite
then so is~%
$S$.
\end{theorem}

\medskip

Eisenbud,
Peeva,
and Sturmfels constructed \cite{eisenbud}
a~GS basis in
$k\langle X\rangle$
by lifting a~commutative Gr\"{o}bner basis for
$k[X]$
and adding all commutators.
Write
$X=\{x_1,x_2,\dots,x_n\}$
and put
$$
S_1=\{h_{ij}=x_ix_j-x_jx_i| \ i>j\}\subset k\langle X \rangle.
$$
Consider the natural map
$\gamma\colon k\langle X \rangle\rightarrow k[X]$
carrying
$x_i$
to
$x_i$
and the lexicographic splitting of
$\gamma$,
which is defined as the~%
$k$-linear map
$$
\delta\colon k[X]\rightarrow k\langle X \rangle, \quad
x_{i_1}x_{i_2}\cdots x_{i_r}\mapsto x_{i_1}x_{i_2}\cdots x_{i_r}
\quad \text{if}\quad i_1\leq i_2\cdots \leq i_r.
$$

Given
$u\in [X]$,
we express it as
$u=x_1^{l_1}x_2^{l_2}\cdots x_n^{l_n}$,
where
$l_i\geq 0$,
using an~arbitrary monomial ordering on
$[X]$.

Following \cite{eisenbud},
define an~ordering on
$X^*$
using the ordering
$x_1<x_2<\cdots<x_n$
as follows:
given
$u,v\in X^*$,
put
$$
u>v\ \
\mbox{ if }\ \ \gamma(u)>\gamma(v)\
\mbox{ in } \ [X]\
\mbox{ or } \ (\gamma(u)=\gamma(v) \
\mbox{and} \ u>_{lex}v).
$$
It is easy to check that
this is a~monomial ordering on
$X^*$
and
$\overline{\delta(s)}=\delta(\bar{s})$
for every
$s\in k[X]$.
Moreover,
$v\geq\delta(u)$
for every
$v\in \gamma^{-1}(u)$.

Consider an~arbitrary ideal~%
$L$
of
$k[X]$
generated by monomials.
Given
$m=x_{i_1}x_{i_2}\cdots x_{i_r}\in L, \ i_1\leq i_2\cdots \leq i_r$,
denote by
$U_L(m)$
the set of all monomials
$u\in [x_{i_1+1},\cdots, x_{i_r-1}]$
such that
neither
$ux_{i_2}\cdots x_{i_r}$
nor
$ux_{i_1}\cdots x_{i_{r-1}}$
lie in~%
$L$.

\begin{theorem}\label{th9}(\cite{eisenbud})
Consider the orderings on
$[X]$
and
$X^*$
defined above.
If~%
$S$
is a~minimal Gr\"{o}bner basis in
$k[X]$
then
$S'=\{\delta(us)  |  s\in S, u\in U_L(\bar{s}) \} \cup S_1$
is a~minimal Gr\"{o}bner--Shirshov basis in
$k\langle X\rangle$,
where~%
$L$
is the monomial ideal of
$k[X]$
generated by~%
$\bar{S}$.
\end{theorem}

%

Jointly with Bokut, Chen, and Chen
\cite {BCC08},
we generalized this result
to lifting a~GS basis
$S\subset k[Y]\otimes k\langle X\rangle$,
see A.~A.~Mikhalev, A.~Zolotykh \cite {MZ},
to a~GS basis of
$
Id(S,[y_i,y_j]
\text{ for all }
(i,j)
$)
of
$k\langle Y\rangle \otimes k\langle X\rangle$.

\medskip

Recall that
for a~prime number~%
$p$
the Gauss ordering on the natural numbers is described as
$s\leq_p t$
whenever
%
%
$\tbinom{t}{s} \not \equiv0 \mod p$.
Let
$\mathord{\leq_0} = \mathord{\leq}$
be the usual
ordering on the natural numbers. A~monomial ideal
$L$
of
$k[X]$
is
called
$p$-Borel-fixed
whenever it satisfies the following condition:
for each monomial generator~%
$m$
of~%
$L$,
if
$m$
is divisible by
$x_j^t$
but no higher power of
$x_j$
then
$(x_i/x_j)^sm\in L$
for all
$i<j$
and
$s\leq_p t$.

Thus,
we have the following Eisenbud--Peeva--Sturmfels lifting theorem.

\begin{theorem}\label{th9.1}(\cite{eisenbud})
Given an~ideal~%
$I$
of
$k[X]$,
take
$L=Id(\bar{f}, \ f\in I)$
and
$J=\gamma^{-1}(I)\subset k\langle X\rangle$.
\begin{enumerate}[\rm(i)]
\item
If~%
$L$
is~%
$0$-Borel-fixed
then a~minimal Gr\"{o}bner--Shirshov basis of~%
$J$
is obtained by applying~%
$\delta$
to a~minimal Gr\"{o}bner basis of~%
$I$
and adding commutators.
\item
If~%
$L$
is~%
$p$-Borel-fixed
for some~%
$p$
then~%
$J$
has a~finite Gr\"{o}bner--Shirshov basis.
\end{enumerate}
\end{theorem}

\noindent {\bf Proof.}
Assume that~%
$L$
is~%
$p$-Borel-fixed
for some~%
$p$.
Take a~generator
$m=x_{i_1}x_{i_2}\cdots x_{i_r}$
of~%
$L$,
where
$x_{i_1}\leq x_{i_2}\leq \cdots \leq x_{i_r}$,
and suppose that
$x_{i_r}^t$
is the highest power of
$x_{i_r}$
dividing~%
$m$.
Since
$t\leq_p t$,
it follows that
$x_l^tm/x_{i_r}^t\in L$
for
$l<i_r$.
This implies that
$x_l^tm/x_{i_r}\in L$
for
$l<i_r$,
and hence,
every monomial in
$U_L(m)$
satisfies
$deg_{x_l}(u)<t$
for
$i_1<l<i_r$.
Thus,
$U_L(m)$
is a~finite set,
and the result follows from Theorem \ref{th9}.
In particular,
if
$p=0$
then
$U_L(m)={1}$.
\hfill
$\blacksquare$

\medskip

In characteristic
$p\geq0$
observe that
if the field~%
$k$
is infinite
then after a~generic change of variables~%
$L$
is~%
$p$-Borel-fixed.
Then Theorems \ref{th9} and \ref{th9.1} imply 

\begin{corollary}\label{co3.2}(\cite{eisenbud})
Consider an~infinite field
$k$
and an~ideal
$I\subset k[X]$.
After a~general linear change of variables,
the ideal
$\gamma^{-1}(I)$
in
$k\langle X\rangle$
has a~finite Gr\"{o}bner--Shirshov basis.
\end{corollary}

\subsection{Composition-Diamond lemma for  modules}

Consider
$S$,
$T\subset k\langle X\rangle$
and
$f$,
$g\in k\langle X\rangle$.
Kang and Lee define \cite{kl1}
the composition of
$f$~%
and~%
$g$
as follows.

\begin{definition}\em (\cite{kl1,kl2})\label{d1}
\begin{enumerate}[\textrm(a)]
\item
If there exist
$a, \ b\in X^*$
such that
$w=\overline{f}a=b\overline{g}$
with
$|w|<|\overline{f}|+|\overline{g}|$
then the intersection composition is defined as
$(f,g)_w=fa-bg$.
\item
If there exist
$a$,
$b\in X^*$
such that
$w=a\overline{f}b=\overline{g}$
then the inclusion composition is defined as
$(f,g)_w=afb-g$.
\item
The composition
$(f,g)_w$
is called \textit{right-justified}
whenever
$w=\overline{f}=a\overline{g}$
for some
$a\in X^*$.
\end{enumerate}
\end{definition}

If
$f-g=\sum \alpha_i a_is_i b_i+\sum \beta_j c_j t_j$,
where
$\alpha_i, \beta_j \in k$,
$a_i, b_i, c_j\in X^*$,
$s_i\in S$,
and
$t_j\in T$
with
$a_i \overline{s}_i b_i<w$
and
$c_j \overline{t}_j<w$
for all
$i$~%
and~%
$j$,
then we call
$f-g$
\textit{trivial} with respect to
$S$~%
and~%
$T$
and write
$f\equiv g \mod(S,T;w)$.

\begin{definition} (\cite{kl1,KL})
A~pair
$(S,T)$
of monic subsets of
$k\langle X\rangle$
is called a~\textit{GS pair}
if~%
$S$
is closed under composition,~%
$T$
is closed under right-justified composition with respect to~%
$S$,
and given
$f\in S$,
$g\in T$,
and
$w\in X^*$
such that if
$(f,g)_w$
is defined,
we have
$(f,g)_w\equiv 0 \mod (S,T;w)$.
In this case,
say that
$(S,T)$
is a~\textit{GS pair}
for the~%
$A$-module
$_AM =_Ak\langle X\rangle/(k\langle X\rangle T + Id(S))$,
where
$A=k\langle X|S\rangle$.
\end{definition}

\begin{theorem}\label{thm2.13}
(Kang and Lee \cite{kl1,KL},
the CD-lemma for cyclic modules)
Consider a~pair
$(S, T)$
of monic subsets of
$k\langle X\rangle$,
the associative algebra
$A = k\langle X|S\rangle$
defined by~%
$S$,
and the left cyclic module
$_AM =_Ak\langle X\rangle/(k\langle X\rangle T + Id(S))$
defined by
$(S, T)$.
Suppose that
$(S, T)$
is a~Gr\"{o}bner--Shirshov pair for the~%
$A$-module
$_AM$
and
$p\in k\langle X\rangle T + Id(S)$.
Then
$\bar{p} = a\bar{s}b$
or
$\bar{p} = c\bar{t}$,
where
$a, b, c \in X^*$,
$s \in S$,
and
$t \in T$.
\end{theorem}

Applications of Theorem \ref{thm2.13}
appeared in
\cite{KangLL02,KangLL04,kl2}.

Take two sets
$X$~%
and~%
$Y$
and consider the free left
$k\langle X\rangle$-module
$\Mod_{k\langle X\rangle} \langle Y\rangle$
with $k\langle X\rangle$-basis~%
$Y$.
Then
$\Mod_{k\langle X\rangle} \langle Y\rangle=\oplus_{y\in Y} k\langle X\rangle y$
is called a~\textit{double-free} module.
We now define  the GS basis in
$\Mod_{k\langle X\rangle} \langle Y\rangle$.
Choose a~monomial ordering~%
$<$
on
$X^*$,
and a~well-ordering~%
$<$
on~%
$Y$.
Put
$X^*Y=\{uy|u\in X^*, \ y\in Y\}$
and define an~ordering
$<$
on
$X^*Y$
as follows:
for any
$w_1=u_1y_1$,
$w_2=u_2y_2\in X^*Y$,
$$
w_1< w_2\Leftrightarrow
u_1<u_2 \quad \mbox{ or } \quad
u_1=u_2, \ y_1<y_2
$$

Given
$S\subset \Mod_{k\langle X\rangle} \langle Y\rangle$
with all
$s\in S$
monic,
define composition in~%
$S$
to be only inclusion composition,
which means that
$\bar{f}=a\bar{g}$
for some
$a\in X^*$,
where
$f,g\in S$.
If
$(f,g)_{\bar f}=f-ag=\sum \alpha_i a_i s_i$,
where
$\alpha_i\in k$,
$a_i \in X^*$,
$s_i \in S$,
and
$a_i \overline{s}_i< {\bar f}$,
then this composition is called trivial modulo
$(S,{\bar f})$.

\begin{theorem}\label{l3.5}
(Chibrikov \cite{Chibrikov04},
see also \cite{CCZ}, the CD-lemma for modules)
\label{l2}
Consider a~non-empty set
$S\subset mod_{k\langle X\rangle} \langle Y\rangle$
with all
$s\in S$
monic
and choose an~ordering~%
$<$
on
$X^*Y$
as before.
The following statements are equivalent:
\begin{enumerate}[\textrm(i)]
\item
$S$
is a~Gr\"{o}bner--Shirshov basis in
$\Mod_{k\langle X\rangle} \langle Y\rangle$.
 \item
If
$0\neq f\in k\langle X\rangle S$
then
$\overline{f}=a\overline{s}$
for some
$a\in X^*$
and
$s\in S$.
\item
$Irr(S)=\{w \in X^*Y|w\neq a\overline{s}, \ a\in X^*, \ s\in S\}$
is a~linear basis for the quotient
$\Mod_{k\langle X\rangle} \langle Y|S\rangle=
\Mod_{k\langle X\rangle} \langle Y\rangle /k\langle
X\rangle S$.
\end{enumerate}
\end{theorem}

\medskip

Outline of the proof.
Take
$u\in X^*Y$
and express it as
$u=u^Xy_u$
with
$u^X\in X^*$
and
$y_u\in Y$.
Put
$$
  \cm(u,v) = a^Xu=b^Xv,\quad \lcm(u,v) = u=d^Xv,
$$
where
$y_u=y_v$.
Up to the order of
$u$~%
and~%
$v$,
we have
$cm(u,v)=c\cdot \lcm(u,v)$.

The composition of two monic elements
$f,g\in \Mod_{k\langle X\rangle}(Y)$
is
$$
(f,g)|_{\lcm(\bar f, \bar g)}=  \lcm(\bar f, \bar g)|_{\bar f\mapsto
f} - \lcm(\bar f, \bar g)|_{\bar g\mapsto g}.
$$
If
$a_1\bar s_1= a_2\bar s_2$
for monic
$s_1$
and
$s_2$
then
$a_1s_1-a_2s_2 = c\cdot (s_1,s_2)_{\lcm(\bar s_1, \bar s_2)}$.
This gives an~analogue of Lemma \ref{L4.1} for modules
and the implication
(i)$\Rightarrow$(ii)
of Theorem \ref{l3.5}.

\medskip

Given
$S\subset k\langle X\rangle$,
put
$A=k\langle X|S\rangle$.
We can regard every left~%
$A$-module
$_AM$
as a
$k\langle X\rangle$-module
in a~natural way:
$fm:=(f+Id(S))m$
for
$f\in k\langle X\rangle$
and
$m\in M$.
Observe that
$_AM$
is an~epimorphic image of some free~%
$A$-module.
Assume now that
$_A M=
\Mod_A\langle Y|T\rangle=\Mod_{A} \langle Y\rangle/AT$,
where
$T\subset \Mod_{A} \langle Y\rangle$.
Put
$$
T_1=\{\sum f_iy_i\in  \Mod_{k\langle X\rangle} \langle Y\rangle|
\sum (f_i+Id(S))y_i\in T\}
$$
and
$R=SX^*Y\cup T_1$.
Then
$_A M= \mod_{k\langle X\rangle}\!\langle Y|R\rangle$
as
$k\langle X\rangle$-modules.

\begin{theorem}
Given a~submodule~%
$I$
of
$\Mod_{k\langle X\rangle} \langle Y\rangle$
and a~monomial ordering~%
$<$
on
$X^*Y$
as above,
there exists a~unique reduced Gr\"{o}bner--Shirshov basis~%
$S$
for~%
$I$.
\end{theorem}

\begin{corollary}(P.M. Cohn)
Every left ideal~%
$I$
of
$k\langle X\rangle$
is a~free left
$k\langle X\rangle$-module.
\end{corollary}

\textbf{Proof: }
Take a~reduced Gr\"{o}bner--Shirshov basis~%
$S$
of~%
$I$
as a~%
$k\langle X\rangle$-submodule
of the cyclic
$k\langle X\rangle$-module.
Then~%
$I$
is a~free left
$k\langle X\rangle$-module
with a~%
$k\langle X\rangle$-basis~%
$S$.
\hfill
$\blacksquare$

As an~application of the CD-lemma for modules,
we give  GS bases for the Verma modules
over the Lie algebras of coefficients
of free Lie conformal algebras.
We find linear bases for these modules.

Let
$\mathcal {B}$
be a~set of symbols.
Take the constant locality function
$N: \mathcal {B}\times \mathcal {B}\rightarrow \mathbb{Z}_+$;
that is,
$N(a,b)\equiv N$
for all
$a, \ b\in \mathcal {B}$.
Put
$X=\{b(n)| \ b\in \mathcal {B}, \ n\in \mathbb{Z}\}$
and consider the Lie algebra
$L=Lie(X|S)$
over a~field~%
$k$
of characteristic~0
generated by~%
$X$
with the relations
\[
S=\Big\{\sum\limits_{s} (-1)^s
%
%
\tbinom{N}{s}
[b(n-s)a(m+s)]=0| \ a, \ b\in \mathcal {B}, \ m, n\in
\mathbb{Z}\Big\}.
\]
For every
$b\in \mathcal {B}$,
put
$\widetilde{b}=\sum \limits_{n} b(n) z^{-n-1}\in L[[z,z^{-1}]]$.
It is well-known that
these elements generate a~free Lie conformal algebra~%
$C$
with data
$(\mathcal{B},N)$
(see \cite{ro99}).
Moreover,
the coefficient algebra of~%
$C$
is just~%
$L$.

Suppose that~%
$\mathcal {B}$
is linearly ordered.
Define an~ordering on~%
$X$
as
$$
a(m)<b(n)\Leftrightarrow m<n \ \mbox{ or } (m=n \ \mbox{and} \ a<b).
$$
We use the deg-lex ordering on
$X^*$.
It is clear that
the leading term of each polynomial in~%
$S$
is
$b(n)a(m)$
with
$$
n-m>N \ \mbox{or} \ (n-m=N \ \mbox{and} \ (b>a \ \mbox{or} \ (b=a \
\mbox{and} \ N \ \mbox{is odd}))).
$$
The following lemma is essentially from \cite{ro99}.

\begin{lemma}(\cite{CCZ})
With the deg-lex ordering on
$X^*$,
the set
$S$
is a~GS basis in
$Lie(X)$.
\end{lemma}

\begin{corollary}(\cite{CCZ})\label{coro2.16}
A~linear basis of the universal enveloping algebra
$\mathcal {U}=\mathcal U(L)$
of~%
$L$
consists of the monomials
$$
a_1(n_1)a_2(n_2)\cdots a_k(n_k)
$$
with
$a_i\in \mathcal {B}$
and
$n_i\in \mathbb{Z}$
such that for every
$1\leq i<k$
we have
\[
n_i-n_{i+1}\leq \left \{
\begin{array}{ll}
N-1 & \mbox{if
$a_i>a_{i+1}$
or ($a_i=a_{i+1}$
\ and \
$N$
is odd)} \\
 N  & \mbox{otherwise.}
 \end{array}
 \right.
 \]
 \end{corollary}

An~%
$L$-module~%
$M$
is called \textit{restricted}
if for all
$a\in C$
and
$v\in M$
there is some integer~%
$T$
such that
$a(n)v=0$
for
$n\geq T$.

An~%
$L$-module~%
$M$
is called a~\textit{highest weight module}
whenever it is generated over~%
$L$
by a~single element
$m\in M$
satisfying
$L_+ m=0$,
where
$L_+$
is the subspace of~%
$L$
generated by
$\{a(n)|a\in C, n\geq0\}$.
In this case~%
$m$
is called a~\textit{highest weight vector}.

Let us now construct a~universal highest weight module~%
$V$
over~%
$L$,
which is often called the \textit{Verma module}.
Take the trivial~%
$1$-dimensional
$L_+$-module
$kI_v$
generated by
$I_v$;
hence,
$a(n)I_v=0$
for all
$a\in \mathcal {B}, \ n\geq0$.
Clearly,
$$
V=Ind_{L_+} ^L kI_v=\mathcal {U}(L)\otimes_{\mathcal {U}(L_+)}
kI_v\cong \mathcal {U}(L)/\mathcal {U}(L)L_+.
$$
Then~%
$V$
has the structure of the highest weight module over~%
$L$
with the action given by multiplication on
$\mathcal {U}(L)/\mathcal{U}(L)L_+$
and a~highest weight vector
$I\in \mathcal {U}(L)$.
In addition,
$V=\mathcal {U}(L)/\mathcal {U}(L)L_+$
is the universal enveloping vertex algebra of~%
$C$
and the embedding
$\varphi: C\rightarrow V$
is given by
$a\mapsto a(-1)I$
(see also \cite{ro99}).

\begin{theorem}(\cite{CCZ})
With the above notions,
a~linear basis of~%
$V$
consists of the elements
$$
a_1(n_1)a_2(n_2)\cdots a_k(n_k), \ a_i\in \mathcal {B}, \ n_i\in
\mathbb{Z}
$$
satisfying the condition in Corollary \ref{coro2.16}
and
$n_k<0$.
\end{theorem}

\textbf{Proof: }
Clearly,
as
$k\langle X\rangle$-modules,
we have
$$
_{\mathcal {U}} V=_{\mathcal {U}} (\mathcal {U}(L)/\mathcal{U}(L)L_+)=
\Mod_{k\langle X\rangle} \langle I| \ S^{(-)}X^*I, \
a(n)I, \ n\geq0\rangle=_{k\langle X\rangle} \langle I| \ S'\rangle,
$$
where
$S'=\{S^{(-)}X^*I, \ a(n)I, \ n\geq0\}$.
In order to show that
$S'$
is a~Gr\"obner--Shirshov basis,
we only need to verify that
$w=b(n)a(m)I$,
where
$m\geq 0$.
Take
\[
f=\sum\limits_{s} (-1)^s
%
%
\tbinom{n}{s} 
 (b(n-s)a(m+s)-a(m+s)b(n-s))I
\quad \text{and}\quad
 g=a(m)I.
\]
Then
$(f,g)_w=f-b(n)a(m)I\equiv 0 \mod(S',w)$
since
$n-m\geq N$,
$m+s\geq 0$,
$n-s\geq 0$,
and
$0\leq s \leq N$.
It follows that
$S'$
is a~Gr\"obner--Shirshov basis.
Now,
the result follows from the CD-lemma for modules.
\hfill
$\blacksquare$

\subsection{Composition-Diamond lemma for categories }

Denote by~%
$X$
an~oriented multi-graph.
A~path
$$
a_n\rightarrow a_{n-1}\rightarrow\cdots\rightarrow a_1\rightarrow a_0,
\quad
n\geq 0,
$$
in~%
$X$
with edges
$x_n,\dots, x_2, x_1$
is a~partial word
$u = x_1x_2\cdots x_n$
on
$X$
with source
$a_n$
and target
$a_0$.
Denote by
$C(X)$
the free category generated by~%
$X$
(the set of all partial words (paths) on~%
$X$
with partial multiplication,
the free `partial path monoid' on~%
$X$).
A~well-ordering on
$C(X)$
is called monomial
whenever it is compatible with partial multiplication.

A~polynomial
$f \in kC(X)$
is a~linear combination of partial words
with the same source and target.
Then
$kC(X)$
is the \textit{partial path algebra} on~%
$X$
(the free associative partial path algebra generated by~%
$X$).

Given
$S\subset kC(X)$,
denote by
$Id(S)$
the minimal subset of
$kC(X)$
that includes~%
$S$
and is closed under the partial operations of addition and multiplication.
The elements of
$Id(S)$
are of the form
$\sum\alpha_ia_is_ib_i$
with
$\alpha_i\in k$,
$a_i,b_i\in C(X)$,
and
$s_i\in S$,
and all~%
$S$-words
have the same source and target.

Both inclusion and intersection compositions are possible.

With these differences,
the statement and proof of the CD-lemma
are the same as for the free associative algebra.

\begin{theorem}\label{lt1}(\cite{BoChLi-CD-category},
the CD-lemma for categories)
Consider a~nonempty set
$S\subset kC(X)$
of monic polynomials
and a~monomial ordering~%
$<$
on
$C(X)$.
Denote by
$Id(S)$
the ideal of
$kC(X)$
generated by~%
$S$.
The following statements are equivalent:
\begin{enumerate}[\textrm(i)]
\item
The set $S$
is a~Gr\"{o}bner--Shirshov basis in
$kC(X)$.
\item
$f\in Id(S)\Rightarrow \bar f =a\bar s b$
for some
$s\in S$
and
$a,b\in C(X)$.
\item
the set $Irr(S)=\{u\in C(X) |u\ne a\bar s b\ a,b\in C(X),\ s\in S \}$
is a~linear basis for
$kC(X)/Id(S)$,
which is denoted by
$kC(X|S)$.
\end{enumerate}
\end{theorem}

\medskip

Outline of the proof.

Define
$w=\lcm(u,v),\ u,v\in C(X)$
and the general composition
$(f,g)_w$
for
$f,g\in kC(X)$
and
$w= \lcm(\bar f, \bar g)$
by the same formulas as above.
Under the conditions of the analogue of Lemma \ref{L4.1},
we again have
$a_1s_1b_1 - a_2s_2b_2 = c(s_1,s_2)_wd \equiv 0 \mod(S,w)$,
where
$w=\lcm(\bar s_1,\bar s_2)$
and
$c,d\in C(X).$
This implies the analogue of Lemma \ref{L4.1}
and the main assertion
(i)$\Rightarrow$(ii)
of Theorem \ref{lt1}.

\medskip

Let us present some applications of CD-lemma for categories.

For each non-negative integer~%
$p$,
denote by
$[p]$
the set
$\{0,1,2,\ldots,p\}$
of integers in their usual ordering.
A~(weakly) monotonic map
$\mu:[q]\rightarrow [p]$
is a~function from
$[q]$
to
$[p]$
such that
$i\leq j$
implies
$\mu (i)\leq \mu (j)$.
The objects
$[p]$
with weakly monotonic maps
as morphisms
constitute the category
$\Delta$
called the simplex category.
It is convenient to use
two special families of monotonic maps,
$$\varepsilon_q^i:[q-1]\rightarrow [q],\quad
\eta_q^i:[q+1]\rightarrow [q]$$
defined for
$i=0,1,\dots q$
(and for
$q>0$
in the case of
$\varepsilon^i$)
by
\begin{equation*}
\varepsilon_q^i(j)=\left\{
\begin{array}{r@{\quad}l}
j\text{ \ \ \ \   }  & \ \ \ \mbox{ if } \  i>j, \\
j+1 & \ \ \ \mbox{ if } \   i\leq j, \\
\end{array}%
\right.
\end{equation*}

\begin{equation*}
\eta_q^i(j)=\left\{
\begin{array}{r@{\quad}l}
j\text{ \ \ \ \   }  & \ \ \ \mbox{ if } \  i\geq j, \\
j-1 & \ \ \ \mbox{ if } \   j>i. \\
\end{array}%
\right.
\end{equation*}

Take the oriented multi-graph
$X=(V(X),E(X))$
with
\begin{gather*}
V(X)=\{[p]\ |\ p\in Z^+\cup\{0 \}\},
\\
E(X)=\{
\varepsilon_p^i:[p-1]\rightarrow [p], \ \eta_q^j:[q+1]\rightarrow
[q]\ |\ p>0, 0\leq i\leq p, 0\leq j \leq q\}.
\end{gather*}
Consider the relation
$S\subseteq C(X)\times C(X)$
consisting of: 
\begin{eqnarray*}
&&f_{_{q+1,q}}:\ \ \
\varepsilon_{q+1}^i\varepsilon_{q}^{j-1}=\varepsilon_{q+1}^j\varepsilon_{q}^{i}
\quad \text{for } j>i;\\
&&g_{_{q,q+1}}:\ \ \
\eta_{q}^{j}\eta_{q+1}^{i}=\eta_{q}^i\eta_{q+1}^{j+1}
\quad \text{for } j\geq i;\\
&& h_{_{q-1,q}}:\ \ \
\eta_{q-1}^{j}\varepsilon_{q}^i=\left\{
\begin{array}{r@{\quad}l}
\varepsilon_{q-1}^i\eta_{q-2}^{j-1} &
\quad \text{for } j>i,\\
1_{q-1} &
\quad \text{for } i=j \text{ or } i=j+1, \\
\varepsilon_{q-1}^{i-1}\eta_{q-2}^{j} &
\quad \text{for } i>j+1.%
\end{array}%
\right.
\end{eqnarray*}
This yields a~presentation
$\Delta=C(X|S)$
of the simplex category~%
$\Delta$.

Order now
$C(X)$
as follows.

Firstly,
for
$\eta_{p}^{i},\eta_{q}^{j}\in \{\eta_{p}^{i}|p\geq 0, 0\leq i\leq p\}$
put
$\eta_{p}^{i}>\eta_{q}^{j}$
iff
$p>q$
or
($p=q$
and
$i<j$).

Secondly,
for
$$
u=\eta_{p_1}^{i_1}\eta_{p_2}^{i_2}\cdots
\eta_{p_n}^{i_n}\in
\{\eta_{p}^{i}|p\geq 0, 0\leq i\leq p\}^*
$$
(these are all possible words on
$\{\eta_{p}^{i}|p\geq 0, 0\leq i\leq p\}$,
including the empty word
$1_v$,
where
$v\in Ob(X)$),
define
$$
\wt(u)=(n,\eta_{p_n}^{i_n},\eta_{p_{n-1}}^{i_{n-1}},\cdots,\eta_{p_1}^{i_1}).
$$
Then,
for
$u,v\in \{\eta_{p}^{i}|p\geq 0, 0\leq i\leq p\}^*$
put
$u>v$
iff
$\wt(u)>\wt(v)$
lexicographically.

Thirdly,
for
$\varepsilon_{p}^{i},\varepsilon_{q}^{j}\in
\{\varepsilon_p^{i},|p\in Z^+, 0\leq i\leq p\}$,
put
$\varepsilon_{p}^{i}>\varepsilon_{q}^{j}$
iff
$p>q$
or
($p=q$
and
$i<j$).

Finally,
for
$u=v_{0}\varepsilon_{p_{1}}^{i_1}v_1\varepsilon_{p_{2}}^{i_2}\cdots\varepsilon_{p_{n}}^{i_n}
v_{n}\in C(X)$,
where
$n\geq 0$,
and
$v_j\in \{\eta_{p}^{i}|p\geq 0, 0\leq i\leq p\}^*$
put
$\wt(u)=(n,v_0,v_1,\cdots,v_n,\varepsilon_{p_{1}}^{i_1},\cdots,\varepsilon_{p_{n}}^{i_n})$.
Then for every
$u,v\in C(X)$,
$$
 u\succ_{_1}v \Leftrightarrow \wt(u)>\wt(v) \quad \mbox{lexicographically}.
$$
It is easy to check that
$\succ_{_1}$
is a~monomial ordering on
$C(X)$.
Then we have 

\begin{theorem}(\cite{BoChLi-CD-category})
For
$X$~%
and~%
$S$
defined above,
with the ordering
$\succ_1$
on
$C(X)$,
the set~%
$S$
is  a~Gr\"{o}bner--Shirshov basis for the simplex partial path algebra
$kC(X|S)$.
\end{theorem}

\begin{corollary}(\cite{Sa})
Every morphism
$\mu:[q]\rightarrow [p]$
of the simplex category
has a~unique expression of the form
$$
\varepsilon_p^{i_1}\dots\varepsilon_{p-m+1}^{i_m}\eta_{q-n}^{j_1}\dots\eta_{q-1}^{j_n}
$$
with
$p\geq i_1>\dots > i_m\geq 0$,
$0\leq j_1< \dots < j_n <q$,
and
$q-n+m=p$.
\end{corollary}

The cyclic category is defined by generators and relations
as follows,
see \cite{Gelfand}.
Take the  oriented (multi) graph
$Y=(V(Y),E(Y))$
with
$V(Y)=\{[p]\ |\ p\in Z^+\cup\{0 \}\}$
and
$$
E(Y)=\{ \varepsilon_p^i:[p-1]\rightarrow [p], \
\eta_q^j:[q+1]\rightarrow [q],\ t_q:[q]\rightarrow [q] |\ p>0, 0\leq
i\leq p, 0\leq j \leq q\}.
$$
Consider the relation
$S\subseteq C(Y)\times C(Y)$
consisting of:
\begin{eqnarray*}
&&f_{_{q+1,q}}:\ \ \
\varepsilon_{q+1}^i\varepsilon_{q}^{j-1}=\varepsilon_{q+1}^j\varepsilon_{q}^{i}
\quad \text{for } j>i;\\
&&g_{_{q,q+1}}:\ \ \
\eta_{q}^{j}\eta_{q+1}^{i}=\eta_{q}^i\eta_{q+1}^{j+1}
\quad\text{for } j\geq i;\\
&& h_{_{q-1,q}}:\ \ \ \eta_{q-1}^{j}\varepsilon_{q}^i=\left\{
\begin{array}{r@{\quad}l}
\varepsilon_{q-1}^i\eta_{q-2}^{j-1} &
\quad \text{for } j>i,\\
1_{q-1} &
\quad \text{for } i=j \text{ or } i=j+1,\\
\varepsilon_{q-1}^{i-1}\eta_{q-2}^{j} &
\quad \text{for } i>j+1,%
\end{array}%
\right.\\
&&\rho_1:\ \ \ \ \ \
t_q\varepsilon_{q}^{i}=\varepsilon_{q}^{i-1}t_{q-1}
\quad \text{for } i=1,\dots,q;\\
&&\rho_2:\ \ \ \ \ \
t_q\eta_{q}^{i}=\eta_{q}^{i-1}t_{q+1}
\quad \text{for } i=1,\dots,q;\\
&&\rho_3:\ \ \ \ \ \
t_q^{q+1}=1_q.
\end{eqnarray*}

The category
$C(Y|S)$
is called the cyclic category
and denoted by~%
$\Lambda$.

Define an~ordering on
$C(Y)$
as follows.

Firstly,
for
$t_p^i$,
$t_q^j\in \{t_q|q\geq 0\}^{*}$
put
$(t_p)^i> (t_q)^j$
iff
$i>j$
or
($i=j$
and
$p>q$).

Secondly,
 for
$\eta_{p}^{i},\eta_{q}^{j}\in
\{\eta_{p}^{i}|p\geq 0, 0\leq i\leq p\}$
put
$\eta_{p}^{i}>\eta_{q}^{j}$
iff
$p>q$
or
($p=q$
and
$i<j$).

Thirdly,
for
$$
u=w_0\eta_{p_1}^{i_1}w_1\eta_{p_2}^{i_2}\cdots w_{n-1} \eta_{p_n}^{i_n}w_n\in
\{t_q,\eta_{p}^{i}| q,p\geq 0, 0\leq i\leq p\}^*,
$$
where
$w_i\in\{t_q|q\geq 0\}^*$,
put
$$
\wt(u)=(n,w_0,w_1,\cdots,w_n,\eta_{p_n}^{i_n},\eta_{p_{n-1}}^{i_{n-1}},\cdots,\eta_{p_1}^{i_1}).
$$
Then for every
$u,v\in$
$\{t_q,\eta_{p}^{i}|q,p\geq 0, 0\leq i\leq p\}^*$
put
$u>v$
iff
$\wt(u)>\wt(v)$
lexicographically.

Fourthly,
for
$\varepsilon_{p}^{i},\varepsilon_{q}^{j}\in$
$\{\varepsilon_p^{i},|p\in Z^+, 0\leq i\leq p\}$,
$\varepsilon_{p}^{i}>\varepsilon_{q}^{j}$
iff
$p>q$
or
($p=q$
and
$i<j$).

Finally,
for
$u=v_{0}\varepsilon_{p_{1}}^{i_1} v_1\varepsilon_{p_{2}}^{i_2}
\cdots\varepsilon_{p_{n}}^{i_n} v_{n}\in C(Y)$
and
$v_j\in \{t_q,\eta_{p}^{i}|q,p\geq 0, 0\leq i\leq p\}^*$
define
$\wt(u)=(n,v_0,v_1,\cdots,v_n,\varepsilon_{p_{1}}^{i_1},\cdots,\varepsilon_{p_{n}}^{i_n}).$

Then for every
$u,v\in C(Y)$
put
$
 u\succ_{_2}v \Leftrightarrow \wt(u)<\wt(v)\ \ \ \mbox{lexicographically}.
$

It is also easy to verify that
$\succ_{_2}$
is a~monomial ordering on
$C(Y)$
which extends
$\succ_{_1}$.
Then we have 

\begin{theorem}(\cite{BoChLi-CD-category})
Consider
$Y$~%
and~%
$S$
defined as the above.
Put
$\rho_4: \ \
t_q\varepsilon_{q}^0=\varepsilon_{q}^{q}$
and
$\rho_5:\ \
t_q\eta_{q}^0=\eta_{q}^{q}t_{q+1}^{2}$.
Then
\begin{enumerate}[\textrm(1)]
\item
With the ordering
$\succ_{_2}$
on
$C(Y)$,
the set
$S\cup \{\rho_4,\rho_5\}$
is a~Gr\"{o}bner--Shirshov basis for the cyclic category
$C(Y|S)$.
\item
Every morphism
$\mu:[q]\rightarrow [p]$
of the cyclic category
$\Lambda=C(Y|S)$
has a~unique expression of the form
$$
\varepsilon_p^{i_1}\dots\varepsilon_{p-m+1}^{i_m}\eta_{q-n}^{j_1}\dots\eta_{q-1}^{j_n}t_{q}^k
$$
with
$p\geq i_1>\dots> i_m\geq 0$,
$0\leq j_1< \dots < j_n <q$,
$0\leq k \leq q$,
and
$q-n+m=p$.
\end{enumerate}
\end{theorem}

\subsection{Composition-Diamond lemma
for associative algebras over  commutative algebras}

Given two well-ordered sets
$X$~%
and~%
$Y$,
put
$$
N=[X]Y^*=\{u=u^Xu^Y| u^X\in [X] \ and \ u^Y\in Y^*\}
$$
and denote by
$kN$
the~%
$k$-space
spanned by~%
$N$.
Define the multiplication of words as
$$
u=u^Xu^Y,v=v^Xv^Y\in N \Rightarrow
uv=u^Xv^Xu^Yv^Y\in N.
$$
This makes
$kN$
an~algebra isomorphic to the tensor product
$k[ X] \otimes k\langle Y\rangle$,
called a~`double free associative algebra'.
It is a~free object in the category of
all associative algebras over all commutative algebras
(over~%
$k$):
every associative algebra
$_KA$
over a~commutative algebra~%
$K$
is isomorphic to
$k[ X] \otimes k\langle Y\rangle/Id(S)$
as a~%
$k$-algebra
and a~%
$k[X]$-algebra.

Choose a~monomial ordering~%
$>$
on~%
$N$.
The following definitions of compositions
and the GS basis are taken from \cite{MZ}.

Take two monic polynomials
$f$~%
and~%
$g$
in
$k[X]\otimes k\langle Y \rangle$
and denote by~%
$L$
the least common multiple of
$\bar{f}^X$
and
$\bar{g}^X$.

$1.$
 Inclusion.
Assume that
$\bar{g}^Y$
is a~subword of
$\bar{f}^Y$,
say,
$\bar{f}^Y=c\bar{g}^Yd$
for some
$c,d\in Y^*$.
If
$\bar{f}^Y=\bar{g}^Y$
then
$\bar{f}^X\geq \bar{g}^X$
and if
$\bar{g}^Y=1$
then we set
$c=1$.
Put
$w=L\bar{f}^Y=Lc\bar{g}^Yd$.
Define the composition
$C_1(f,g,c)_w=
\frac{L}{\bar{f}^X}f-\frac{L}{\bar{g}^X}cgd$.

$2.$
 Overlap.
Assume that
a~non-empty beginning of
$\bar{g}^Y$
is a~non-empty ending of
$\bar{f}^Y$,
say,
$\bar{f}^Y=cc_0$,
$\bar{g}^Y=c_0d$,
and
$\bar{f}^Yd=c\bar{g}^Y$
for some
$c,d,c_0\in Y^*$
and
$c_0\neq1$.
Put
$w=L\bar{f}^Yd=Lc\bar{g}^Y$.
Define the composition
$C_2(f,g,c_0)_w=
\frac{L}{\bar{f}^X}fd-\frac{L}{\bar{g}^X}cg$.

$3.$
 External.
Take a~(possibly empty) associative word
$c_0\in Y^*$.
In the case that
the greatest common divisor of
$\bar{f}^X$
and
$\bar{g}^X$
is non-empty
and both
$\bar{f}^Y$
and
$\bar{g}^Y$
are non-empty,
put
$w=L\bar{f}^Yc_0\bar{g}^Y$
and define the composition
$
C_3(f,g,c_0)_w=\frac{L}{\bar{f}^X}fc_0\bar{g}^Y-\frac{L}{\bar{g}^X}\bar{f}^Yc_0g
$.

A~monic subset~%
$S$
of
$k[X]\otimes k\langle Y \rangle$
is called a~GS basis
whenever
all compositions of elements of~%
$S$,
say
$(f,g)_w$,
are trivial modulo
$(S,w)$:
$$
(f,g)_w=\sum_i\alpha_ia_is_ib_i,
$$
where
$a_i,b_i\in N$,
$s_i\in S$,
$\alpha_i\in k$,
and
$a_i\bar{s_i}b_i<w$
for all~%
$i$.

\begin{theorem}
(Mikhalev-Zolotykh \cite{MZ,MikhalevZ98},
\label{thm2.22}
the CD-lemma for associative algebras over commutative algebras)
Consider a monic subset
$S\subseteq k[X]\otimes k\langle Y \rangle$
and a~monomial ordering~%
$<$
on~%
$N$.
The following statements are equivalent:
\begin{enumerate}[\textrm(i)]
\item
The set $S$
is a~Gr\"{o}bner--Shirshov basis in
$k[X]\otimes k\langle Y \rangle$.
\item
For every element
$f\in Id(S)$,
the monomial~%
$\bar{f}$
contains
$\bar{s}$
as its subword for some
$s\in S$.
\item
The set $Irr(S)=\{w \in N|w\neq a\overline{s}b, \ a,b\in N, \ s\in S\}$
is a~linear basis for the quotient
$k[X]\otimes k\langle Y \rangle$.
\end{enumerate}
\end{theorem}

\medskip

Outline of the proof.
For
$$
w = \lcm(u,v) = \lcm(u^X,v^X)\lcm(u^Y,v^Y)
$$
the general composition is
$$
(s_1,s_2)_w=(\lcm(u^X,v^X)/u^X)w|_{u\mapsto s_1}-
(\lcm(u^X,v^X)/v^X)w|_{v\mapsto s_2},
$$
where
$s_1,s_2\in k[X]\langle Y\rangle$
are~%
$k$-monic with
$u=\bar s_1$
and
$v=\bar s_2$.
Moreover,
$(s_1,s_2)_w \equiv 0 \mod(\{s_1,s_2\},w)$
whenever
$w=u^Xv^Xu^Yc^Yv^Y$
with
$c^Y\in Y^*$,
that is,~%
$w$
is a~trivial least common multiple relative to both~%
$X$-words
and~%
$Y$-words.
%
This implies the analog of Lemma \ref{L4.1} 
and the claim
(i)$\Rightarrow$(ii)
in Theorem~\ref{thm2.22}.

We apply this lemma in Section 4.3.

\subsection{PBW-theorem for Lie algebras}

Consider a~Lie algebra
$(L,[\ ])$
over a~field~%
$k$
with a~well-ordered linear basis
$X=\{x_i|i\in I\}$
and multiplication table
$S=\{[x_ix_j]=[|x_ix_j|]|i>j,\ i,j\in I\}$,
where for every
$i,j\in I$
we write
$[|x_ix_j|]=\Sigma_t\alpha_{ij}^tx_t$
with
$\alpha_{ij}^t\in k$.
Then
$
U(L)=k\langle X|S^{(-)}\rangle
$
is called the universal enveloping associative algebra of~%
$L$,
where
$S^{(-)}=\{x_ix_j-x_jx_i=[|x_ix_j|]|i>j,\ i,j\in I\}$.

\begin{theorem} (PBW Theorem)
In the above notation
and with the deg-lex ordering on
$X^*$,
the set
$S^{(-)}$
is a~Gr\"{o}bner--Shirshov basis of
$k\langle X\rangle$.
Then by the CD-lemma for associative algebras,
the set
$Irr(S^{(-)})$
consists of the elements
$$
x_{i_1}\dots x_{i_n}
\text{ with } i_1\leq\dots\leq i_n,\ i_1,\dots,i_n\in I,\
n\geq0,
$$
and constitutes a~linear basis of
$U(L)$.
\end{theorem}

\begin{theorem}
\label{thm2.23}
(the PBW Theorem in Shirshov's form)
Consider
$L=Lie(X|S)$
with
$S\subset Lie(X)\subset k\langle X\rangle$
and
$U(L)=k\langle X|S^{(-)}\rangle$.
The following statements are equivalent.
\begin{enumerate}[\textrm(i)]
\item
For the deg-lex ordering,~%
$S$
is a~GS basis of
$Lie(X)$.
\item
For the deg-lex ordering,
$S^{(-)}$
is a~GS basis of
$k\langle X\rangle$.
\item
A~linear basis of
$U(L)$
consists of the words
$u=u_1\cdots u_n$,
where
$u_1\preceq\dots \preceq u_n$
in the lex ordering,
$n\geq0$,
and every
$u_i$
is an~%
$S^{(-)}$-irreducible
associative Lyndon--Shirshov word in~%
$X$.

\item
A~linear basis of~%
$L$
is the set of all~%
$S$-irreducible
Lyndon--Shirshov Lie monomials
$[u]$
in~%
$X$.

\item
A~linear basis of
$U(L)$
consists of the polynomials
$u=[u_1]\cdots [u_n]$,
where
$u_1\preceq\dots \preceq u_n$
in the lex ordering,
$n\geq0$,
and every
$[u_i]$
is an~%
$S$-irreducible non-associative
Lyndon--Shirshov word in~%
$X$.
\end{enumerate}
\end{theorem}

The PBW theorem,
Theorem \ref{thm4.9},
the CD-lemmas for associative and Lie algebras,
Shirshov's factorization theorem,
and property ({\bf VIII}) of Section 4.2
imply that
every LS-subword of~%
$u$
is a~subword of some~%
$u_i$.

L.~Makar-Limanov gave \cite{MakarLimanov94}
an~interesting form of the PBW theorem
for a~finite dimensional Lie algebra.

\subsection{Drinfeld--Jimbo algebra
$U_q(A)$,
Kac--Moody enveloping algebra
$U(A)$,
and the PBW basis of
$U_q(A_N)$}

Take an~integral symmetrizable
$N \times N$
Cartan matrix
$A=(a_{ij})$.
Hence,
$a_{ii}=2$,
$a_{ij} \leq0$
for
$i \neq j$,
and there exists a~diagonal matrix~%
$D$
with diagonal entries
$d_i$,
which are nonzero integers,
such that the product
$DA$
is symmetric.
Fix a~nonzero element~%
$q$
of~%
$k$
with
$q^{4d_i} \neq 1$
for all~%
$i$.
Then the Drinfeld--Jimbo quantum enveloping algebra is
$$
U_q(A)=k \langle X \cup H \cup Y|S^+ \cup K \cup T \cup S^- \rangle,
$$
where
\begin{eqnarray*}
X&=&\{x_i\},\ \ H=\{h_i^{\pm1}\},\ \
Y=\{y_i\},\\
S^+&=&\{ \sum_{\nu=0}^{1-a_{ij}}(-1)^\nu \left(
\begin{array}{c}
1-a_{ij} \\
\nu
\end{array}
 \right)_tx_i^{1-a_{ij}-\nu}x_jx_i^\nu,
\text{ where }
i \neq j, \ t=q^{2d_i}\},\\
S^-&=& \{ \sum\limits_{\nu=0}^{1-a_{ij}}(-1)^\nu\left(
\begin{array}{c}
1-a_{ij} \\
\nu
\end{array}
\right)_ty_i^{1-a_{ij}-\nu}y_jy_i^\nu,
\text{ where }
i \neq j, \ t=q^{2d_i}\},\\
K&=&\{h_ih_j-h_jh_i, h_ih_i^{-1}-1, h_i^{-1}h_i-1,
x_jh_i^{\pm1}-q^{\mp1}
d_ia_{ij}h^{\pm1}x_j, h_i^{\pm1}y_j-q^{\mp1}y_jh^{\pm1}\},  \\
T&=&\{x_iy_j-y_jx_i- \delta_{ij}
\frac{h_i^2-h_i^{-2}}{q^{2d_i}-q^{-2d_i}}\},
\end{eqnarray*}
and
$$
\left(
\begin{array}{c}
m \\
 n \\
 \end{array}
 \right)_t = \left\{
 \begin{array}{cc}
 \prod\limits_{i=1}^n \frac{t^{m-i+1}-t^{i-m-1}}{t^i-t^{-i}} &
\quad(\text{for } m>n>0),\\
1 &
\quad(\text{for } n=0 \text{ or } m=n).
 \end{array}
 \right.
$$

\begin{theorem} (\cite{BoMa})
For every symmetrizable Cartan matrix~%
$A$,
with the deg-lex ordering on
$\{X\cup H\cup Y\}^*$,
the set
$S^{+c}\cup T \cup K \cup S^{-c}$
is a~Gr\"{o}bner--Shirshov basis of the Drinfeld--Jimbo algebra
$U_q(A)$,
where
$S^{+c}$
and
$S^{-c}$
are the Shirshov completions of
$S^{+}$
and
$S^{-}$.
\end{theorem}

\begin{corollary}(Rosso \cite{Rosso89}, Yamane \cite{Ya})
For every symmetrizable Cartan matrix~%
$A$
we have the triangular decomposition
$$
U_q(A)=U_q^+(A)\otimes k[H]\otimes U_q^-(A)
$$
with
$U_q^+(A)=k\langle X|S^+\rangle$
and
$U_q^-(A)=k\langle Y|S^-\rangle$.
\end{corollary}

Similar results are valid for the Kac--Moody Lie algebras
$g(A)$
and their universal enveloping algebras
$$
U(A)=k\langle  X\cup H\cup Y|S^{+}\cup H \cup K \cup S^{-}\rangle,
$$
where
$S^+,\ S^-$
are the same as for
$U_q(A)$,
$$
K=\{h_ih_j-h_jh_i,\
x_jh_i-h_ix_j+ d_ia_{ij}x_i, h_iy_i-y_ih_i+ d_ia_{ij}y_j\},
$$
and
$T=\{x_iy_j-y_jx_i- \delta_{ij}h_i\}$.

\begin{theorem} (\cite{BoMa})
For every symmetrizable Cartan matrix~%
$A$,
the set
$S^{+c}\cup T \cup K \cup S^{-c}$
is a~Gr\"{o}bner--Shirshov basis of the universal enveloping algebra
$U(A)$
of the Kac--Moody Lie algebra
$g(A)$.
\end{theorem}

The PBW theorem in Shirshov's form implies

\begin{corollary}(Kac \cite{Kac})
For every symmetrizable Cartan matrix~%
$A$,
we have the triangular decomposition
$$
U(A)=U^+(A)\otimes k[H]\otimes U^-(A),
\quad
g(A)=g^+(A)\oplus k[H]\oplus g^-(A).
$$
\end{corollary}

E.~N.~Poroshenko \cite{poroshenkoC-D,poroshenkoA}
found GS bases for the Kac--Moody algebras of types
$\widetilde{A_n}$,
$\widetilde{B_n}$,
$\widetilde{C_n}$,
and
$\widetilde{D_n}$.
He used the available linear bases of the algebras \cite{Kac}.

\medskip

Consider now
$$
A=A_N=\left(
  \begin{array}{ccccc}
    2 & -1 & 0 & \cdots & 0 \\
    -1 & 2 & -1 & \cdots & 0 \\
    0 & -1 & 2 & \cdots & 0 \\
    \cdot & \cdot & \cdot & \cdot & \cdot \\
    0 & 0 & 0 & \cdots & 2 \\
  \end{array}
\right)
$$
and assume that
$q^8 \neq1$.
Introduce new variables,
defined by Jimbo (see \cite{Ya}),
which generate
$U_q(A_N)$:
$$
\widetilde{X}=\{x_{ij}, 1 \leq i<j \leq N+1\},
$$
where
$$
x_{ij}=\left\{
\begin{array}{cc}
x_i & \ \ \  j=i+1,\\
qx_{i,j-1}x_{j-1,j}-q^{-1}x_{j-1,j}x_{i,j-1} & \ \ \ j>i+1.
\end{array}
\right.
$$
Order the set
$\widetilde{X}$
as follows:
$
x_{mn}>x_{ij} \Longleftrightarrow(m,n)>_{lex}(i,j).
$
Recall from Yamane \cite{Ya} the notation
\begin{eqnarray*}
C_1&=&\{((i,j),(m,n))|i=m<j<n\},\
C_2=\{((i,j),(m,n))|i<m<n<j\},\\
C_3&=&\{((i,j),(m,n))|i<m<j=n\},\
C_4=\{((i,j),(m,n))|i<m<j<n\},\\
C_5&=&\{((i,j),(m,n))|i<j=m<n\},\ C_6=\{((i,j),(m,n))|i<j<m<n\}.
\end{eqnarray*}
Consider the set
$\widetilde{S}^+$
consisting of Jimbo's relations:
\begin{eqnarray*}
x_{mn}x_{ij}&-&q^{-2}x_{ij}x_{mn}  \ \ \ \ \ \ \ \ \ \ \ \ \ \ \ \
\ \ \ \ \ \ \ \ \ ((i,j),(m,n)) \in C_1 \cup C_3,\\
x_{mn}x_{ij}&-&x_{ij}x_{mn} \ \ \ \ \ \ \ \ \ \ \ \ \ \ \ \ \ \ \ \
\ \ \ \ \ \ \ \ \ ((i,j),(m,n)) \in C_2
\cup C_6,\\
x_{mn}x_{ij}&-&x_{ij}x_{mn}+(q^2-q^{-2})x_{in}x_{mj} \ \ \
((i,j),(m,n)) \in C_4,\\
x_{mn}x_{ij}&-&q^2x_{ij}x_{mn}+qx_{in} \ \ \ \ \ \ \ \ \ \ \ \ \ \ \
\ \ ((i,j),(m,n)) \in C_5.
\end{eqnarray*}
It is easy to see that
$U^+_q(A_N)=k\langle
\widetilde{X}|\widetilde{S^+}\rangle$.

A~direct proof \cite{ChShaoShum} shows that
$\widetilde{S}^+$
is a~GS basis for
$k\langle \widetilde{X}|\widetilde{S^+}\rangle=U^+_q(A_N)$
(\cite{BoMa}).
The proof is different from the argument of
L.~A.~Bokut and~P.~Malcolmson \cite{BoMa}.
This yields

\begin{theorem} (\cite{BoMa})
In the above notation
and with the deg-lex ordering on
$\{\widetilde{X}\cup H\cup \widetilde{Y}\}^*$,
the set
$\widetilde{S}^+\cup T \cup K \cup \widetilde{S}^-$
is a~Gr\"{o}bner--Shirshov basis of
$$
U_q(A_N)=k \langle \widetilde{X}\cup H\cup
\widetilde{Y}|\widetilde{S}^+\cup T \cup K \cup
\widetilde{S}^-\rangle.
$$
\end{theorem}

\begin{corollary}(\cite{Rosso89,Ya})
For
$q^8\neq 1$,
a~linear basis of
$U_q(A_n)$
consists of 
$$
y_{m_1n_1}\cdots y_{m_ln_l}h_1^{s_1}\cdots h_N^{s_N}x_{i_1j_1}\cdots
x_{i_kj_k}
$$
with
$(m_1,n_1)\leq\dots \leq(m_l,n_l)$,
$(i_1,j_1)\leq\dots \leq(i_k,j_k)$,
$k,l\geq0$
and
$s_t\in \mathbb{Z}$.
\end{corollary}

\section{Gr\"obner--Shirshov bases for groups and semigroups}

In this section
we apply the method of GS bases
for braid groups in different sets of generators,
Chinese monoids,
free inverse semigroups,
and plactic monoids in two sets of generators
(row words and column words).

\medskip

Given a~set~%
$X$
consider
$S\subseteq X^*\times X^*$
the congruence
$\rho(S)$
on
$X^*$
generated by~%
$S$,
the quotient semigroup
$$
A=\sgp\langle X|S\rangle=X^*/\rho(S),
$$
and the semigroup algebra
$k(X^*/\rho(S))$.
Identifying the set
$\{u=v|(u,v)\in S\}$
with~%
$S$,
it is easy to see that
$$
\sigma\colon k\langle X|S\rangle\rightarrow k(X^*/\rho(S)),
\quad
\sum\alpha_iu_i+Id(S)\mapsto \sum\alpha_i\overline {u_i}
$$
is an~algebra isomorphism.

The Shirshov completion
$S^c$
of
$S$
consists of semigroup relations,
$S^c=\{u_i-v_i,\ i\in I\}$.
Then
$Irr(S^c)$
is a~linear basis of
$k\langle X|S\rangle$,
and so
$\sigma(Irr(S^c))$
is a~linear basis of
$k(X^*/\rho(S))$.
This shows that
$Irr(S^c)$
consists precisely of
the normal forms of the elements of the semigroup
$\sgp\langle X|S\rangle$.

Therefore,
in order to find the normal forms of the semigroup
$\sgp\langle X|S\rangle$,
it suffices to find a~GS basis
$S^c$
in
$k\langle X|S\rangle$.
In particular,
consider the group
$G=gp\langle X|S\rangle$,
where
$S=\{(u_i,v_i)\in F(X)\times F(X)|i\in I\}$
and
$F(X)$
is the free group on a~set~%
$X$.
Then~%
$G$
has a~presentation
$$
G=\sgp\langle X\cup
X^{-1}|S,x^{\varepsilon}x^{-\varepsilon}=1,\varepsilon=\pm1,x\in
X\rangle,\ \ X\cap X^{-1}=\emptyset
$$
as a~semigroup.

\subsection{Gr\"obner--Shirshov bases for braid groups}

Consider the Artin braid group
$B_{n}$
of type
${{\bf A}_{n-1}}$
(Artin \cite{Artin26}).
We have
$$
B_{n}=gp\langle \sigma_1,\dots, \sigma_n\ |\
\sigma_j\sigma_i=\sigma_i\sigma_j\ (j-1>i),\
\sigma_{i+1}\sigma_i\sigma_{i+1}=\sigma_i\sigma_{i+1}\sigma_i,\
1\leq i\leq n-1\rangle.
$$

\subsubsection{Braid groups in the Artin--Burau generators}

Assume that
$X=Y\dot{\cup}Z$
with
$Y^*$
and~%
$Z$
well-ordered
and that
the ordering on
$Y^*$
is monomial.
Then every word in~%
$X$
has the form
$u=u_0z_1u_1\cdots z_{k}u_{k}$,
where
$k\geq 0$,
$u_i\in Y^*$,
and
$z_{i}\in Z$.
Define the inverse weight of the word
$u\in X^*$
as
$$
\inwt(u)=( k, u_{k}, z_{k}, \cdots,u_1, z_{1}, u_{0} )
$$
and the inverse weight lexicographic ordering as
$$
u>v\Leftrightarrow \inwt(u)>\inwt(v).
$$
Call this ordering the inverse tower ordering for short.
Clearly,
it is a~monomial ordering on~%
$X^*$.

When
$X=Y\dot{\cup}Z$,
$Y=T\dot{\cup} U$,
and
$Y^*$
is endowed with the inverse tower ordering,
define the inverse tower ordering on
$X^*$
with respect to the presentation
$
X=(T\dot{\cup} U)\dot{\cup} Z.
$
In general,
for
$$
X=(\cdots
(X^{(n)}\dot{\cup}X^{(n-1)})\dot{\cup}\cdots)\dot{\cup}X^{(0)}
$$
with
$X^{(n)}$-words
equipped with a~monomial ordering
we can define the inverse tower ordering of~%
$X$-words.

Introduce a~new set of generators for the braid group
$B_{n}$,
called the Artin--Burau generators.
Put
\begin{eqnarray*}
&&s_{i,i+1}=\sigma_{i}^{2},\ \ \
s_{i,j+1}=\sigma_{j}\cdots\sigma_{i+1}\sigma_{i}^{2}\sigma_{i+1}^{-1}\cdots\sigma_{j}^{-1},
\ \ \ \  1\leq i<j\leq n-1;\\
&&\sigma_{i,j+1}=\sigma_{i}^{-1}\cdots\sigma_{j}^{-1},\ \ \ 1\leq
i\leq j\leq n-1;\ \ \ \sigma_{ii}=1,\ \ \ \{a,b\}=b^{-1}ab.
\end{eqnarray*}
Form the sets
$$
S_{j}=\{s_{i,j},s_{i,j}^{-1}, \ 1\leq i,j<n\} \ \mbox{ and } \
\Sigma^{-1}=\{\sigma_{1}^{-1},\cdots\sigma_{n-1}^{-1}\}.
$$
Then the set
$$
S=S_{n}\cup S_{n-1}\cup \cdots \cup S_{2}\cup\Sigma^{-1}
$$
generates
$B_{n}$
as a~semigroup.

Order now the alphabet~%
$S$
as 
$$
S_{n}<S_{n-1}<\cdots < S_{2}<\Sigma^{-1},
$$
and
$$
s_{1,j}^{-1}< s_{1,j}< s_{2,j}^{-1}<\cdots<s_{j-1,j} \ , \ \ \
\sigma_{1}^{-1}<\sigma_{2}^{-1}<\cdots\sigma_{n-1}^{-1}.
$$
Order
$S_{n}$-words
by the \textit{deg-inlex} ordering;
that is,
first compare words by length
and then by the inverse lexicographic ordering
starting from their last letters.
Then we use the inverse tower ordering of~%
$S$-words.

\begin{lemma}
(Artin \cite{Artin47}, Markov \cite{Markov1945})
The following Artin--Markov relations
hold in the braid group
$B_{n}$:
\begin{eqnarray}\label{e77}
&&\sigma_{k}^{-1}s_{i,j}^{\delta}=s_{i,j}^{\delta}\sigma_{k}^{-1}
\text{for } k\neq i-1,i,j-1,j,\\
&&\sigma_{i}^{-1}s_{i,i+1}^{\delta}=s_{i,i+1}^{\delta}\sigma_{1}^{-1},\\
&&\sigma_{i-1}^{-1}s_{i,j}^{\delta}=s_{i-1,j}^{\delta}\sigma_{i-1}^{-1},\\
&&\sigma_{i}^{-1}s_{i,j}^{\delta}=\{s_{i+1,j}^{\delta},s_{i,i+1}\}\sigma_{i}^{-1},\\
&&\sigma_{j-1}^{-1}s_{i,j}^{\delta}=s_{i,j-1}^{\delta}\sigma_{j-1}^{-1},\\
&&\sigma_{j}^{-1}s_{i,j}^{\delta}=\{s_{i,j+1}^{\delta},s_{j,j+1}\}\sigma_{j}^{-1},
\end{eqnarray}
where
$\delta=\pm1$;
\begin{eqnarray}
&&s_{j,k}^{-1}s_{k,l}^{\varepsilon}=\{s_{k,l}^{\varepsilon},s_{j,l}^{-1}\}s_{j,k}^{-1},\\
&&s_{j,k}s_{k,l}^{\varepsilon}=\{s_{k,l}^{\varepsilon},s_{j,l}s_{k,l}\}s_{j,k},\\
&&s_{j,k}^{-1}s_{j,l}^{\varepsilon}=\{s_{j,l}^{\varepsilon},s_{k,l}^{-1}s_{j,l}^{-1}\}s_{j,k}^{-1},\\
&&s_{j,k}s_{j,l}^{\varepsilon}=\{s_{j,l}^{\varepsilon},s_{k,l}\}s_{j,k},\\
&&s_{i,k}^{-1}s_{j,l}^{\varepsilon}=\{s_{j,l}^{\varepsilon},s_{k,l}s_{i,l}s_{k,l}^{-1}s_{i,l}^{-1}\}s_{i,k}^{-1},\\
&&s_{i,k}s_{j,l}^{\varepsilon}=\{s_{j,l}^{\varepsilon},s_{i,l}^{-1}s_{k,l}^{-1}s_{i,l}s_{k,l}\}s_{i,k},
\end{eqnarray}
where
$i<j<k<l$
and
$\varepsilon = \pm1$;
\begin{eqnarray}
&&s_{i,k}^{\delta}s_{j,l}^{\varepsilon}=s_{j,l}^{\varepsilon}s_{i,k}^{\delta},\\
&&\sigma_{j}^{-1}\sigma_{k}^{-1}=\sigma_{k}^{-1}\sigma_{j}^{-1}
\quad\text{for } j<k-1\\
&&\sigma_{j,j+1}\sigma_{k,j+1}=\sigma_{k,j+1}\sigma_{j-1,j}
\quad\text{for } j<k,\\
&&\sigma_{i}^{-2}=s_{i,i+l}^{-1},\\
\label{e87}&&s_{i,j}^{\pm1}s_{i,j}^{\mp1}=1,
\end{eqnarray}
where
$j<i<k<l$
or
$i<k<j<l$,
and
$\varepsilon$,
$\delta=\pm 1$.
\end{lemma}

\begin{theorem}(\cite{Bo-Ch-Shum})
The Artin--Markov relations
(\ref{e77})-(\ref{e87})
form a~Gr\"{o}bner--Shirshov basis of the braid group
$B_{n}$
in terms of the Artin--Burau generators
with respect to the inverse tower ordering of words.
\end{theorem}

It is claimed in \cite{Bo-Ch-Shum} that
some compositions are trivial.
Processing all compositions explicitly,
\cite{ChMo-GSB-Braid-Artin}
supported the claim.

\begin{corollary}
(Markov--Ivanovskii \cite{Artin47})
The  following words are normal forms of the braid group
$B_n$:
$$
f_nf_{n-1}\dots f_2\sigma_{i_nn}\sigma_{i_{n-1}n-1}\dots\sigma_{i_22},
$$
where all
$f_j$
for
$2\leq j\leq n$
are free irreducible words in
$\{s_{ij},\ i<j\}$.
\end{corollary}

\subsubsection{Braid groups in the Artin--Garside generators}

The Artin--Garside generators of the braid group
$B_{n+1}$
are
$\sigma_i,\ 1\leq i\leq n,\ \bigtriangleup,\ \bigtriangleup^{-1}$
(Garside \cite{Garside1969} 1969),
where
$\bigtriangleup=\Lambda_1\cdots\Lambda_n$
with
$\Lambda_i=\sigma_1\cdots\sigma_i$.

Putting
$\bigtriangleup^{-1}<\bigtriangleup< \sigma_1<\dots<\sigma_n$,
order
$\{\bigtriangleup^{-1},\bigtriangleup,\sigma_1,\dots,\sigma_n\}^*$
by the deg-lex ordering.

Denote by
$V(j, i)$,
$W(j, i), \dots$
for
$j\leq i$
positive words in the letters
$\sigma_j,\ \sigma_{j+1},\dots,\sigma_i$,
assuming that
$V(i + 1,i) = 1$,
$W(i + 1,i) = 1,\dots$.

Given
$V = V(1, i )$,
for
$1\leq k\leq n-i$
denote by
$V^{(k)}$
the result of shifting the indices of all letters in~%
$V$
by
$k$:
$\sigma_1\mapsto \sigma_{k+1},\dots, \sigma_i\mapsto\sigma_{k+i}$,
and put
$V'=V^{(1)}$.
Define
$\sigma_{ij}=\sigma_i\sigma_{i-1}\dots\sigma_{j}$
for
$j\leq i-1$,
while
$\sigma_{ii}=\sigma_i$
and
$\sigma_{ii+1}=1$.

\begin{theorem}(\cite{Bo08,BFKS02})\
A~Gr\"{o}bner--Shirshov basis~%
$S$
of
$B_{n+1}$
in the Artin--Garside generators
consists of the following relations:
\begin{eqnarray*}
&&\sigma_{i+1}\sigma_{i}V(1,i-1)W(j, i )\sigma_{i+1j} =
\sigma_{i}\sigma_{i+1}\sigma_{i}V(1,i-1)\sigma_{ij}W( j, i )',\\
&&\sigma_{s}\sigma_{k}=\sigma_{k}\sigma_{s}
\quad\text{for } s-k\geq 2,\\
&&\sigma_{1}V_1\sigma_{2}\sigma_{1}V_2\cdots
V_{n-1}\sigma_{n}\cdots\sigma_{1}= \bigtriangleup
V_1^{(n-1)}V_2^{(n-2)}\cdots V_{(n-1)}',\\
&&\sigma_l\bigtriangleup=\bigtriangleup\sigma_{n-l+1}
\quad\text{ for } 1\leq l\leq n,\\
&&\sigma_l\bigtriangleup^{-1}=\bigtriangleup^{-1}\sigma_{n-l+1}
\quad\text{for } 1\leq l\leq n,\\
&&\bigtriangleup\bigtriangleup^{-1}=1,\
\bigtriangleup^{-1}\bigtriangleup=1,
\end{eqnarray*}
where
$1\leq i \leq n-1$
and
$1\leq j \leq i+1$;
moreover,
$W(j,i)$
begins with
$\sigma_i$
unless it is empty,
and
$V_i = V_i (1, i )$.
\end{theorem}

There are corollaries.

\begin{corollary}
The~%
$S$-irreducible
normal form of each word of
$B_{n+1}$
coincides with its Garside normal form \cite{Garside1969}.
\end{corollary}

\begin{corollary}
(Garside \cite{Garside1969})
The semigroup
$B^+_{n+1}$
of positive braids
can be embedded into a~group.
\end{corollary}

\subsubsection{Braid groups in the Birman--Ko--Lee generators}

Recall that the Birman--Ko--Lee generators
$\sigma_{ts}$
of the braid
group
$B_n$
are 
$$
\sigma_{ts}=(\sigma_{t-1}\sigma_{t-2}\dots\sigma_{s+1})\sigma_{s}
(\sigma_{s+1}^{-1}\cdots\sigma_{t-2}^{-1}\sigma_{t-1}^{-1})
$$
and we have the presentation
\begin{eqnarray*}
B_n&=&gp\langle \sigma_{ts},\ n\geq
t>s\geq1|\sigma_{ts}\sigma_{rq}=\sigma_{rq}\sigma_{ts},\
(t-r)(t-q)(s-r)(s-q)>0,\\
&&\ \ \
\sigma_{ts}\sigma_{sr}=\sigma_{tr}\sigma_{ts}=\sigma_{sr}\sigma_{tr},
n\geq t>s>r\geq1\rangle.
\end{eqnarray*}
Denote by
$\delta$
the Garside word,
$\delta=\sigma_{nn-1}\sigma_{n-1n-2}\cdots\sigma_{21}$.

Define the order as
$\delta^{-1}<\delta<\sigma_{ts}<\sigma_{rq}$
iff
$(t,s) < (r,q)$
lexicographically.
Use the deg-lex ordering on
$\{\delta^{-1},\delta,\sigma_{ts},\ n\geq t>s\geq1\}^*$.

Instead of
$\sigma_{ij}$,
we write simply
$(i, j)$
or
$( j, i)$.
We also set
$$
(t_m, t_{m-1}, \dots, t_1)=
(t_m, t_{m-1})(t_{m-1},t_{m-2})\dots(t_2, t_1),
$$
where
$t_j\neq t_{j+1},\ 1\leq j\leq m-1$.
In this notation,
we can write the defining relations of
$B_n$
as
\begin{gather*}
(t_3, t_2, t_1)=(t_2, t_1, t_3)=(t_1, t_3, t_2)
\quad\text{for } t_3>t_2>t_1,
\\
(k, l)(i, j)=(i, j)(k, l)
\quad\text{for } k>l,\ i > j,\ k > i,
\end{gather*}
where either
$k > i > j >l$
or
$k >l > i > j$.

Denote by
$V_{[t_2,t_1]}$,
where
$n\geq t_2>t_1\geq1$,
a~positive word in
$(k, l)$
satisfying
$t_2\geq k > l \geq t_1$.
We can use any capital Latin letter with indices
instead of~%
$V$,
and appropriate indices
(for instance,
$t_3$
and
$t_0$
with 
$t_3 >t_0$)
instead of
$t_2$
and
$t_1$.
Use also the following equalities in
$B_n$:
$$
V_{[t_2-1,t_1]}(t_2, t_1)=(t_2, t_1)V'_{[t_2-1,t_1]}
$$
for
$t_2 >t_1$,
where
$
V'_{[t_2-1,t_1]}=(V_{[t_2-1,t_1]})
|_{(k,l)\mapsto (k,l), \text{ if } l\neq t_1; \ (k,t_1)\mapsto(t_2,k)};
$
$$
W_{[t_2-1,t_1]}(t_1, t_0)=(t_1, t_0)W^\star_{[t_2-1,t_1]}
$$
for
$t_2>t_1>t_0$,
where
$
W^\star_{[t_2-1,t_1]}=(W_{[t_2-1,t_1]})
|_{(k,l)\mapsto (k,l), \text{ if } l\neq t_1; \ (k,t_1)\mapsto(k,t_0)}.
$

\begin{theorem}(\cite{Bo09})
A~Gr\"{o}bner--Shirshov basis of the braid group
$B_{n}$
in
the Birman--Ko--Lee generators consists of the following relations:
\begin{eqnarray*}
&&(k, l)(i, j) = (i, j)(k, l)
\quad\text{for } k >l > i > j,\\
&&(k, l)V_{[j-1,1]}(i,j) = (i, j)(k, l)V_{[j-1,1]}
\quad\text{for } k > i > j >l,\\
&&(t_3, t_2)(t_2, t_1)=(t_2, t_1)(t_3, t_1),\\
&&(t_3, t_1)V_{[t_2-1,1]}(t_3, t_2)=(t_2, t_1)(t_3, t_1)V_{[t_2-1,1]},\\
&&(t,s)V_{[t_2-1,1]}(t_2, t_1)W_{[t_3-1,t_1]}(t_3, t_1)=
(t_3, t_2)(t,s)V_{[t_2-1,1]}(t_2, t_1)W'_{[t_3-1,t_1]},\\
&&(t_3,s)V_{[t_2-1,1]}(t_2, t_1)W_{[t_3-1,t_1]}(t_3, t_1)=
(t_2,s)(t_3,s)V_{[t_2-1,1]}(t_2, t_1)W'_{[t_3-1,t_1]},\\
&&(2,1)V_{2[2,1]}(3,1)\dots V_{n-1[n-1,1]}(n,1)=\delta V'_{2[2,1]}\dots V'_{n-1[n-1,1]},\\
&&(t,s)\delta=\delta(t+1,s+1),\ \
(t,s)\delta^{-1}=\delta^{-1}(t-1,s-1)
\text{ with } t\pm1,\ s\pm1 \pmod n,\\
&&\delta\delta^{-1}=1,\ \delta^{-1}\delta=1,
\end{eqnarray*}
where
$V_{[k,l]}$
means,
as above,
a~word in
$(i, j)$
satisfying
$k\geq i> j\geq l$,
$t >t_3$,
and
$t_2 > s$.
\end{theorem}

There are two corollaries.

\begin{corollary}(Birman--Ko--Lee \cite{Birman-Ko-Lee})
The semigroup
$B^+_n$
of positive braids in the Birman--Ko--Lee generators
embeds into a~group.
\end{corollary}

\begin{corollary}(Birman--Ko--Lee \cite{Birman-Ko-Lee})
The~%
$S$-irreducible normal form of a~word in
$B_n$
in the Birman--Ko--Lee generators
coincides with the Birman--Ko--Lee--Garside normal form
$\delta^kA$,
where
$A\in B^+_n$.
\end{corollary}

\subsubsection{Braid groups in the Adjan--Thurston generators}

The symmetric group
$S_{n+1}$
has the presentation
\begin{eqnarray*}
S_{n+1}=gp\langle s_1,\dots, s_n\ |\ s_i^2=1, s_js_i=s_is_j\
(j-1>i), s_{i+1}s_is_{i+1}=s_is_{i+1}s_i\rangle.
\end{eqnarray*}
L.~A.~Bokut and L.-S.~Shiao \cite{bs}
found the normal form for
$S_{n+1}$
in the following statement:
the set
$N=\{s_{1i_1}s_{2i_2}\cdots s_{ni_n}|\ i_j\leq j+1\}$
is a~Gr\"{o}bner--Shirshov normal form for
$S_{n+1}$
in the generators
$s_i=(i, i+1)$
relative to the deg-lex ordering,
where
$s_{ji}=s_js_{j-1}\cdots s_i$
for
$j\geq i$
and
$s_{jj+1}=1$.

Take
$\alpha\in S_{n+1}$
with the normal form
$\overline{\alpha}=s_{1i_1}s_{2i_2}\cdots s_{ni_n}\in N$.
Define the length of~%
$\alpha$
as
$|\overline{\alpha}|=l(s_{1i_1}s_{2i_2}\cdots s_{ni_n})$
and write
$\alpha\perp\beta$
whenever
$|\overline{\alpha\beta}|=|\overline{\alpha}|+|\overline{\beta}|$.
Moreover,
every
$\overline{\alpha}\in N$
has a~unique expression
$\overline{\alpha}=s_{_{l_1i_{l_1}}}s_{_{l_2i_{l_2}}}\cdots
s_{_{l_ti_{l_t}}}$
with all
$s_{_{l_ji_{l_j}}}\neq1$.
The number~%
$t$
is called the \textit{breadth} of~%
$\alpha$.

Now put
$$
B'_{n+1}=gp\langle r(\overline{\alpha}),\ \alpha\in
S_{n+1}\setminus\{1\}\ |\
r(\overline{\alpha})r(\overline{\beta})=r(\overline{\alpha\beta}),\
\alpha \perp \beta \rangle,
$$
where
$r(\overline{\alpha})$
stands for a~letter with index~%
$\overline{\alpha}$.

Then for the braid group with~%
$n$
generators
we have
$B_{n+1}\cong B'_{n+1}$.
Indeed,
define
\begin{gather*}
\theta\colon B_{n+1}\rightarrow B'_{n+1}, \quad
\sigma_i\mapsto r(s_i),\\
\theta'\colon B'_{n+1}\rightarrow B_{n+1}, \quad
r(\overline{\alpha})\mapsto \overline{\alpha}|_{s_i\mapsto\sigma_i}.
\end{gather*}
These mappings are homomorphism satisfying
$\theta\theta'={l}_{B'_{n+1}}$
and
$\theta'\theta={l}_{B_{n+1}}$.
Hence,
$$
B_{n+1}=gp\langle r(\overline{\alpha}),\ \alpha\in
S_{n+1}\setminus\{1\}\ |\
r(\overline{\alpha})r(\overline{\beta})=r(\overline{\alpha\beta}),\
\alpha \perp \beta \rangle.
$$

Put
$X=\{r(\overline{\alpha}),\ \alpha\in S_{n+1}\setminus\{1\}\}$.
These generators of
$B_{n+1}$
are called the Adjan--Thurston generators.

Then the positive braid semigroup generated by~%
$X$
is
$$
B_{n+1}^{+}=\sgp\langle X\ |\
r(\overline{\alpha})r(\overline{\beta})=r(\overline{\alpha\beta}),\
\alpha \perp \beta \rangle.
$$

Assume that
$s_1<s_2<\cdots<s_n$.
Define
$r(\overline{\alpha})<r(\overline{\beta})$
if and only if
$|\overline{\alpha}|>|\overline{\beta}|$
or
$|\overline{\alpha}|=|\overline{\beta}|$
and
$\overline{\alpha}<_{lex}\overline{\beta}$.
Clearly,
 this is a~well-ordering on~%
$X$.
We will use the deg-lex ordering on
$X^*$.

\begin{theorem}(\cite{ChZhong-Braid})
\label{t2.4}
The Gr\"{o}bner--Shirshov basis of
$B_{n+1}^{+}$
in the Adjan--Thurston generator~%
$X$
relative to the deg-lex ordering on
$X^*$
consists of the relations
$$
r(\overline{\alpha})r(\overline{\beta})=
r(\overline{\alpha\beta})
\quad\text{for } \alpha \perp \beta;
\qquad
r(\overline{\alpha})r(\overline{\beta\gamma})=
r(\overline{\alpha\beta})r(\overline{\gamma})
\quad\text{for } \alpha \perp \beta \perp \gamma.
$$
\end{theorem}

\begin{theorem}(\cite{ChZhong-Braid})
The Gr\"{o}bner--Shirshov basis of
$B_{n+1}$
in the Adjan--Thurston generator~%
$X$
with respect to the deg-lex ordering on
$X^*$
consists of the relations
\begin{eqnarray*}
&&(1)\ \ \
r(\overline{\alpha})r(\overline{\beta})=r(\overline{\alpha\beta})
\quad\text{for } \alpha \perp \beta, \\
&&(2)\ \ \
r(\overline{\alpha})r(\overline{\beta\gamma})=r(\overline{\alpha\beta})r(\overline{\gamma})
\quad\text{for } \alpha \perp \beta \perp \gamma,\\
&&(3)\ \ \
r(\overline{\alpha})\Delta^{\varepsilon}=\Delta^{\varepsilon}r(\overline{\alpha}')
\quad\text{for } \overline{\alpha}'=\overline{\alpha}|_{s_i\mapsto s_{n+1-i}},\\
&&(4)\ \ \
r(\overline{\alpha\beta})r(\overline{\gamma\mu})=
\Delta r(\overline{\alpha}')r(\overline{\mu})
\quad\text{for } \alpha \perp \beta \perp \gamma \perp \mu
\quad\text{with } r(\overline{\beta\gamma})=\Delta,\\
&&(5)\ \ \
\Delta^{\varepsilon}\Delta^{-\varepsilon}=1.
\end{eqnarray*}
\end{theorem}

\begin{corollary}(Adjan--Thurston)
The normal forms for
$B_{n+1}$
are
$\Delta^k
r(\overline{\alpha_1})\cdots r(\overline{\alpha_s})$
for
$k\in \mathbb{Z}$,
where
$r(\overline{\alpha_1})\cdots r(\overline{\alpha_s})$
is minimal in the deg-lex ordering.
\end{corollary}

\subsection{Gr\"obner--Shirshov basis for the Chinese monoid}

The Chinese monoid
$CH(X,<)$
over a~well-ordered set
$(X, <)$
has the presentation
$
CH(X)=\sgp\langle X| S\rangle
$,
where
$X=\{x_i|i\in I\}$
and
$S$
consists of the relations
\begin{eqnarray*}
&& x_ix_jx_k=x_ix_kx_j= x_jx_ix_k \quad\text{for } i>j>k,\\
 &&x_ix_jx_j= x_jx_ix_j,\ x_ix_ix_j= x_ix_jx_i \quad\text{for } i>j.
\end{eqnarray*}

\begin{theorem}(\cite{Chen-Qiu})\
With the deg-lex ordering on
$X^*$,
the following relations (1)-(5) constitute
a~Gr\"{o}bner--Shirshov basis of the Chinese monoid
$CH(X)$:
\begin{enumerate}[\textrm(1)]
\item
$x_ix_jx_k-x_jx_ix_k$,
\item
$x_ix_kx_j-x_jx_ix_k$,
\item
$x_ix_jx_j-x_jx_ix_j$,
\item
$x_ix_ix_j-x_ix_jx_i$,
\item
$x_ix_jx_ix_k-x_ix_kx_ix_j$,
\end{enumerate}
where
$x_i, x_j, x_k \in X$
and
$i>j>k$.
\end{theorem}

Denote by~%
$\Lambda$
the set consistsing of the words on~%
$X$
of the form
$u_n=w_1w_2\cdots w_n$
with
$n\geq0$,
where
\begin{eqnarray*}
w_1 &=& x_1^{t_{11}}\\
w_2 &=& (x_2x_1)^{t_{21}}x_2^{t_{22}}\\
w_3 &=&  (x_3x_1)^{t_{31}}(x_3x_2)^{t_{32}}x_3^{t_{33}}\\
& & \hspace{0.4cm} \cdots\\
w_n &=&
(x_nx_1)^{t_{n1}}(x_nx_2)^{t_{n2}}\cdots(x_nx_{n-1})^{t_{n(n-1)}}
x_n^{t_{nn}}
\end{eqnarray*}
for
$x_i\in X$
with
$x_1<x_2< \cdots <x_n$,
and all exponents are non-negative.

\begin{corollary}(\cite{Jc01})\
This~%
$\Lambda$
is a set of normal forms of elements of the Chinese monoid
$CH(X)$.
\end{corollary}

\subsection{Gr\"obner--Shirshov basis for free inverse semigroup}

Consider a~semigroup~%
$S$.
An~element
$s\in S$
is called an~inverse of
$t\in S$
whenever
$sts = s$
and
$tst = t$.
An~inverse semigroup is a~semigroup in which
every element~%
$t$
has a~unique inverse,
denoted by
$t^{-1}$.

Given a~set~%
$X$,
put
$X^{-1}=\{x^{-1}|x\in X\}$.
On assuming that
$X\cap X^{-1}=\varnothing$,
denote
$X\cup X^{-1}$
by~%
$Y$.
Define the formal inverses of the elements of
$Y^*$
as
$$
1^{-1}=1,\ (x^{-1})^{-1}=x\ (x\in X),
$$
$$
(y_1y_2\cdots y_n)^{-1}=
y_n^{-1}\cdots y_2^{-1}y_1^{-1}\ (y_1,\ y_2,\cdots,\ y_n\in Y).
$$

It is well known that
$$
\mathcal {FI}(X)=\sgp\langle Y|\ aa^{-1}a=a,\
aa^{-1}bb^{-1}=bb^{-1}aa^{-1}, \ a,b\in Y^* \rangle
$$
is the free inverse semigroup
(with identity)
generated by~%
$X$.

Introduce the notions of a~formal idempotent,
a~(prime) canonical idempotent,
and an~ordered (prime) canonical idempotent in
$Y^*$.
Assume that~%
$<$
is a~well-ordering on~%
$Y$.
\begin{enumerate}[(i)]
\item
 The empty word 1 is an
idempotent.
\item
If~%
$h$
is an~idempotent and
$x\in Y$
then
$x^{-1}hx$
is
both an~idempotent and a~prime idempotent.
\item
If
$e_1, e_2,\cdots, e_m$,
where
$m> 1$,
are prime idempotents
then
$e=e_1e_2\cdots e_m$
is an~idempotent.
\item
An~idempotent
$w\in Y^*$
is called \textit{canonical}
whenever~%
$w$
avoids subwords of the form
$x^{-1}exfx^{-1}$,
where
$x\in Y$,
both
$e$~%
and~%
$f$
are idempotents.
\item
A~canonical idempotent
$w\in Y^*$
is called \textit{ordered}
if every subword
$e=e_1e_2\cdots e_m$
of~%
$w$
with
$m>2$
and
$e_i$
being idempotents
satisfies
$\fir(e_1)<\fir(e_2)<\cdots<\fir(e_m)$,
where
$\fir(u)$
is the first letter of
$u\in Y^*$.
\end{enumerate}

\begin{theorem}(\cite{BoChZhao-inverse-sg})
Denote by~%
$S$
the subset of
$k\langle Y\rangle$
consisting two kinds of polynomials:
\begin{itemize}
\item
$ef-fe$,
where
$e$~%
and~%
$f$
are ordered prime canonical idempotents with
$ef>fe$;
\item
$x^{-1}e'xf'x^{-1}-f'x^{-1}e'$,
where
$x\in Y$,
$x^{-1}e'x$,
and
$xf'x^{-1}$
are ordered prime canonical idempotents.
\end{itemize}
Then,
with the deg-lex ordering on
$Y^*$,
the set~%
$S$
is a~Gr\"{o}ber--Shirshov basis of the free inverse semigroup
$\sgp\langle Y|S\rangle$.
\end{theorem}

\begin{theorem}(\cite{BoChZhao-inverse-sg})
The normal forms of elements of the free inverse semigroup
$\sgp\langle Y|S\rangle$
are
$$
u_0e_1u_1\cdots e_mu_m\in Y^*,
$$
where
$m\geq0$,
$u_1,\cdots, u_{m-1}\neq1$
and
$u_0u_1\cdots u_m$
avoids subwords of the form
$yy^{-1}$
for
$y\in Y$,
while
$e_1,\cdots, e_m$
are ordered canonical idempotents
such that
the first
(respectively last)
letter of 
$e_i$,
for
$1\leq i\leq m$
is not equal to the first
(respectively last)
letter of
$u_i$
(respectively
$u_{i-1}$).
\end{theorem}

The above normal form is analogous to the semi-normal forms
of Poliakova and Schein \cite{Schein}, 2005.

\subsection{Approaches to plactic monoids
via Gr\"obner--Shirshov bases in row and column generators}

Consider the set
$X=\{x_1,\dots,x_n\}$
of~%
$n$
elements
with the ordering
$x_1<\dots<x_n$.
M.~P.~Sch\"utzenberger called
$P_n=\sgp\langle X|T\rangle$
a~plactic monoid
(see also M. Lothaire \cite{Lothaire02}, Chapter 5),
where~%
$T$
consists of the Knuth relations
\begin{eqnarray*}&&x_ix_kx_j= x_kx_ix_j
\quad\text{for } x_i\leq x_j<x_k,\\
&&x_jx_ix_k= x_jx_kx_i
\quad\text{for } x_i< x_j\leq x_k.
\end{eqnarray*}

A~nondecreasing word
$R\in X^*$
is called a~row
and a~strictly decreasing word
$C\in X^*$
is called a~column;
for example,
$x_1x_1x_3x_5x_5x_5x_6$
is a~row and
$x_6x_4x_2x_1$
is a~column.

For two rows
$R,S\in A^*$
say that~%
$R$
dominates~%
$S$
whenever
$|R|\leq|S|$
and every letter of~%
$R$
is greater than
the corresponding letter of~%
$S$,
where
$|R|$
is the length of~
$R$.

A~(semistandard) Young tableau on~%
$A$
(see \cite{M.L})
is a~word
$w=R_1R_2 \cdots R_t$
in
$U^*$
such that
$R_{i}$
dominates
$R_{i+1}$
for all
$i=1,\dots,t-1$.
For example,
$$
x_4x_5x_5x_6\cdot x_2x_2x_3x_3x_5x_7\cdot x_1x_1x_1x_2x_4x_4x_4
$$
is a~Young tableau.

A.~J.~Cain, R.~Gray, and A.~Malheiro \cite{Portugal}
use the Schensted--Knuth normal form
(the set of (semistandard) Young tableaux)
to prove that
the multiplication table of column words,
$uv=u'v'$,
forms a~finite GS basis of
the finitely generated plactic monoid.
Here the Young tableaux
$u'v'$
is the output of the column Schensted algorithm applied to
$uv$,
but
$u'v'$
is not made explicit.

In this section
we give new explicit formulas
for the multiplication tables of row and column words.
In addition,
we give independent proofs that
the resulting sets of relations
are GS bases in row and column generators respectively.
This yields two new approaches to plactic monoids via their GS bases.

\subsubsection{Plactic monoids in the row generators}

Consider the plactic monoid
$P_n=\sgp\langle X|T\rangle$,
where
$X=\{1,2,\dots,n\}$
with
$1<2<\dots<n$.
Denote by~%
$\mathbb{N}$
the set of non-negative integers.
It is convenient
to express the rows
$R\in X^*$
as
$R=(r_1,r_2,\dots, r_n)$,
where
$r_i$
for
$i=1,2,\dots,n$
is the number of occurrences of the letter~%
$i$.
For example,
$R=111225=(3,2,0,0,1,0,\dots,0)$.

Denote by~%
$U$
the set of all rows in
$X^*$
and order
$U^*$
as follows.
Given
$R=(r_1,r_2,\dots, r_n)\in U$,
define the length
$|R|=r_1+\dots+r_n$
of~%
$R$
in~%
$X^*$.

Firstly,
order~%
$U$:
for every
$R, S\in U$,
put
$R<S$
if and  only if
$|R|<|S|$
or
$|R|=|S|$
and
$(r_1,r_2,\dots, r_n)>(s_1,s_2,\dots,s_n)$
lexicographically.
Clearly,
this is a~well-ordering on~%
$U$.
Then,
use the deg-lex ordering on~%
$U^*$.

\begin{lemma}(\cite{ChLiJing-plactic})\label{l3.18}
Take
$\Phi = (\phi_1,\dots,\phi_n)\in U$.
For
$1 \leq p \leq n$
put
$$
\Phi_p=\sum_{i=1}^{p}\phi_i,
$$
where
$\phi_i$
($w_i$,
$z_i$,
$w'_i$,
and
$z'_i$,
see below)
stands for a~lowercase symbol,
and
$\Phi_p$
(
$W_p$,
$Z_p$,
$W'_p$,
and
$Z'_p$,
see below)
for the corresponding uppercase symbol.
Take
$W=(w_1,w_2,\dots,w_n)$
and
$Z=(z_1,z_2,\dots,z_n)$
in~%
$U$.
Put
$W'=(w_1',w_2',\dots,w_n')$
and
$Z'=(z_1',z_2',\dots,z_n')$,
where
\begin{equation*}
w_1'=0,
\quad
w_p'=\min(Z_{p-1}-W'_{p-1},w_p),
\quad
z_q'=w_q+z_q-w'_q
\eqno{(*)}
\label{row-Schensted}
\end{equation*}
for
$n\geq p\geq2$
and
$n\geq q\ge 1$.

Then
$W\cdot Z=W'\cdot Z'$
in
$P_n=\sgp\langle X|T\rangle$
and
$W'\cdot Z'$
is a~Young tableau on~%
$X$,
which could have only one row,
that is,
$Z'= (0,0,\dots,0)$.
Moreover,
$$
P_n=\sgp\langle X|T\rangle\cong \sgp\langle U|\Gamma\rangle,
$$
where
$\Gamma=\{W\cdot Z=W'\cdot Z'\ |\ W,Z\in V\}$.
\end{lemma}

We should emphasize that
$(\ast)$
gives explicitly the product of two rows
obtained by the Schensted row algorithm.

Jointly with our students Weiping Chen and Jing Li
we proved \cite{ChLiJing-plactic},
independently of Knuth's normal form theorem \cite{Knuth},
that~%
$\Gamma$
is a~GS basis of the plactic monoid algebra in row generators
with respect to the deg-lex ordering.
In particular,
this yields a~new proof of Knuth's theorem.

\subsubsection{Plactic monoids in the column generators}

Consider the plactic monoid
$P_n=\sgp\langle X|T\rangle$,
where
$X=\{1,2,\dots,n\}$
with
$1<2<\dots<n$.
Every Young tableaux is a~product of columns.
For example,
$$
4556\cdot 223357\cdot 1112444 =
(421) (521)(531)(632)(54)(74)(4)
$$
is a~Young tableau.

Given a~column
$C\in X^*$,
denote by
$c_i$
the number of occurrences of the letter~%
$i$
in~%
$C$.
Then
$c_i\in \{0,1\}$
for
$i=1,2,\dots,n$.
We write
$C=(c_1;c_2;\dots;c_n)$.
For example,
$C=6421=(1;1;0;1;0;1;0;\dots;0).$

Put
$V=\{C\ |\ C \mbox{ is a~column in }X^*\}$.
For
$R=(r_1;r_2;\dots; r_n)\in V$
define
$\wt(R)=(|R|,r_1,\dots,r_n)$.
Order~%
$V$
as follows:
for
$R, S\in V$,
put
$R<S$
if and only if
$\wt(R)>\wt(S)$
lexicographically.
Then,
use the deg-lex ordering on~%
$V^*$.

For
$\Phi = (\phi_1;\dots;\phi_n)\in V$,
put
$\Phi_p=\sum_{i=1}^{p}\phi_i,\ 1\leq p\leq n$,
where~%
$\phi$
stands for some lowercase symbol defined above
and~%
$\Phi$
stands for the corresponding uppercase symbol.

\begin{lemma}(\cite{ChLiJing-plactic})\label{l3.21}
Take
$W=(w_1;w_2;\dots;w_n)$,
$Z=(z_1;z_2;\dots;z_n)\in V$.
Define
$W'=(w_1';w_2';\dots;w_n')$
and
$Z'=(z_1';z_2';\dots;z_n')$,
where
\begin{equation*}
z_1'=\min(w_1,z_1),
\quad
z_p'=\min(W_p-Z'_{p-1},z_p),
\quad
w_q'=w_q+z_q-z'_q
\eqno{(**)}
\end{equation*}
for
$n\geq p\ge 2$
and
$n\geq q\ge 1$.
Then
$W',Z'\in V$
and
$W\cdot Z=W'\cdot Z'$
in
$P_n=\sgp\langle X|T\rangle$,
and
$W'\cdot Z'$
is a~Young tableau on~%
$X$.
Moreover,
$$
P_n=\sgp\langle X|T\rangle\cong \sgp\langle V|\Lambda\rangle,
$$
where
$\Lambda=\{W\cdot Z=W'\cdot Z'\ |\ W,Z\in V\}$.
\end{lemma}

Eq.~$(\ast\ast)$
gives explicitly the product of two columns
obtained by the Schensted column algorithm.

Jointly with our students Weiping Chen and Jing Li
we proved \cite{ChLiJing-plactic},
independently of Knuth's normal form theorem \cite{Knuth},
that~%
$\Lambda$
is a~GS basis of the plactic monoid algebra in column generators
with respect to the deg-lex ordering.
In particular,
this yields another new proof of Knuth's theorem.
Previously
Cain, Gray, and Malheiro \cite{Portugal} established the same result
using Knuth's theorem,
and they did not find~%
$\Lambda$
explicitly.

\medskip

\noindent{\bf Remark}:
All results of \cite{ChLiJing-plactic} are valid for
every plactic monoid,
not necessarily finitely generated.

\section{Gr\"{o}bner--Shirshov bases for Lie algebras}

In this section
we first give a~different approach to the LS basis and the Hall basis
of a~free Lie algebra
by using Shirshov's CD-lemma for anti-commutative algebras.
Then,
using the LS basis,
we construct the classical theory of GS bases
for Lie algebras over a~field.
Finally,
we mention GS bases for Lie algebras over a~commutative algebra
and give some interesting applications.

\subsection{Lyndon--Shirshov  basis and Lyndon--Shirshov words
in anti-commutative algebras}

A~linear space~%
$A$
equipped with a~bilinear product
$x\cdot y$
is called an~anti-commutative algebra
if it satisfies the identity
$x^2=0$,
and so
$x\cdot y=-y\cdot x$
for every
$x,y\in A$.

Take a~well-ordered set~%
$X$
and denote by
$X^{**}$
the set of all non-associative words.
Define three orderings 
$\succ_{lex}$,
$>_{_{deg-lex}}$,
and
$>_{_{n-deg-lex}}$
(non-associative deg-lex)
on
$X^{**}$.
For
$(u),(v)\in X^{**}$
put
\begin{itemize}
\item
$(u)=((u_1)(u_2))\succ_{lex}(v)=((v_1)(v_2))$
(here
$(u_2)$
or
$(v_2)$
is empty when
$|(u)|=1$
or
$|(v)|=1$)
iff one of the following holds:
\begin{enumerate}[(a)]
\item
$u_1u_2>v_1v_2$
in the lex ordering;
\item
$u_1u_2=v_1v_2$
and
$(u_1)\succ_{lex}(v_1)$;
\item
$u_1u_2=v_1v_2$,
$(u_1)=(v_1)$,
and
$(u_2)\succ_{lex}(v_2)$;
\end{enumerate}
\item
$(u)=((u_1)(u_2))>_{deg-lex}(v)=((v_1)(v_2))$
iff one of the following holds:
\begin{enumerate}[(a)]
\item
$u_1u_2>v_1v_2$
in the deg-lex ordering;
\item
$u_1u_2=v_1v_2$
and
$(u_1)>_{deg-lex}(v_1)$;
\item
$u_1u_2=v_1v_2$,
$(u_1)=(v_1)$,
and
$(u_2)>_{deg-lex}(v_2)$;
\end{enumerate}
\item
$(u)>_{_{n-deg-lex}}(v)$
iff one of the following holds:
\begin{enumerate}[(a)]
\item
$|(u)|>|(v)|$;
\item
if
$|(u)|=|(v)|$,
$(u)=((u_1)(u_2))$,
and
$(v)=((v_1)(v_2))$
then
$(u_1)>_{_{n-deg-lex}}(v_1)$
or
($(u_1)=(v_1)$
and
$(u_2)>_{_{n-deg-lex}}(v_2)$).
\end{enumerate}
\end{itemize}

\medskip

Define regular words
$(u)\in X^{**}$
by induction on
$|(u)|$:

(i)
$x_i\in X$
is a~regular word.

(ii)
$(u)=((u_1)(u_2))$
is regular if both
$(u_1)$
and
$(u_2)$
are
regular and
$(u_1)\succ_{lex}(u_2)$.

Denote
$(u)$
 by
$[u]$
whenever
$(u)$
is regular.

The set
$N(X)$
of all regular words on~%
$X$
constitutes a~linear basis of the free anti-commutative algebra
$AC(X)$
on~%
$X$.

The following result gives an~alternative approach to
the definition of LS words as the radicals of associative supports~%
$u$
of the normal words~%
$[u]$.

\begin{theorem} (\cite{BoChLi10})
Suppose that
$[u]$
is a~regular word of the anti-commutative algebra
$AC(X)$.
Then
$u=v^m$,
where~%
$v$
is a~Lyndon--Shirshov word in~%
$X$
and
$m\geq1$.
Moreover,
the set of associative supports of the words in
$N(X)$
includes the set of all Lyndon--Shirshov words in~%
$X$.
\end{theorem}

Fix an~ordering
$>_{deg-lex}$
on
$X^{**}$
and choose monic polynomials
$f$~%
and~%
$g$
in
$AC(X)$.
If there exist
$a,b\in {X^*}$
such that
$[w]=[\bar f] =[a[\bar g] b]$
then the inclusion composition of
$f$~%
and~%
$g$
is defined as
$(f,g)_{[w]}=f-[a[g]b]$.

A~monic subset~%
$S$
of
$AC(X)$
is called a~GS basis in
$AC(X)$
if every inclusion composition
$(f,g)_{[w]}$
in~%
$S$
is trivial modulo
$(S,[w])$.

\begin{theorem}\label{cd anti-comm}
(Shirshov's CD-lemma for anti-commutative algebras, cf.~\cite{Sh62a})
Consider a~nonempty set
$S\subset AC(X)$
of monic polynomials with the ordering
$>_{deg-lex}$
on
$X^{**}$.
The following statements are equivalent:
\begin{enumerate}[\textrm(i)]
\item
The set~%
$S$
is a~Gr\"{o}bner--Shirshov basis in
$AC(X)$.
\item
If
$f\in Id(S)$
then
$[\bar f] =[a[\bar s ]b]$
for some
$s\in S\
and\ a,b\in X^*$,
where
$[as b]$
is a~normal
$S$-word.
\item
The set
$$
Irr(S)=\{[u]\in N(X) |\ [u]\ne [a[\bar s] b]\ a,b\in X^*,\ s\in S
\mbox{ and }
[as b] \mbox{ is a~normal } S\mbox{-word}\}
$$
is a~linear basis of the algebra
$AC(X|S)=AC(X)/Id(S)$.
\end{enumerate}
\end{theorem}

Define the subset
$S_1$
the free anti-commutative algebra
$AC(X)$
as 
\begin{eqnarray*}
S_1&=&\{([u][v])[w]-([u][w])[v]-[u]([v][w]) \ |
\\
&&\ \ \ [u],[v],[w]\in N(X)\ \mbox{\ and \ }
[u]\succ_{lex}[v]\succ_{lex}[w]\}.
\end{eqnarray*}

It is easy to prove that the free Lie algebra
admits a~presentation as an~anti-commutative algebra:
$Lie(X)=AC(X)/Id(S_1)$.

The next result gives an~alternating approach to
the definition of the LS basis of a~free Lie algebra
$Lie(X)$
as a~set of irreducible non-associative words
for an~anti-commutative GS basis in
$AC(X)$.

\begin{theorem}(\cite{BoChLi10})
Under the ordering
$>_{deg-lex}$,
the subset
$S_1$
of
$AC(X)$
is an~anti-commutative Gr\"{o}bner--Shirshov basis in
$AC(X)$.
Then
$Irr(S_1)$
is the set of all non-associative LS words in~%
$X$.
So,
the LS monomials
constitute a~linear basis of the free Lie algebra
$Lie(X)$.
\end{theorem}

\begin{theorem}(\cite{BoChLi09})
Define
$S_2$
by analogy with
$S_1$,
but using
$>_{_{n-deg-lex}}$
instead of
$\succ_{lex}$.
Then with the ordering
$>_{_{n-deg-lex}}$
the subset
$S_2$
of
$AC(X)$
is also an~anti-commutative GS basis.
The set
$Irr(S_2)$
amounts to the set of all Hall words in~%
$X$
and forms a~linear basis of a~free Lie algebra
$Lie(X)$.
\end{theorem}

\subsection{Composition-Diamond lemma for Lie algebras over a~field}

We start with some concepts and results from the literature
concerning the theory of GS bases for the free Lie algebra
$Lie(X)$
generated by~%
$X$
over a~field~%
$k$.

Take a~well-ordered set
$X=\{x_i |i\in I\}$
with
$x_i>x_t$
whenever
$i>t$,
for all
$i,t\in I$.
Given
$u=x_{i_1}x_{i_2} \cdots x_{i_m} \in X^*$,
define the length (or degree) of~%
$u$
to be~%
$m$
and denote it by
$|u|=m$
or
$deg(u)=m$,
put
$\fir(u)=x_{i_1}$,
and introduce
\begin{gather*}
x_\beta=
\min(u)=
\min\{x_{i_1},x_{i_2}, \cdots ,x_{i_m}\},\\
X'(u)=\{x_i^j=x_i\underbrace{x_\beta \cdots x_\beta}_{j}| i>\beta, j\geq0\}.
\end{gather*}
Order the new alphabet
$X'(u)$
as follows:
$$
x_{i_1}^{j_1}>x_{i_2}^{j_2} \ \Leftrightarrow\
i_1>i_2 \text{ or } i_1=i_2 \text{ and } j_2>j_1.
$$
Assuming that
$$
u=x_{r_1}\underbrace{x_\beta \cdots x_\beta}_{m_1} \cdots
x_{r_t}\underbrace{x_\beta \cdots x_\beta}_{m_t},
$$
where
$r_i>\beta$,
define the Shirshov elimination
$$
u'=x_{r_1}^{m_1} \cdots x_{r_t}^{m_t} \in (X'{(u)})^*.
$$

We use two linear orderings on
$X^*$:
\begin{enumerate}[(i)]
\item
the lex ordering (or lex-antideg ordering):
$1\succ v$
if
$v\neq1$
and,
by induction,
if
$u=x_{i}u_1$
and
$v=x_{j}v_1$
then
$u\succ v$
if and only if
$x_{i}>x_{j}$
or
$x_{i}=x_{j}$
and
$u_1\succ v_1$;
\item
the deg-lex ordering:
$u >v$
if
$|u|>|v|$
or
$|u|=|v|$
and
$u\succ v$.
\end{enumerate}

\noindent{\bf Remark}
In commutative algebras,
the lex ordering
is understood to be the lex-deg ordering with the condition
$v>1$
for
$v\neq1$.

\medskip

We cite some useful properties of ALSWs and NLSWs
(see below)
following Shirshov \cite{Sh58,Sh58b,Sh62b},
see also \cite{Shirshov-Selected}.
Property~({\bf X}) was given by
Shirshov \cite{Sh58b} and Chen-Fox-Lyndon \cite{CFL}.
Property ({\bf VIII})
was implicitly used in Shirshov \cite{Sh62b},
see also Chibrikov \cite{Chibrikov12}.

We regard
$Lie(X)$
as the Lie subalgebra of the free associative algebra
$k \langle X \rangle$
generated by~%
$X$
with the Lie bracket
$[u,v]=uv-vu$.
Below  we prove that
$Lie(X)$
is the free Lie algebra generated by~%
$X$
for every commutative ring~%
$k$
(Shirshov \cite{Sh58}).
For a~field,
this follows from the PBW theorem
because the free Lie algebra
$Lie(X)=Lie(X|\emptyset)$
has the universal enveloping associative algebra
$k \langle X \rangle=k \langle X|\emptyset \rangle$.

Given
$f \in k \langle X \rangle$,
denote by~%
$\bar{f}$
the leading word of~%
$f$
with respect to the deg-lex ordering
and write
$f=\alpha_{\bar f}\bar f-r_{_f}$
with
$\alpha_{\bar f}\in k$.

\begin{definition}(\cite{Lyndon,Sh58})
Refer to
$w\in X^*\setminus\{1\}$
as an~associative Lyndon--Shirshov word,
or ALSW for short,
whenever
$$
(\forall u,v\in X^*, u,v\neq1) \ w=uv\Rightarrow w>vu.
$$
Denote the set of all ALSWs on~%
$X$
by
$ALSW(X)$.
\end{definition}

Associative Lyndon--Shirshov words enjoy the following properties
(Lyndon \cite{Lyndon},
Chen--Fox--Lyndon \cite{CFL},
Shirshov \cite{Sh58,Sh58b}).

({\bf I})
Put
$x_{\beta} = \min(uv)$.
If
$\fir(u)\neq x_{\beta}$
and
$\fir(v)\neq x_{\beta}$
then
$$
u \succ v\ (\mbox{in the lex ordering on } X^*) \Leftrightarrow
u'\succ v'\ (\mbox{in the lex ordering on } (X'_{uv})^*).
$$

({\bf II})
(Shirshov's key property of ALSWs)
A~word~%
$u$
is an~ALSW in
$X^*$
if and only if
$u'$
is an~ALSW in
$(X'{(u)})^*$.

Properties ({\bf I}) and ({\bf II})
enable us to prove the properties of ALSWs and NLSWs (see below)
by induction on length.

({\bf III})
(down-to-up bracketing)
$u\in ALSW(X)\Leftrightarrow
(\exists k)\ |u^{(k)}|_{_{(X(u))^{(k)}}}=1$,
where
$u^{(k)}=(u')^{(k-1)}$
and
$(X(u))^{(k)}=(X'(u))^{(k-1)}$.
In the process
$u\rightarrow u'\rightarrow u''\rightarrow\cdots$
we use the algorithm of joining the minimal letters of
$u$,
$u'\dots$
to the previous words.

({\bf IV})
If
$u,v\in ALSW(X)$
then
$uv\in ALSW(X)\Leftrightarrow u\succ v$.

({\bf V})
$w\in ALSW(X)\Leftrightarrow$
(for every
$u,v\in X^*\setminus\{1\}$
and
$w=uv\Rightarrow w\succ v$).

({\bf VI})
If
$w\in ALSW(X)$
then an~arbitrary proper prefix of~%
$w$
cannot be a~suffix of~%
$w$
and
$wx_{_{\beta}}\in ALSW(X)$
if
$x_{_{\beta}} = \min(w)$.

({\bf VII})
(Shirshov's factorization theorem)
Every associative word
$w$
can be uniquely represented as
$w=c_{1}c_{2} \ldots c_{n}$,
where
$c_{1}, \ldots ,c_{n}\in ALSW(X)$
and
$c_{1} \preceq c_{2} \preceq \ldots \preceq c_{n}$.

Actually,
if we apply to
$w$
the algorithm of joining the minimal letter to the previous word
using the Lie product,
$w\rightarrow w'\rightarrow w''\rightarrow\cdots$,
then after finitely many steps
we obtain
$w^{(k)}=[c_{1}][c_{2}] \ldots [c_{n}]$,
with
$c_{1} \preceq c_{2} \preceq \ldots \preceq c_{n}$,
and
$w=c_{1}c_{2} \ldots c_{n}$
would be the required factorization
(see an~example in the Introduction).

({\bf VIII})
If an~associative word~%
$w$
is represented as in ({\bf VII})
and~%
$v$
is a LS subword of~%
$w$
then~%
$v$
is a~subword of one of the words
$c_{1}$,
$c_{2},\ldots,c_{n}$.

({\bf IX})
If
$u_1u_2$
and
$u_2u_3$
are ALSWs
then so is
$u_1u_2u_3$
provided that
$u_2\neq1$.

({\bf X})
If
$w=uv$
is an~ALSW
and~%
$v$
is its longest proper ALSW ending,
then~%
$u$
is an~ALSW as well
(Chen--Fox--Lyndon \cite{CFL}, Shirshov \cite{Sh58b}).

\begin{definition}
(down-to-up bracketing of ALSW,
Shirshov \cite{Sh58})
For an~ALSW
$w$,
there is the down-to-up bracketing
$w\rightarrow w'\rightarrow w''\rightarrow\cdots\rightarrow w^{(k)}=[w]$,
where each time we join the minimal letter of the previous word
using Lie multiplication.
To be more precise,
we use the induction
$[w]=[w']_{_{x_i^j\mapsto [[x_ix_{\beta}]\cdots x_{\beta}]}}$.
\end{definition}

\begin{definition}
(up-to-down bracketing of ALSW,
Shirshov \cite{Sh58b},
Chen--Fox--Lyndon \cite{CFL})
For an~ALSW~%
$w$,
we define the up-to-down Lie bracketing
$[[w]]$
by the induction
$[[w]]=[[[u]][[v]]]$,
where
$w=uv$
as in ({\bf X}).
\end{definition}

({\bf XI})
If
$w\in ALSW(X)$
then
$[w]=[[w]]$.

({\bf XII})
Shirshov's definition of a NLSW
(non-associative LS word)
$(w)$
below
is the same as
$[w]$
and
$[[w]]$;
that is,
$(w)=[w]=[[w]]$.
Chen, Fox, and Lyndon \cite{CFL}
used
$[[w]]$.

\begin{definition} (Shirshov \cite{Sh58})
A~non-associative word
$(w)$
in~%
$X$
is a  NLSW
if
\begin{enumerate}[\rm(i)]
\item
$w$
is an~ALSW;
\item
if
$(w)=((u)(v))$
then both
$(u)$
and
$(v)$
are NLSWs
(then {\rm ({\bf IV})} implies that
$u\succ v$);
\item
if
$(w)=(((u_1)(u_2))(v))$
then
$u_2 \preceq v$.
\end{enumerate}
Denote the set of all NLSWs on
$X$
by
$NLSW(X)$.
\end{definition}

({\bf XIII})
If
$u\in ALSW(X)$
and
$[u]\in NLSW(X)$
then
$\overline{[u]}=u$
in
$k \langle X \rangle$.

({\bf XIV})
The set
$NLSW(X)$
is linearly independent in
$Lie(X)\subset k \langle X \rangle$
for every commutative ring~%
$k$.

({\bf XV})
$NLSW(X)$
is a~set of linear generators in every Lie algebra generated by~%
$X$
over an~arbitrary commutative ring~%
$k$.

({\bf XVI})
$Lie(X)\subset k \langle X \rangle$
is the free Lie algebra over the commutative ring~%
$k$
with the~%
$k$-basis
$NLSW(X)$.

({\bf XVII})
(Shirshov's special bracketing \cite{Sh58})
Consider
$w=aub$
with
$w,u\in ALSW(X)$.
Then
\begin{enumerate}[(i)]
\item
$
[w]=[a[uc]d],
$
where
$b=cd$
and possibly
$c=1$.
\item
Express~%
$c$
in the form
$
c=c_{1}c_{2} \ldots c_{n},
$
where
$c_{1}, \ldots ,c_{n}\in ALSW(X)$
and
$c_{1} \preceq  c_{2} \preceq \ldots \preceq c_{n}$.
Replacing
$[uc]$
by
$[\ldots[[u][c_{1}]] \ldots [c_{n}]]$,
we obtain the word
$$
[w]_{u}=[a[\ldots[[[u][c_{1}]][c_2]] \ldots [c_{n}]]d]
$$
which is called the Shirshov special bracketing of~%
$w$
relative to~%
$u$.
\item
$[w]_u=a[u]b+\sum\limits_{i}\alpha_ia_i[u]b_i$
in
$k \langle X \rangle$
with
$\alpha_i\in{k}$
and
$a_i,b_i\in X^*$
satisfying
$a_iub_i<aub$,
and hence
$\overline{[w]}_{u}=w$.
\end{enumerate}

Outline of the proof.
Put
$x_{\beta} = \min(w)$.
Then
$w'=a'(ux_{\beta}^m)'(b_1)'$
in
$(X(w)')^*$,
where
$b=x_{\beta}^mb_1$
and
$ux_{\beta}^m$
is an~ALSW.
Claim~(i) follows from ({\bf II}) by induction on length.
The same applies to claim (iii).

\medskip

({\bf XVIII})
(Shirshov's Lie elimination of the leading word)
Take two monic Lie polynomials
$f$~%
and~%
$s$
with
$\bar f=a\bar sb$
for some
$a,b\in X^*$.
Then
$f_1=f-[asb]_{\bar s}$
is a~Lie polynomial with smaller leading word,
and so
$\bar f_1<\bar f$.

({\bf XIX})
(Shirshov's double special bracketing)
Assume that
$w=aubvc$
with
$w,u,v\in ALSW(X)$.
Then there exists a~bracketing
$[w]_{u,v}$
such that
$[w]_{u,v}=[a[u]b[v]c]_{u,v}$
and
$\overline{[w]_{u,v}}=w$.

More precisely,
$[w]_{u,v} = [a[up]_uq[vr]_vs]$
if
$[w]=[a[up]q[vr]s]$,
and
$$
[w]_{u,v}= [a[\ldots[\ldots[[u][c_1]]\ldots[c_i]_v]\ldots[c_n]]p]
$$
if
$[w]=[a[uc]p]$,
where
$c=c_1\dots c_n$
is the Shirshov factorization of~%
$c$
and~%
$v$
is a~subword of~%
$c_i$.
In both cases
$[w]_{u,v} = a[u]b[v]d + \sum \alpha_i a_i[u]b_i[v]d_i$
in
$k \langle X \rangle$,
where
$a_iub_ivd_i < w$.

({\bf XX})
(Shirshov's algorithm for recognizing Lie polynomials,
cf.\ the Dynkin--Specht--Wever and Friedrich algorithms).
Take
$s\in Lie(X)\subset k\langle X\rangle$.
Then
$\bar s$
is an~ALSW and
$s_1=s-\alpha_{\bar s}[\bar s]$
is a~Lie polynomial with a~smaller maximal word
(in the deg-lex ordering),
$\bar s_1<\bar s$,
where
$s=\alpha_{\bar s}[\bar s]+\dots$.
Then
$s_2=s_1-\alpha_{\bar
s_1}[\bar s_1],\ \overline{s_2}<\overline{s_1}$.
Consequently,
$s\in Lie(X)$
if and only if
after finitely many steps we obtain
$$
s_{m+1}=s-\alpha_{\bar s}[\bar s]-\alpha_{\bar s_1}[\bar
s_1]-\dots-\alpha_{\bar s_m}[\bar s_m]=0.
$$
Here~%
$k$
can be an~arbitrary commutative ring.

\begin{definition}
Consider
$S\subset Lie(X)$
with all
$s\in S$
monic.
Take
$a,b\in{X^*}$
and
$s\in S$.
If
$a\bar{s}b$
is an~ALSW
then we call
$[asb]_{\bar{s}}=[a\bar{s}b]_{\bar{s}}|_{[\bar{s}]\mapsto{s}}$
a~special normal~%
$S$-word
(or a~special normal
$s$-word),
where
$[a\bar{s}b]_{\bar{s}}$
is defined in {\rm({\bf XVII})~(ii)}.
A~Lie~%
$S$-word
$(asb)$
is called a~normal~%
$S$-word
whenever
$\overline{(asb)}=a\overline{s}b$.
Every special normal~%
$s$-word
is a~normal~%
$s$-word by {\rm({\bf XVII})~(iii)}.

For
$f,\ g\in S$
there are two kinds of Lie compositions:
\begin{enumerate}[\rm(i)]
\item
If
$w=\bar{f}=a\bar{g}b$
for some
$a,b\in X^*$
then the polynomial
$\langle f,g\rangle_w=f - [agb]_{\bar{g}}$
is called the inclusion composition of
$f$~%
and~%
$g$
with respect to~%
$w$.
\item
If~%
$w$
is a~word satisfying
$w=\bar{f}b=a\bar{g}$
for some
$a,b\in X^*$
with
$deg(\bar{f})+deg(\bar{g})>deg(w)$
then the polynomial
$\langle f,g\rangle_w=[fb]_{\bar{f}}-[ag]_{\bar{g}}$
is called the intersection composition of
$f$~%
and~%
$g$
with respect to~%
$w$,
and~%
$w$
is an~ALSW by~{\rm({\bf IX})}.
\end{enumerate}

Given a~Lie polynomial~%
$h$
and
$w\in X^*$,
say that~%
$h$
is trivial modulo
$(S,w)$
and write
$h\equiv_{Lie}0 \  mod(S,w)$
whenever
$h=\sum_{i}\alpha_i(a_is_ib_i)$,
where each
$\alpha_i\in k,\ (a_is_ib_i)$
is a~normal~%
$S$-word
and
$a_i\overline{s_i}b_i<w$.

A~set~%
$S$
is called a~GS basis in
$Lie(X)$
if every composition
$(f,g)_w$
of polynomials
$f$~
and~
$g$
in~%
$S$
is trivial modulo~%
$S$
and~%
$w$.
\end{definition}

\medskip

({\bf XXI})
If
$s\in Lie(X)$
is monic
and
$(asb)$
is a~normal~%
$S$-word
then
$(asb)=asb +\sum_{i}\alpha_ia_isb_i$,
where
$a_i\overline{s}b_i<a\overline{s}b$.

A~proof of ({\bf XXI}) follows from the CD-lemma for associative algebras
since
$\{s\}$
is an~associative GS basis by ({\bf IV}).

({\bf XXII})
Given two monic Lie polynomials
$f$~%
and~%
$g$,
we have
$$
\langle f,g\rangle_w-(f,g)_w\equiv_{ass}0 \ \mod(\{f,g\},w).
$$

\noindent{\bf Proof.}
If
$\langle f,g\rangle_w$
and
$(f,g)_w$
are intersection compositions,
where
$w=\bar{f}b=a\bar{g}$,
then ({\bf XIII}) and ({\bf XVII}) yield
$$
\langle f,g\rangle_w=[fb]_{\bar{f}}-[ag]_{\bar{g}}
=fb+\sum\limits_{I_1}\alpha_ia_ifb_i-ag-\sum\limits_{I_2}\beta_ja_jgb_j,
$$
where
$a_i\bar{f}b_i,a_j\bar{g}b_j < \bar{f}b=a\bar{g}=w$.
Hence,
$$
\langle f,g\rangle_w-(f,g)_w\equiv_{ass}0 \mod(\{f,g\},w).
$$
In the case of inclusion compositions
we arrive at the same conclusion.
\hfill
$\blacksquare$

\begin{theorem}\label{thm4.9}
(
PBW Theorem in Shirshov's form
\cite{BoMa97,BoMa99}, see Theorem \ref{thm2.23})
A~nonempty set
$S\subset{Lie(X)}\subset{k\langle X\rangle}$
of monic Lie polynomials
is a~Gr\"{o}bner--Shirshov basis in
$Lie(X)$
if and only if~%
$S$
is a~Gr\"{o}bner--Shirshov basis in
$k\langle X\rangle$.
\end{theorem}

\noindent{\bf Proof.}
Observe that,
by definition,
for any
$f,g\in S$
the composition lies in
$Lie(X)$
if and only if
it lies
$k\langle X\rangle$.

Assume that~%
$S$
is a~GS basis in
$Lie(X)$.
Then we can express every composition
$\langle f,g\rangle_w$
as
$
\langle f,g\rangle_w=\sum\limits_{I_1}\alpha_i(a_is_ib_i),
$
where
$(a_is_ib_i)$
are normal
$S$-words and
$a_i\bar{s_i}b_i<w$.
By ({\bf XXI}),
we have
$
\langle f,g\rangle_w=\sum\limits_{I_2}\beta_jc_js_jd_j
$
with
$c_j\bar{s_j}d_j<w$.
Therefore,
({\bf XXII}) yields
$
(f,g)_w\equiv_{ass}0 \mod(S,w).
$
Thus,
$S$
is a~GS basis in
$k\langle X\rangle$.

Conversely,
assume that~%
$S$
is a~GS basis in
$k\langle X\rangle$.
Then the CD-lemma for associative algebras implies that
$\overline{\langle f,g\rangle_w}=a\overline{s}b<w$
for some
$a,b\in X^*$
and
$s\in S$.
Then
$h=\langle f,g\rangle_w-\alpha[asb]_{\bar s}\in Id_{ass}(S)$
is a~Lie polynomial
and
$\bar h<\overline{\langle f,g\rangle_w}$.
Induction on
$\overline{\langle f,g\rangle_w}$
yields
$
\langle
f,g\rangle_w\equiv_{Lie}0 \mod(S,w).
$
 \hfill
$\blacksquare$

\begin{theorem}\label{cdLie}
(The CD-lemma for Lie algebras over a~field)
Consider a~nonempty set
$S\subset{Lie(X)}\subset{k\langle X\rangle}$
of monic Lie polynomials
and denote by
$Id(S)$
the ideal of
$Lie(X)$
generated by~%
$S$.
The following statements are equivalent:
\begin{enumerate}[\rm(i)]
\item
The set~%
$S$
is a~Gr\"{o}bner--Shirshov basis in
$Lie(X)$.
\item
If
$f\in{Id(S)}$
then
${\bar{f}=a\bar{s}b}$
for some
$s\in{S}$
and
$a,b\in{X^*}$.
\item
The set
$$
Irr(S)=\{[u]\in NLSW(X) \ | \  u\neq{a\bar{s}b}, \
s\in{S},\ a,b\in{X^*}\}
$$
is a~linear basis for
$Lie(X|S)$.
\end{enumerate}
\end{theorem}

\noindent{\bf Proof.}
(i)$\Rightarrow$(ii).
Denote by
$Id_{ass}(S)$
and
$Id_{Lie}(S)$
the ideals of
$k\langle X\rangle$
and
$Lie(X)$
generated by~%
$S$ respectively.
Since
$Id_{Lie}(S)\subseteq Id_{ass}(S)$,
Theorem \ref{thm4.9} and the CD-lemma for associative algebras
imply the claim.

(ii)$\Rightarrow$(iii).
Suppose that
$\sum\alpha_i[u_i]=0$
in
$Lie(X|S)$
with
$[u_i]\in Irr(S)$
and
$u_1>u_2>\cdots$,
that is,
$\sum\alpha_i[u_i]\in{Id_{Lie}(S)}$.
Then all
$\alpha_i$
must vanish.
Otherwise we may assume that
$\alpha_1\neq0$.
Then
$\overline{\sum\alpha_i[u_i]}=u_1$
and (ii) implies that
$[u_1]\not\in{Irr(S)}$,
which is a~contradiction.
On the other hand,
by the next property ({\bf XXIII}),
$Irr(S)$
generates
$Lie(X|S)$
as a~linear space.

(iii)$\Rightarrow$(i).
This part follows from ({\bf XXIII}).
\hfill
$\blacksquare$

\ \

The next property is similar to Lemma \ref{L4.2}.

({\bf XXIII})
Given
$S\subset Lie(X)$,
we can express every
$f\in Lie(X)$
as
$$
f=\sum \alpha_i[u_i]+\sum\beta_j[a_js_jb_j]_{\bar{s_j}}
$$
with
$\alpha_i,\beta_j\in k$,
$[u_i]\in Irr(S)$
satisfying
$\overline{[u_i]}\leq \overline{f}$,
and
$[a_js_jb_j]_{\bar{s_j}}$
are special normal $S$-word satisfying
$\overline{[a_js_jb_j]_{\bar{s_j}}}\leq \overline{f}$.

\medskip

({\bf XXIV})
Given a~normal~%
$s$-word
$(asb)$,
take
$w=a\overline{s}b$.
Then
$(asb)\equiv[asb]_{\bar{s}}  \mod(s,w)$.
It follows that
$h\equiv_{Lie}0 \mod(S,w)$
is equivalent to
$h=\sum_{i}\alpha_i[a_is_ib_i]_{\bar{s_i}}$,
where
$[a_is_ib_i]_{\bar{s_i}}$
are special normal~%
$S$-words
with
$a_i\bar{s_i}b_i<w$.

\noindent{\bf Proof.}
Observe that
for every  monic Lie polynomial~%
$s$,
the set
$\{s\}$
is a~GS basis in
$Lie(X)$.
Then ({\bf XVIII}) and the CD-lemma for Lie algebras yield
$(asb)\equiv[asb]_{\bar{s}}  \mod(s,w)$.
\hfill
$\blacksquare$

\medskip

Summary of the proof of Theorem \ref{cdLie}.

Given two ALSWs
$u$~%
and~%
$v$,
define the ALSW-$\lcm(u,v)$
(or
$\lcm(u,v)$
for short)
as follows:
\begin{multline*}
w=\lcm(u,v) \in \{aucvb \text{ (an ALSW)},
a,b,c\in X^* \text{ (a trivial $\lcm$)};
\\
u=avb, a,b \in X^* \text{ (an inclusion $\lcm$)};
\\
ub=av, a,b\in X^*, deg(ub) <deg(u) + deg(v)
\text{ (an intersection $\lcm$)}\}.
\end{multline*}

Denote by
$[w]_{u,v}$
the Shirshov double special bracketing of~%
$w$
in the case that~%
$w$
is a~trivial
$\lcm(u,v)$,
by
$[w]_u$
and
$[w]_v$
the Shrishov special bracketings of~%
$w$
if~%
$w$
is an~inclusion or intersection
$\lcm$
respectively.
Then we can define a~general Lie composition for monic Lie polynomials
$f$
and
$g$
with
$\bar f=u$
and
$\bar g = v$
as
$$
(f,g)_w = [w]_{u,v}|_{[u]\mapsto f} -  [w]_{u,v}|_{[v]\mapsto g}
$$
if~%
$w$
is a~trivial
$\lcm(u,v)$
(it is
$0 \mod(\{f,g\}, w)$),
and
$$
(f,g)_w = [w]_u|_{[u]\mapsto f} -  [w]_v|_{[v]\mapsto g}
$$
if~%
$w$
is an~inclusion or intersection
$\lcm(u,v)$.

If
$S\subset Lie(X)\subset k\langle X\rangle$
is a~Lie GS basis
then~%
$S$
is an~associative GS basis.
This follows from property ({\bf XVII}) (iii)
and justifies the claim
(i)$\Rightarrow$(ii)
of Theorem \ref{cdLie}.

Shirshov's original proof of
(i)$\Rightarrow$(ii)
in Theorem \ref{cdLie},
(see \cite{Sh62b,Shirshov-Selected}),
rests on an~analogue of Lemma \ref{L4.1} for Lie algebras.

\begin{lemma}(\cite{Sh62b,Shirshov-Selected})
If
$(a_1s_1b_1), (a_2s_2b_2)$
are normal~%
$S$-words
with equal leading associative words,
$w=a_1\bar {s_1}b_1 = a_2\bar {s_2}b_2$,
then they are equal
$\mod(S,w)$,
that is,
$(a_1s_1b_1) - (a_2s_2b_2)\equiv 0 \mod(S,w).$
\end{lemma}

Outline of the proof.
We have
$w_1=cwd$
and
$w=\lcm(\bar {s_1},\bar {s_2})$.
Shirshov's (double) special bracketing lemma yields
$$
  [w_1]_w = [c[[w]d_1]d_2] = c[w]d + \sum \alpha_ia_i[w]b_i
$$
with
$a_iwb_i < w_1$.
The ALSW~%
$w$
includes
$u=\bar {s_1}$
and
$v=\bar
{s_2}$
as subwords,
and so there is a~bracketing
$\{w\}\in \{[w]_{u,v}, [w]_u, [w]_v\}$
such that
$$
[a_1s_1b_1] = [c \{w\}|_{[u]\mapsto s_1} d], \quad
[a_2s_2b_2] = [c \{w\}|_{[v]\mapsto s_2}d]
$$
are normal
$s_1$- and
$s_2$- words
with the same leading associative word
$w_1$.
Then
$$
[a_1s_1b_1] - [a_2s_2b_2] = [c(s_1,s_2)_wd] \equiv 0 \mod(S, w_1).
$$

Now it is enough to prove that
two normal Lie~%
$s$-words
with the same leading associative words,
say
$w_1$,
are equal
$\mod(s,w_1)$:
$$
f=(asb) - [asb] \equiv_{Lie} 0 \mod(s,w_1)
\quad\text{ provided that } \bar{f}<w_1.
$$
Since
$f\in Id_{ass}(s)$,
we have
$\bar f=c_1\bar{s}d_1$
by the CD-lemma for associative algebras
with one Lie polynomial relation
$s$.
Then
$f-\alpha [c_1sd_1]_{\bar s}$
is a~Lie polynomial with the leading associative word smaller than
$w_1$.
Induction on
$w_1$
finishes the proof.

\subsubsection{Gr\"{o}bner--Shirshov basis for
the Drinfeld--Kohno Lie algebra}

In this section we give a~GS basis for the Drinfeld--Kohno Lie algebra
$\textbf{L}_n$.
\begin{definition}
Fix an~integer
$n>2$.
The Drinfeld--Kohno Lie algebra
$\textbf{L}_{n}$
over
$\mathbb{Z}$
is defined by generators
$t_{ij}=t_{ji}$
for distinct indices
$1 \leq i,\ j \leq n-1$
satisfying the relations
$[t_{ij}t_{kl}]=0$
and
$[t_{ij}(t_{ik}+t_{jk})]= 0$
for distinct
$i$,
$j$,
$k$,
and~%
$l$.
\end{definition}

Therefore,
we have the presentation
$\textbf{L}_{n}=Lie_{\mathbb{Z}}(T| S)$,
where
$T=\{t_{ij}|\ 1 \leq i<j\leq n-1\}$
and~%
$S$
consists of the following relations:
\begin{eqnarray}\label{a}
&&[t_{ij}t_{kl}]=0 \quad\mbox{if } k<i<j,\ k<l,\ l\neq\ i,\ j;\\
&&\label{b} [t_{jk}t_{ij}]+[t_{ik}t_{ij}]=0 \quad\mbox{if } i<j<k;\\
&&\label{c} [t_{jk}t_{ik}]-[t_{ik}t_{ij}]=0 \quad\mbox{if} i<j<k.
\end{eqnarray}

Order
$T$
by setting
$t_{ij}<t_{kl}$
if either
$i<k$
or
$i=k$
and
$j<l$.
Let~%
$<$
be the deg-lex ordering on
$T^*$.

\begin{theorem}(\cite{ChenLiTang})\label{t1}
With
$S=\{(\ref{a}), (\ref{b}), (\ref{c})\}$
as before
and the deg-lex ordering~%
$<$
on
$T^*$,
the set~%
$S$
is a~Gr\"{o}bner--Shirshov basis of
$\textbf{L}_{n}$.
\end{theorem}

\begin{corollary}\label{cor1}
The Drinfeld--Kohno Lie algebra
$\textbf{L}_{n}$
is a~free~%
$\mathbb{Z}$-module
with
$\mathbb{Z}$-basis
$\cup_{i=1}^{n-2} NLSW(T_i)$,
where
$T_i=\{t_{ij}\ |\ i<j\leq\ n-1\}$
for
$i=1,\dots, n-2$.
\end{corollary}

\begin{corollary}(\cite{[Et]})
The Drinfeld--Kohno Lie algebra
$\textbf{L}_{n}$
is an~iterated semidirect product of free Lie algebras
$A_i$
generated by
$T_i=\{t_{ij}\ |\ i<j\leq\ n-1\}$,
for
$i=1,\dots, n-2$.
\end{corollary}

\subsubsection{Kukin's example of a~Lie algebra
with undecidable word problem}

A.~A.~Markov \cite{Markov},
E.~Post \cite{Post46},
A.~Turing \cite{Turing},
P.~S.~Novikov \cite{Novikov},
and W.~W.~Boone \cite{Boon}
constructed finitely presented semigroups and groups
with undecidable word problem.
For groups this also follows from
Higman's theorem \cite{Higman} asserting that
every recursively presented group
embeds into a~finitely presented group.
A~weak analogue of Higman's theorem for Lie algebras
was proved in \cite{Bo72},
which was enough for the existence of a~finitely presented Lie algebra
with undecidable word problem.
In this section
we give Kukin's construction \cite{Kukin}
of a~Lie algebra
$A_P$
for every semigroup~%
$P$
such that
if~%
$P$
has undecidable word problem
then so does
$A_P$.

Given a~semigroup
$P=\sgp\langle x,y|u_i=v_i, \ i\in I\rangle$,
consider the Lie algebra
$$
A_P=Lie(x, \hat{x}, y, \hat{y}, z|S)
$$
with~%
$S$
consisting of the relations
\begin{enumerate}[(1)]
\item
$[\hat{x}x]=0,\ [ \hat{x}y]=0,\ [ \hat{y}x]=0,\ [\hat{y}y]=0$;
\item
$[\hat{x}z]=-[zx],\  [\hat{y}z]=-[zy]$;
\item
$\lfloor zu_i\rfloor=\lfloor zv_i\rfloor,\ i\in  I$.
\end{enumerate}
Here,
$\lfloor zu\rfloor$
stands for the left normed bracketing.

Put
$\hat{x}>  \hat{y}> z>x>y$
and denote by~%
$>$
the deg-lex ordering on the set
$\{\hat{x},  \hat{y}, x,y,z\}^*$.
Denote by~%
$\rho$
the congruence on
$\{x,y\}^*$
generated by
$\{(u_i,v_i), \ i\in I\}$.
Put
\begin{enumerate}
\item[$(3')$]
$\lfloor zu\rfloor=\lfloor zv\rfloor,\ (u,v)\in \rho$
with
$u>v$.
\end{enumerate}

\begin{lemma}\label{L4.12}(\cite{ChenLiTang})
In this notation,
the set
$S_1=\{(1), (2), (3')\}$
is a~GS basis in
$Lie( \hat{x}, \hat{y}, x,y,z)$.
\end{lemma}

\noindent{\bf Proof:}
For every
$u\in \{x,y\}^*$,
we can show that
$\overline{\lfloor zu\rfloor}=zu$
by induction on
$|u|$.
All possible compositions in
$S_1$
are the intersection compositions of (2) and
$(3')$,
and the inclusion compositions of
$(3')$
and
$(3')$.

For
$(2)\wedge(3')$,
we take
$f=[\hat{x}z]+[zx]$
and
$g=\lfloor zu\rfloor-\lfloor zv\rfloor$.
Therefore,
$w=\hat{x}zu$
with
$(u,v)\in \rho$
and
$u>v$.
We have
\begin{eqnarray*}
&&\langle[ \hat{x}z]+[zx],\lfloor zu\rfloor-\lfloor
zv\rfloor\rangle_w=[fu]_{\bar f}-[\hat{x}g]_{\bar g}\\
&\equiv&\lfloor ([ \hat{x}z]+[zx])u\rfloor
-[ \hat{x}(\lfloor zu\rfloor-\lfloor zv\rfloor)]\\
&\equiv &\lfloor zxu\rfloor+\lfloor  \hat{x}zv\rfloor\equiv\lfloor
zxu\rfloor-\lfloor zxv\rfloor\equiv0 \mod(S_1,w).
\end{eqnarray*}
For
$(3')\wedge(3')$,
we use
$w=zu_1=zu_2e$,
where
$e\in \{x,y\}^*$
and
$(u_i,v_i)\in \rho$
with
$u_i>v_i$
for
$i=1,2$.
We have
\begin{eqnarray*}
&&\langle\lfloor zu_1\rfloor-\lfloor zv_1\rfloor,\lfloor
zu_2\rfloor-\lfloor zv_2\rfloor\rangle_w\equiv
(\lfloor zu_1\rfloor-\lfloor zv_1\rfloor)-\lfloor(\lfloor zu_2\rfloor-\lfloor zv_2\rfloor)e\rfloor\\
&\equiv &\lfloor\lfloor zv_2\rfloor e\rfloor-\lfloor
zv_1\rfloor\equiv\lfloor zv_2e\rfloor-\lfloor zv_1\rfloor\equiv0 \mod(S_1,w).
\end{eqnarray*}

Thus,
$S_1=\{(1), (2), (3')\}$
is a~GS basis in
$Lie(\hat{x}, \hat{y}, x,y,z)$.
\hfill
$\blacksquare$

\begin{corollary}(Kukin \cite{Kukin})\label{t4.12}
For
$u,v\in \{x,y\}^*$
we have
$$
u=v\ \mbox{ in the semigroup}\ P\Leftrightarrow\lfloor
zu\rfloor=\lfloor zv\rfloor\ \mbox{ in the Lie algebra }\ A_P.
$$
\end{corollary}

\noindent{\bf Proof:}
Assume that
$u=v$
in the semigroup~%
$P$.
Without loss of generality we may assume that
$u=au_1b$
and
$v=av_1b$
for some
$a,b\in \{x,y\}^*$
and
$(u_1,v_1)\in \rho$.
For every
$r\in \{x,y\}$
relations (1) yield
$[\hat{x}r]=0$;
consequently,
$\lfloor zxc\rfloor=\lfloor[z \hat{x}]c\rfloor=[\lfloor zc\rfloor
\hat{x}]$
and
$\lfloor zyc\rfloor=[\lfloor zc\rfloor\hat{y}]$
for every
$c\in\{x,y\}^*$.
This implies that in
$A_P$
we have
$$
\lfloor
zu\rfloor=\lfloor zau_1b\rfloor=\lfloor\lfloor zau_1\rfloor
b\rfloor=\lfloor\lfloor zu_1\widehat{\overleftarrow{a}}\rfloor
b\rfloor=\lfloor zu_1\widehat{\overleftarrow{a}}b\rfloor=\lfloor
zv_1\widehat{\overleftarrow{a}}b\rfloor=\lfloor
zav_1b\rfloor=\lfloor zv\rfloor,
$$
where for every
$x_{i_1}x_{i_2}\cdots x_{i_n}\in \{x,y\}^*$
we put
$$
\overleftarrow{x_{i_1}x_{i_2}\cdots x_{i_n}}
:=x_{i_n}x_{i_{n-1}}\cdots x_{i_1},
\quad
\widehat{x_{i_1}x_{i_2}\cdots x_{i_n}}
:=\widehat{x_{i_1}}\widehat{x_{i_{2}}}\cdots
\widehat{x_{i_n}}.
$$
Moreover,
$(3')$
holds in
$A_P$.

Suppose that
$\lfloor zu\rfloor=\lfloor zv\rfloor\ \mbox{ in the Lie
algebra }\ A_P$.
Then both
$\lfloor zu\rfloor$
and
$\lfloor zv\rfloor$
 have the same normal form in
$A_P$.
Since
$S_1$
is a~GS basis in
$A_P$,
we can reduce both
$\lfloor zu\rfloor$
and
$\lfloor zv\rfloor$
to the same normal form
$\lfloor zc\rfloor$
for some
$c\in \{x,y\}^*$
using only relations
$(3')$.
This implies that
$u=c=v$
in
$P$.
\hfill
$\blacksquare$

\medskip

By the corollary,
if the semigroup~%
$P$
has undecidable word problem
then so does the Lie algebra
$A_P$.

\subsection{Composition-Diamond lemma for Lie algebras
over commutative algebras}

For a~well-ordered set
$X=\{x_i |i\in I\}$,
consider the free Lie algebra
$Lie(X)\subset k\langle X \rangle$
with the Lie bracket
$[u,v]=uv-vu$.

Given a~well-ordered set
$Y=\{y_j|j\in J\}$,
the free commutative monoid
$[Y]$
generated by
$Y$
is a~linear basis of
$k[Y]$.
Regard
$$
Lie_{k[Y]}(X)\cong {k}[Y]\otimes Lie(X)
$$
as a~Lie subalgebra of the free associative algebra
$k[Y] \langle X \rangle\cong k[Y]\otimes k \langle X \rangle$
generated by~%
$X$
over the polynomial algebra
$k[Y]$,
equipped with the Lie bracket
$[u,v]=uv-vu$.
Then
$NLSW(X)$
constitutes a
$k[Y]$-basis
of
$Lie_{k[Y]}(X)$.
Put
$[Y]X^*=\{\beta t|\beta\in [Y], \ t\in X^*\}$.
For
$u=\beta t\in [Y]X^*$,
put
$u^X=t$
and
$u^Y=\beta$.

Denote the deg-lex orderings on
$[Y]$
and
$X^*$
by
$>_{_Y}$
and
$>_{_X}$.
Define an~ordering
$>$
on
$[Y]X^*$
as follows:
for
$u,v\in [Y]X^*$,
put
$$
u> v  \ \mbox{if}  \ (u^X>_{_X}v^X)  \ \mbox{or}\  (u^X=v^X  \
\mbox{and} \ u^Y>_{_Y}v^Y ).
$$

We can express every element
$f\in Lie_{k[Y]}(X)$
as
$f=\sum\alpha_i\beta_i[u_i]$,
where
$\alpha_i\in k$,
$\beta_i\in [Y]$,
and
$[u_i]\in NSLW(X).$

Then
$f=\sum\alpha_i\beta_i[u_i]=\sum g_j(Y)[u_j]$,
where
$g_j(Y)\in k[Y]$
are polynomials in the~%
$k$-algebra
$k[Y]\langle X\rangle$.
The leading word~%
$\bar{f}$
of~%
$f$
in
$k[Y]\langle X\rangle$
is of the form
$\beta_1 u_1$
with
$\beta_1\in [Y]$
and
$u_1\in ALSW(X)$.
The polynomial~%
$f$
is called monic
(or~%
$k$-monic)
if the coefficient of~%
$\bar f$
is equal to~1,
that is,
$\alpha_1 = 1$.
The notion of
$k[Y]$-monic
polynomials
is introduced similarly:
$\alpha_1 = 1$
and
$\beta_1 = 1$.

Recall that
every ALSW~%
$w$
admits a~unique bracketing
such that
$[w]$
is a NLSW.

Consider a~monic subset
$S\subset Lie_{k[Y]}(X)$.
Given a~non-associative word
$(u)$
on~%
$X$
with a~fixed occurrence of some
$x_i$
and
$s\in S$,
call
$(u)_{x_i\mapsto s}$
an~%
$S$-word.
Define
$|u|$
to be the~%
$s$-length
of
$(u)_{x_i\mapsto s}$.
Every~%
$S$-word
is of the form
$(asb)$
with
$a,b\in X^*$
and
$s\in S$.
If
$a\bar{s}^Xb\in ALSW(X)$
then we have the special bracketing
$[a\bar{s}^Xb]_{\bar{s}^X}$
of
$a\bar{s}^Xb$
relative to
$\bar{s}^X$.
Refer to
$[asb]_{\bar{s}}=[a\bar{s}^Xb]_{\bar{s}^X}|_{[\bar{s}^X]\mapsto{s}}$
as a~special normal
$s$-word
(or special normal
$S$-word).

An~%
$S$-word
$(u)=(asb)$
is a~normal~%
$s$-word,
denoted by
$\lfloor u\rfloor_{{s}}$,
whenever
$\overline{(asb)}^X=a\overline{s}^Xb$.
The following condition is sufficient.
\begin{enumerate}[(i)]
\item
The~%
$s$-length of
$(u)$
is~1,
that is,
$(u)=s$;
\item
if
$\lfloor u\rfloor_s$
is a~normal~%
$S$-word
of~%
$s$-length~%
$k$
and
$[v]\in NLSW(X)$
satisfies
$|v|=l$
then
$[v]\lfloor u\rfloor_s$
whenever
$v>\overline{\lfloor u\rfloor}_s^X$
and
$\lfloor u\rfloor_s[v]$
whenever
$v<\overline{\lfloor u\rfloor}_s^X$
are normal~%
$S$-words
of~%
$s$-length
$k+l$.
\end{enumerate}

Take two monic polynomials
$f$~%
and~%
$g$
in
$Lie_{k[Y]}(X)$
and put
$L=\lcm(\bar{f}^Y,\bar{g}^Y)$.

There are four kinds of compositions.
\begin{enumerate}
\item[$C_1$:]
Inclusion composition.
If
$\bar{f}^X=a\bar{g}^Xb$
for some
$a,b\in X^*$,
then
$$
C_1\langle f,g\rangle_w=
\frac{L}{\bar{f}^Y}f-\frac{L}{\bar{g}^Y}[agb]_{\bar{g}},
\quad\text{where } w=L\bar{f}^X=La\bar{g}^Xb.
$$
\item[$C_2$:]
Intersection composition.
If
$\bar{f}^X=aa_0$
and
$\bar{g}^X=a_0b$
with
$a,b,a_0\neq1$
then
$$
C_2\langle f,g\rangle_w=
\frac{L}{\bar{f}^Y}[fb]_{\bar{f}}-\frac{L}{\bar{g}^Y}[ag]_{\bar{g}},
\quad\text{where } w=L\bar{f}^Xb=La\bar{g}^X.
$$
\item[$C_3$:]
External composition.
If
$gcd(\bar{f}^Y,\bar{g}^Y)\neq1$
then for all
$a,b,c\in X^*$
satisfying
$$
w=La\bar{f}^Xb\bar{g}^Xc\in
T_A=\{\beta t|\beta\in [Y], \ t\in ALSW(X)\}
$$
we have
$$
C_3\langle f,g\rangle_w=
\frac{L}{\bar{f}^Y}[afb\bar{g}^Xc]_{\bar{f}}-
\frac{L}{\bar{g}^Y}[a\bar{f}^Xbgc]_{\bar{g}}.
$$
\item[$C_4$:]
Multiplication composition.
If
$\bar{f}^Y\neq1$
then for every special normal~%
$f$-word
$[afb]_{\bar{f}}$
with
$a,b\in X^*$
we have
$$
C_4\langle f\rangle_w=[a\bar{f}^Xb][afb]_{\bar{f}},
\quad\text{where } w=a\bar{f}^Xba\bar{f}b.
$$
\end{enumerate}

Given a
$k$-monic subset
$S\subset Lie_{k[Y]}(X)$
and
$w\in[Y]X^*$,
which is not necessarily in
$T_A$,
an~element
$h\in Lie_{k[Y]}(X)$
is called
trivial modulo
$(S,w)$
if
$h$
can be expressed as a~%
$k[Y]$-linear combination
of normal~%
$S$-words
with leading words smaller than~%
$w$.
The set~%
$S$
is a~Gr\"{o}bner--Shirshov basis in
$Lie_{k[Y]}(X)$
if all possible compositions 
in
$S$
are trivial.

\begin{theorem}\label{cdL2.1.4} (\cite{BoChChen-Liecomm},
the CD-lemma for Lie algebras over commutative algebras)
Consider a~nonempty set
$S\subset{Lie_{k[Y]}(X)}$
of monic polynomials
and denote by
$Id(S)$
the ideal of
$Lie_{k[Y]}(X)$
generated by~%
$S$.
The following statements are equivalent:
\begin{enumerate}[\rm(i)]
\item
The set~%
$S$
is a~Gr\"{o}bner--Shirshov basis in
$Lie_{k[Y]}(X)$.
\item
If
$f\in{Id(S)}$
then
${\bar{f}=a\bar{s}b\in T_A}$
for some
$s\in{S}$
and
$a,b\in{[Y]X^*}$.
\item
The set
$Irr(S)=\{[u] \ | \ [u]\in T_N, \ u\neq{a\bar{s}b},
\mbox{ for } s\in{S} \text{ and } a,b\in{[Y]X^*}\}$
is a~linear basis for
$Lie_{{\bf k}[Y]}(X|S)=(Lie_{k[Y]}(X))/Id(S)$.
\end{enumerate}
Here
$
T_A=\{\beta t\ |\ \beta\in [Y], \ t\in ALSW(X)\}
$
and
$
T_N=\{\beta[ t]\ |\ \beta\in [Y], \ [t]\in NLSW(X)\}.
$
\end{theorem}

\medskip

Outline of the proof.

Take
$u,v \in[Y] ALSW(X)$
and write
$u=u^Yu^X$
and
$v=v^Yv^X$.
Define the ALSW-$\lcm(u,v)$
(or
$\lcm(u,v)$
for short)
as
$w=w^Yw^X=\lcm(u^Y,v^Y)\lcm(u^X,v^X)$,
where
\begin{multline*}
\lcm(u^X,v^X) \in \{au^Xcv^Xb \ (an\ ALSW), a,b,c\in X^*;
\\
u^X=av^Xb, a,b \in X^* ;\ u^Xb=av^X, a,b\in X^*, deg(u^Xb) <deg(u^X) +
deg(v^X) \}.
\end{multline*}
Six
$\lcm(u,v)$
are possible:
\begin{itemize}
\item
($Y$-trivial, $X$-trivial)
(a~trivial
$\lcm(u,v))$;
\item
($Y$-trivial, $X$-inclusion);
\item
($Y$-trivial, $X$-intersection);
\item
($Y$-nontrivial, $X$-trivial);
\item
($Y$-nontrivial, $X$-inclusion);
\item
($Y$-nontrivial, $X$-intersection).
\end{itemize}
In accordance with
$\lcm(u,v)$,
six general compositions are possible.

Denote by
$[w^X]_{u^X,v^X}$
the Shirshov double special bracketing of
$w^X$
whenever
$w^X$
is a~%
$X$-trivial
$\lcm(u^X,v^X)$,
by
$[w^X]_{u^X}$
and
$[w^X]_{v^X}$
the Shirshov special bracketings of
$w^X$
whenever
$w^X$
is a lcm of~%
$X$-inclusion or~%
$X$-intersection
respectively.

Define general Lie compositions for~%
$k$-monic
Lie polynomials
$f$~%
and~%
$g$
with
$\bar f=u$
and
$\bar g = v$
as 
\begin{eqnarray*}
(f,g)_w &=&(\lcm(u^Y,v^Y)/u^Y) [w^X]_{u^X,v^X}|_{[u]\mapsto f}
-(\lcm(u^Y,v^Y)/v^Y) [w^X]_{u^X,v^X}|_{[v]\mapsto g},\\
(f,g)_w &=&(\lcm(u^Y,v^Y)/u^Y) [w^X]_u|_{[u]\mapsto f} -
(\lcm(u^Y,v^Y)/v^Y)[w^X]_v|_{[v]\mapsto g}.
\end{eqnarray*}

\begin{lemma}(\cite{BoChChen-Liecomm})
The general composition
$(f,g)_w$
of~%
$k$-monic
Lie polynomials
$f$~%
and~%
$g$
with
$\bar f=u$
and
$\bar g=v$,
where~%
$w$
is a~%
($Y$-trivial, $X$-trivial)
$\lcm(u,v)$,
is
$0\mod(\{f,g\},w).$
\end{lemma}

\noindent{\bf Proof:}
By ({\bf XIX}),
we have
\begin{eqnarray*}
(f,g)_w &=& v^Y[afb[v^X]d] -u^Y[a[u^X]bgd] = [afb[v]d] - [aubgd]\\
&=& [afb([v]-g)d] -[a([u]-f)bgd] \equiv 0 \mod(\{f,g\},w).
\end{eqnarray*}

The proof is complete.
\hfill
$\blacksquare$

\medskip

A~Lie GS basis
$S\subset Lie_{k[Y]}(X)\subset k[Y]\langle X\rangle$
need not be an~associative GS basis
because the PBW-theorem is not valid
for Lie algebras over a~commutative algebra
(Shirshov \cite{Sh53b}).
Therefore,
the argument for
$Lie_k(X)$
above
(see Section 4.2)
fails for
$Lie_{k[Y]}(X)$.

Moreover,
Shirshov's original proof of the CD-lemma fails
because the singleton
$\{s\}\in Lie_{k[Y]}(X)$
is not a~GS basis in general.
The reason is that
there exists a~nontrivial composition
$(s,s)_w$
of type
($Y$-nontrivial, $X$-trivial).

There is another obstacle.
For
$Lie_k(X)$,
every~%
$s$-word
is a~linear combination of normal~%
$s$-words.
For
$Lie_{k[Y]}(X)$
this is not the case.
Hence,
we must use a~multiplication composition
$[u^X]f$
such that
$\bar f=u=u^Yu^X$.

\begin{lemma}(\cite{BoChChen-Liecomm})\label{l4.17}
If every multiplication composition
$[{\bar s}^X]s$,
$s\in S$,
is trivial modulo
$(S, w=[u^X]u)$,
where
$u=\bar s$,
then every~%
$S$-word
is a~linear combination of normal~%
$S$-words.
\end{lemma}

In our paper with Yongshan Chen (\cite{BoChChen-Liecomm}),
we use the following definition of triviality of a~polynomial~%
$f$
modulo
$(S,w)$:
$$
f\equiv0 \mod(S,w) \Leftrightarrow
f =\sum\alpha _i e^Y_i[a_i^Xs_ib_i^X],
$$
where
$[a_i^X[\bar {s_i}^X]b_i^X]$
is the Shirshov special bracketing of the ALSW
$a_i^X\bar {s_i}^Xb_i^X$
with an~ALSW
$\bar {s_i}^X$.

The previous definition of triviality modulo
$(S,w)$
is equivalent to the usual definition
by Lemma \ref{l4.18},
which is key in the proof of
the CD-lemma for Lie algebras over a~commutative algebra.

\begin{lemma}(\cite{BoChChen-Liecomm})\label{l4.18}
Given a~monic set~%
$S$
with trivial multiplication compositions,
take a~normal
$s$-word
$(asb)$
and a~special normal
$s$-word
$[asb]$
with the same leading monomial
$w=a\bar sb$.
Then they are equal modulo
$(s,w)$.
\end{lemma}

Lemmas \ref{l4.17} and \ref{l4.18} imply

\begin{lemma}(\cite{BoChChen-Liecomm})
Given a~monic set~%
$S$
with trivial multiplication compositions,
every element of the ideal generated by~%
$S$
is a~linear combination of special normal~%
$S$-words.
\end{lemma}

On the other hand,
({\bf XVII}) and ({\bf XIX})
imply the following analogue of Lemma~\ref{L4.1} for
$Lie_{k[Y]}(X)$.

\begin{lemma}(\cite{BoChChen-Liecomm})
Given two~%
$k$-monic special normal~%
$S$-words
$e_1^Y[{a_1}^Xs_1{b_1}^X]$
and
$e_2^Y[{a_2}^Xs_2{b_2}^X]$
with the same leading associative word
$w_1$,
their difference is equal to
$[a(s_1,s_2)_wb]$,
where
$w=\lcm(\bar {s_1}, \bar {s_2})$,
$w_1=awb$,
and
$[a(s_1,s_2)_wb] = [w_1]_w|_{[w]\mapsto (s_1,s_2)_w}$.
Hence,
if~%
$S$
is a~GS basis
then the previous special normal~%
$S$-words
are equal modulo
$(S,w_1)$.
\end{lemma}

Now the claim
(i)$\Rightarrow$(ii)
of the CD-lemma for
$Lie_{k[Y]}(X)$
follows.

\medskip

For every Lie algebra
$\mathcal{L}=Lie_K(X|S)$
over the commutative algebra
$K=k[Y|R]$,
$$
U(\mathcal{L})=K\langle
X|S^{(-)}\rangle=k[Y]\langle X|S^{(-)}, RX\rangle,
$$
where
$S^{(-)}$
is just~%
$S$
with all commutators
$[uv]$
replaced with
$uv-vu$,
is the universal enveloping  associative algebra of~%
$\mathcal{L}$.

A~Lie algebra~%
$\mathcal{L}$
over a~commutative algebra~%
$K$
is called \textit{special}
whenever it embeds into its universal enveloping associative algebra.
Otherwise it is called \textit{non-special}.

A.~I.~Shirshov (1953) and P.~Cartier (1958)
gave classical examples
of non-special Lie algebras over commutative algebras over
$GF(2)$,
justified using ad hoc methods.
P.~M.~Cohn (1963)
suggested another non-special Lie algebra
over a~commutative algebra over a~field of positive characteristic.

\begin{example}(Shirshov (1953))
Take
$k=GF(2)$
and
$$
K=k[y_i,i=0,1,2,3|y_0y_i=y_i \ (i=0,1,2,3), \
y_iy_j=0 \ (i,j\neq0)].
$$
Consider
$\mathcal{L}=Lie_K(x_i, 1\leq i\leq 13|S_1, S_2)$,
where
\begin{eqnarray*}
S_1&=&\{[x_2x_1]=x_{11}, \ [x_3x_1]=x_{13}, \ [x_3x_2]=x_{12}, \\
&&\ \ [x_5x_3]=[x_6x_2]=[x_8x_1]=x_{10},  \ [x_ix_j]=0 \ (i>j)\};\\
S_2&=&\{y_0x_i=x_i  \ (i=1,2,\ldots,13),\\
&&\ \ y_1x_1=x_4, y_1x_2=x_5,  y_1x_3=x_6,  y_1x_{12}=x_{10},\\
&&\ \ y_2x_1=x_5,  y_2x_2=x_7,  y_2x_3=x_8,  y_2x_{13}=x_{10}, \\
&& \ \ y_3x_1=x_6,  y_3x_2=x_8,  y_3x_3=x_9,  y_3x_{11}=x_{10}, \\
&&\ \ y_1x_k=0 \ (k=4,5,\ldots,11,13), \\
&&\ \ y_2x_t=0 \ (t=4,5,\ldots,12),\\
&&\ \ y_3x_l=0 \ (l=4,5,\ldots,10,12,13)\}.
\end{eqnarray*}
Then
$\mathcal{L}=Lie_K(X|S_1, S_2)=Lie_{k[Y]}(X|S_1,S_2,RX)$
and
$$
S=S_1\cup S_2\cup RX\cup
\{ y_1x_2=y_2x_1, \ y_1x_3=y_3x_1, \ y_2x_3=y_3x_2\}
$$
is a~GS basis in
$Lie_{k[Y]}(X)$,
which implies that
$x_{10}$
belongs to the linear basis of~%
$\mathcal{L}$
by  Theorem \ref{cdL2.1.4},
that is,
$x_{10}\neq0$
in
$\mathcal{L}$.

On the other hand,
the universal enveloping algebra of~%
$\mathcal{L}$
has the presentation
$$
U_K(\mathcal{L})=K\langle X|S_1^{(-)},S_2\rangle\cong{\bf
k}[Y]\langle X|S_1^{(-)},S_2,RX\rangle.
$$
However,
the GS completion
(see Mikhalev and Zolotykh \cite{MZ})
of
$S_1^{(-)}\cup S_2\cup RX$
in
$k[Y]\langle X\rangle$
is
$$
S^C=S_1^{(-)}\cup S_2\cup RX\cup
\{y_1x_2=y_2x_1, \ y_1x_3=y_3x_1,  \ y_2x_3=y_3x_2, \  x_{10}=0\}.
$$

Thus,
$\mathcal{L}$
is not special.
\end{example}

\begin{example} (Cartier (1958))
Take
$k=GF(2)$
and
$$
K=k[y_1,y_2,y_3|y_i^2=0,\ i=1,2,3].
$$
Consider
$\mathcal{L}=Lie_{K}(x_{ij},1\leq i\leq j\leq3|S)$,
where
$$
S=\{[x_{ii}x_{jj}]=x_{ji} \ (i>j), [x_{ij}x_{kl}]=0, \
y_3x_{33}=y_2x_{22}+y_1x_{11}\}.
$$
Then~%
$\mathcal{L}$
is not special over~%
$K$.
\end{example}

\noindent{\bf Proof.}
The set
$S'=S\cup \{y_i^2x_{kl}=0\ (\forall i,k,l)\}\cup S_1$
is a~GS basis in
$Lie_{k[Y]}(X)$,
where
\begin{eqnarray*}
S_1&=&\{y_3x_{23}=y_1x_{12}, \  y_3x_{13}=y_2x_{12}, \ y_2x_{23}=y_1x_{13}, \ y_3y_2x_{22}=y_3y_1x_{11}, \\
&&\ \ y_3y_1x_{12}=0, \ y_3y_2x_{12}=0,  \ y_3y_2y_1x_{11}=0,   \
y_2y_1x_{13}=0\}.
\end{eqnarray*}
Then,
$y_2y_1x_{12}\in Irr(S')$
and so
$y_2y_1x_{12}\neq0$
in
$\mathcal{L}$.

However,
in
$$
U_K(\mathcal{L})=K\langle X|S^{(-)}\rangle\cong
{\bf k}[Y]\langle X|S^{(-)},y_i^2x_{kl}=0\ (\forall i,k,l)\rangle
$$
we have
\begin{eqnarray*}
0=y_3^2x_{33}^2=(y_2x_{22}+y_1x_{11})^2=y_2^2x_{22}^2+y_1^2x_{11}^2+y_2y_1[x_{22},x_{11}]
 = y_2y_1x_{12}.
\end{eqnarray*}
Thus,
$\mathcal{L}\not \hookrightarrow U_K(\mathcal{L})$.

\medskip

\noindent{\bf Conjecture} (Cohn \cite{Cohn63})
Take the algebra
$K=k[y_1,y_2,y_3|y_i^p=0, i=1,2,3]$
of truncated polynomials over a~field~%
$k$
of characteristic
$p>0$.
The algebra
$$
\mathcal{L}_p=Lie_{K}(x_1,x_2,x_3 \ | \ y_3x_3=y_2x_2+y_1x_1),
$$
called Cohn's Lie algebra,
is not special.

In
$U_K(\mathcal{L}_p)$
we have
$$
0= (y_3x_3)^p = (y_2x_2)^p + \Lambda_p(y_2x_2, y_1x_1) + (y_1x_1)^p=
\Lambda_p(y_2x_2, y_1x_1),
$$
where
$\Lambda_p$
is a~Jacobson--Zassenhaus Lie polynomial.
P.~M.~Cohn conjectured that
$\Lambda_p(y_2x_2, y_1x_1)\neq 0$
in
$\mathcal{L}_p$.
To prove this,
we must know a~GS basis of
$\mathcal{L}_p$
up to degree~%
$p$
in~%
$X$.
We found it for
$p=2,3,5$.
For example,
$\Lambda_2=[y_2x_2,y_1x_1]=y_2y_1[x_2x_1]$
and a~GS basis of
$\mathcal{L}_2$
up to degree~%
$2$
in~%
$X$
is
\begin{eqnarray*}
&&y_3x_3=y_2x_2+y_1x_1, \ y_i^2x_j=0\ (1\leq i,j\leq3), \  y_3y_2x_2=y_3y_1x_1, \ y_3y_2y_1x_1=0, \\
&&y_2[x_3x_2]=y_1[x_3x_1],  \ y_3y_1[x_2x_1]=0, \ y_2y_1[x_3x_1]=0.
\end{eqnarray*}
Therefore,
$y_2y_1[x_2x_1]\in Irr(S^C)$.

Similar though much longer computations show that
$\Lambda_3\neq0$
in
$\mathcal{L}_3$
and
$\Lambda_5\neq0$
in
$\mathcal{L}_5$.
Thus,
we have 

\begin{theorem}(\cite{BoChChen-Liecomm})\
Cohn's Lie algebras
$\mathcal{L}_2$,
$\mathcal{L}_3$,
and
$\mathcal{L}_5$
are non-special.
\end{theorem}

\begin{theorem}\label{t4.5}(\cite{BoChChen-Liecomm})
Given a~commutative~%
${k}$-algebra
$K={k}[Y|R]$,
if~%
$S$
is a~Gr\"{o}bner--Shirshov basis in
$Lie_{{k}[Y]}(X)$
such that
every
$s\in S$
is
${k}[Y]$-monic
then
$\mathcal{L}=Lie_{K}(X|S)$
is special.
\end{theorem}

\begin{corollary}\label{co4.6}(\cite{BoChChen-Liecomm})\
Every Lie~%
$K$-algebra
$L_K=Lie_K(X|f)$
with one monic defining relation
$f=0$
is special.
\end{corollary}

\begin{theorem}(\cite{BoChChen-Liecomm})
Suppose that~%
$S$
is a~finite homogeneous subset of
$Lie_{{k}}(X)$.
Then the word problem of
$Lie_{K}(X|S)$
is solvable
for every finitely generated commutative~%
${k}$-algebra~%
$K$.
\end{theorem}

\begin{theorem}(\cite{BoChChen-Liecomm})
Every finitely or countably generated Lie~%
$K$-algebra
embeds into a~two-generated Lie~%
$K$-algebra,
where~%
$K$
is an~arbitrary commutative~%
$k$-algebra.
\end{theorem}

\section{Gr\"obner-Shirshov bases for
$\Omega$-algebras and operads}

\subsection{CD-lemmas for $\Omega$-algebras}

Some new CD-lemmas for~%
$\Omega$-algebras
have appeared:
for associative conformal algebras \cite{BFK04}
and~%
$n$-conformal algebras \cite{BoChZhang-n-conf},
for the tensor product of free algebras \cite{BCC08},
for metabelian Lie algebras \cite{CC12},
for associative~%
$\Omega$-algebras \cite{BoChQiu-CD-Omega},
for color Lie superalgebras and  Lie~%
$p$-superalgebras \cite{Mikhalev89,Mikhalev92},
for Lie superalgebras \cite{Mikhalev96},
for associative differential algebras \cite{ChChLi-GSB-diff},
for associative Rota--Baxter algebras \cite{BCD08},
for~%
$L$-algebras \cite{BCH13},
for  dialgebras \cite{BCL08},
for pre-Lie algebras \cite{BoChLi-GSB-rightsym},
for semirings \cite{BCM13},
for commutative integro-differential algebras \cite{Guoli2013},
for difference-differential modules
and difference-differential dimension polynomials \cite{Zhou-Winkler},
for~%
$\lambda$-differential associative~%
$\Omega$-algebras \cite{ChQiu-Cd-diff},
for commutative associative Rota--Baxter algebras \cite{Qiu},
for algebras with differential type operators \cite{Guoli2012}.

V.~N.~Latyshev studied general versions of GS (or standard) bases
\cite{Latyshev98,Latyshev00}.

Let us state the CD-lemma for pre-Lie algebras,
see
\cite{BoChLi-GSB-rightsym}.

A~non-associative algebra~%
$A$
is called a~pre-Lie
(or a~right-symmetric)
algebra if~%
$A$
satisfies the identity
$(x,y,z)=(x,z,y)$
for the associator
$(x,y,z)=(xy)z-x(yz)$.
It is a~Lie admissible algebra in the sense that
$A^{(-)}=(A,[xy]=xy-yx)$
is a~Lie algebra.

Take a~well-ordered set
$X=\{x_i|i\in I \}$.
Order
$X^{**}$
by induction on the lengths of the words
$(u)$
and
$(v)$:
\begin{enumerate}[(i)]
\item
When
$|((u)(v))|=2$
put
$(u)=x_i > (v)=x_j$
if and only if
$i>j$.
\item
When
$|((u)(v))|>2$
put
$(u)>(v)$
if and only if
one of the following holds:
\begin{enumerate}[(a)]
\item
$|(u)|>|(v)|$;
\item
if
$|(u)|=|(v)|$
with
$(u)=((u_1)(u_2))$
and
$(v)=((v_1)(v_2))$
then
$(u_1)>(v_1)$
or
$(u_1)=(v_1)$
and
$(u_2)>(v_2)$.
\end{enumerate}
\end{enumerate}

We now quote the definition of good words (see \cite{Se94})
by induction on length:
\begin{enumerate}[(1)]
\item
$x$
is a~good word for any
$x\in X$;
\item
a~non-associative word
$((v)(w))$
is called a~good word if
\begin{enumerate}[(a)]
\item
both
$(v)$
and
$(w)$
are good words and
\item
if
$(v)=((v_1)(v_2))$
then
$(v_2)\leq(w)$.
\end{enumerate}
\end{enumerate}

Denote
$(u)$
by
$[u]$
whenever
$(u)$
is a~good word.

Denote by~%
$W$
the set of all good words in the alphabet~%
$X$
and by
$RS\langle X\rangle$
the free right-symmetric algebra over a~field~%
$k$
generated by~%
$X$.
Then~%
$W$
forms a~linear basis of
$RS\langle X\rangle$,
see \cite{Se94}.
D.~Kozybaev, L.~Makar-Limanov, and U.~Umirbaev \cite{KMLU08}
proved that
the deg-lex ordering on~%
$W$
is monomial.

Given a~set
$S\subset RS\langle X\rangle$
of monic polynomials
and
$s\in S$,
an~%
$S$-word
$(u)_s$
is called a~normal~%
$S$-word
whenever
$(u)_{\bar s }=(a\bar s b)$
is a~good word.

Take
$f,g\in S$,
$[w]\in W$,
and
$a,b\in X^{*}$.
Then there are two kinds of compositions.
\begin{enumerate}[(i)]
\item
If
$\bar{f}=[a\bar{g}b]$
then
$
(f,g)_{\bar{f}}=f-[agb]
$
is called the inclusion composition.

\item
If
$(\bar{f}[w])$
is not good
then
$
f\cdot [w]
$
is called the right multiplication composition.
\end{enumerate}

\begin{theorem}\label{t2.4.1} (\cite{BoChLi-GSB-rightsym},
the CD-lemma for pre-Lie algebras)
Consider a~nonempty set
$S\subset RS\langle X\rangle$
of monic polynomials
and the ordering~%
$<$
defined above.
The following statements are equivalent:
\begin{enumerate}[\rm(i)]
\item
The set~%
$S$
is a~Gr\"{o}bner-Shirshov basis in
$RS\langle X\rangle$.

\item
If
$f\in Id(S)$
then
$\bar f =[a\bar s b]$
for some
$s\in S$
and
$a,b\in X^*$,
where
$[as b]$
is a~normal~%
$S$-word.

\item
The set
$Irr(S)=\{[u]\in W |[u]\ne [a\bar s b]\ a,b\in X^*,\ s\in S \mbox{ and }
[as b] \mbox{ is a~normal } S\mbox{-word}\}$
is a~linear basis of the algebra
$RS\langle X | S\rangle=RS\langle X\rangle/Id(S)$.
\end{enumerate}
\end{theorem}

As an~application,
we have a~GS basis for
the universal enveloping pre-Lie algebra of a~Lie algebra.

\begin{theorem}\label{t2.4.2}(\cite{BoChLi-GSB-rightsym})\
Consider a~Lie algebra
$({\cal{L}},[\ ])$
with a~well-ordered linear basis
$X=\{e_i|\ i\in I\}$.
Write
$
[e_ie_j]=\sum\limits_{m}\alpha_{ij}^me_m
$
with
$\alpha_{ij}^m\in k$.
Denote
$\sum\limits_{m}\alpha_{ij}^me_m$
by
$\{e_ie_j\}$.
Denote by
$$
U({\cal{L}})=RS\langle \{e_i\}_I| \ e_ie_j-e_je_i=\{e_ie_j\}, \ i,j
\in I\rangle
$$
the universal enveloping pre-Lie algebra of
${\cal{L}}$.
The set
\begin{eqnarray*}
S&=&\{f_{ij}=e_ie_j-e_je_i-\{e_ie_j\},\ i,j \in I \ \mbox{ and } \
i>j \}
\end{eqnarray*}
is a~Gr\"{o}bner-Shirshov basis in
$RS\langle X\rangle$.
\end{theorem}

Theorems \ref{t2.4.1} and \ref{t2.4.2}
directly imply the following PBW theorem
for Lie algebras and pre-Lie algebras.

\begin{corollary}
(D.~Segal \cite{Se94})
A~Lie algebra
${\cal{L}}$
embeds into its universal enveloping pre-Lie algebra
$U({\cal{L}})$
as a~subalgebra of
$U({\cal{L}})^{(-)}$.
\end{corollary}

\medskip

Recently the CD-lemmas mentioned above and other combinatorial methods
yielded many applications:
for groups of Novikov--Boone type \cite{Kalorkoti06,Kalorkoti09,Kalorkoti11}
(see also \cite{bokut66,bokut67,ChChLuo,Kalorkoti82},
for Coxeter groups \cite{bs,Lee},
for center-by-metabelian Lie algebras \cite{Umirbaev84},
for free metanilpotent Lie algebras,
Lie algebras and associative algebras
\cite{GUmirbaev99,MSUmirbaev04,Umirbaev93a,Umirbaev93b},
for Poisson algebras \cite{MakarLimanovU11},
for quantum Lie algebras and related problems
\cite{Kharchenko02,Kharchenko10},
for PBW-bases
\cite{Kharchenko99,Kharchenko08,MakarLimanov94},
for extensions of groups and associative algebras \cite{Chen08a,Chen09a},
for (color) Lie ($p$)-superalgebras
\cite{Bakhturin92,BokutKLM,Chibrikov06a,Chibrikov06b,GerdtKornyak95,GerdtKornyak96,GerdtKornyak97,MikhalevZ95,MikhalevZ97,MikhalevZ98},
for Hecke algebras and Specht modules \cite{KangLL02},
for representations of Ariki--Koike algebras \cite{KangLL04},
for the linear algebraic approach to GS bases \cite{kl2},
for HNN groups \cite{ChenZhong08},
for certain one-relator groups \cite{ChenZhong11},
for embeddings of algebras \cite{BokutChenMo-embed,ChenMo-AMS},
for free partially commutative Lie algebras
\cite{ChMo-Partialcomm,poroshenko},
for quantum groups of type
$D_n$,
$E_6$,
and
$G_2$
\cite{Abdu3,Abdu1,Abdu2,Abdu4},
for calculations of homogeneous GS bases \cite{Scala-Leva},
for Picard groups,
Weyl groups,
and Bruck--Reilly extensions of semigroups
\cite{Ates1,Karpuz,Ates3,Ates4,Ates2}.

\medskip

\subsection{CD-lemma for operads}

Following Dotsenko and Khoroshkin (\cite{DK10}, Proposition 3), linear bases for a
symmetric operad and  a shuffle operad are the same provided both of them
are defined by the same generators and defining relations. It means that we need CD-lemma
for shuffle operads only (and we define a GS basis for a symmetric
operad as a GS basis of the corresponding shuffle operad).

We express the elements of the free shuffle operad using planar trees.

Put
$\mathscr{V}=\bigcup_{n=1}^{\infty}\mathscr{V}_{n}$,
where
$\mathscr{V}_{n}=\{\delta_i^{(n)}|i\in I_n\}$
is the set of~%
$n$-ary
operations.

Call a planar tree with~%
$n$
leaves \textit{decorated}
whenever the leaves are labeled by
$[n]=\{1,2,3,\ldots,n\}$
for
$n\in \mathbb{N}$
and every vertex is labeled by an~element of~%
$\mathscr{V}$.

A~decorated tree is called a~tree monomial
whenever for each vertex the minimal value
on the leaves of the left subtree
is always less than that of the right subtree.

Denote by
$\mathscr{F}_\mathscr{V}(n)$
the set of all tree monomials with~%
$n$
leaves
and put
$T=\cup_{n\geq1}\mathscr{F}_\mathscr{V}(n)$.
Given
$\alpha=\alpha(x_1,\dots, x_n)\in \mathscr{F}_\mathscr{V}(n)$
and
$\beta\in \mathscr{F}_\mathscr{V}(m)$,
define the shuffle composition
$\alpha \circ_{i,\sigma} \beta$
as 
$$
\alpha(x_1,
\ldots,x_{i-1},\beta(x_i,x_{\sigma(i+1)},\ldots,x_{\sigma(i+m-1)}),x_{\sigma(i+m)},\ldots,x_{\sigma(m+n-1)}),
$$
which lies in
$\mathscr{F}_\mathscr{V}(n+m-1)$,
where
$1\leq i\leq n$
and the bijection
$$
\sigma:\{i+1,\ldots,m+n-1\}\rightarrow\{i+1,\ldots,m+n-1\}
$$
is an
$(m-1,n-i)$-shuffle,
that is,
\begin{eqnarray*}
&&\sigma(i+1)<\sigma(i+2)<\dots<\sigma(i+m-1),\\
&&\sigma(i+m)<\sigma(i+m+1)<\dots<\sigma(n+m-1).
\end{eqnarray*}

The set~%
$T$
is freely generated by~%
$\mathscr{V}$
with the shuffle composition.

Denote by
$\mathscr{F}_\mathscr{V}=kT$
the~%
$k$-linear space
spanned by~%
$T$.
This space
with the shuffle compositions
$\circ_{i,\sigma}$
is called the free shuffle operad.

Take a~homogeneous subset~%
$S$
of
$\mathscr{F}_\mathscr{V}$.
For
$s\in S$,
define an~%
$S$-word
$u|_s$
as before.

A~well ordering~%
$>$
on~%
$T$
is called monomial
(admissible)
whenever
$$
\alpha>\beta\Rightarrow u|_\alpha>u|_\beta \  \mbox{ for any } u\in
T.
$$

Assume that
$T$
is equipped with a~monomial ordering.
Then each~%
$S$-word
is a~normal~%
$S$-word.

For example,
the following ordering~%
$>$
on~%
$T$
is monomial,
see Proposition~5 of~\cite{DK10}.

Every
$\alpha=\alpha(x_1,\dots, x_n)\in \mathscr{F}_\mathscr{V}(n)$
has a~unique expression
$$
\alpha=(\path(1),\dots, \path(n),[i_1\dots i_n]),
$$
where
$\path(r)\in\mathscr{V}^*$
for
$1\leq r\leq n$
is the unique path from the root to the leaf~%
$r$
and the permutation
$[i_1\dots i_n]$
lists the labels of the leaves of the underlying tree
in the order determined by the planar structure,
from left to right.
In this case define
$$
\wt(\alpha)=(n,\path(1),\dots, \path(n),[i_1\dots i_n]).
$$

Assume that
$\mathscr{V}$
is a~well-ordered set
and use the deg-lex ordering on~%
$\mathscr{V}^*$.
Take the order on the permutations in reverse lexicographic order:
$i>j$
if and only if~%
$i$
is less than~%
$j$
as numbers.

Now,
given
$\alpha,\beta\in T$,
define
$$
\alpha>\beta \Leftrightarrow \wt(\alpha)>\wt(\beta)
\quad\text{lexicographically}.
$$

An~element of
$\mathscr{F}_\mathscr{V}$
is called \textit{homogeneous}
whenever all tree monomials
occurring in this element with nonzero coefficients
have the same arity degree
(but not necessarily the same operation degree).

For two tree monomials
$\alpha$~%
and~%
$\beta$,
say that~%
$\alpha$
is divisible by~%
$\beta$
whenever there exists a~subtree of the underlying tree of~%
$\alpha$
for which the corresponding tree monomial
$\alpha'$
is equal to~%
$\alpha$.

A~tree monomial~%
$\gamma$
is called a~common multiple of two tree monomials
$\alpha$~%
and~%
$\beta$
whenever it is divisible by both
$\alpha$~%
and~%
$\beta$.
A~common multiple~%
$\gamma$
of two tree monomials
$\alpha$~%
and~%
$\beta$
is called a~least common multiple
and denoted by
$\gamma=\lcm(\alpha,\beta)$
whenever
$|\alpha|+|\beta|>|\gamma|$,
where
$|\delta|=n$
for
$\delta\in \mathscr{F}_\mathscr{V}(n)$.

Take two monic homogeneous elements
$f$~%
and~%
$g$
of
$\mathscr{F}_\mathscr{V}$.
If
$\bar{f}$~%
and~%
$\bar{g}$
have a~least common multiple~%
$w$
then
$(f, g)_w=w_{_{\bar f\mapsto f}}-w_{_{\bar g\mapsto g}}$.

\begin{theorem}(\cite{DK10},
the CD-lemma for shuffle operads)
In the above notation,
consider a~nonempty set
$S\subset\mathscr{F}_\mathscr{V}$
of monic homogeneous elements
and a~monomial ordering~%
$<$
on~%
$T$.
The following statements are equivalent:
\begin{enumerate}[\rm(i)]
\item
The set~%
$S$
is a~Gr\"{o}bner-Shirshov basis in
$\mathscr{F}_\mathscr{V}$.
\item
If
$f\in Id(S)$
then
$\bar{f}=u|_{\bar{s}}$
for some~%
$S$-word~%
$u|_s$.
\item
The set
$Irr(S)=\{u\in T| u \neq
v|_{\bar{s}}
 \mbox{ for all } \ S\mbox{-word }\ v|_s\}$
is a~%
$k$-linear basis of
$\mathscr{F}_\mathscr{V}/Id(S)$.
\end{enumerate}
\end{theorem}

As applications,
the authors of \cite{DK10}
calculate Gr\"{o}bner-Shirshov bases for some well-known operads:
the operad Lie of Lie algebras,
the operad As of associative algebras,
and the operad PreLie of pre-Lie algebras.


\begin{thebibliography}{99}

\bibitem{Adjan}Adjan, S.I.: Algorithmic undecidability of certain decision problems of group theory.
 Trudy Moscow Mat. Ob. {\bf6}, 231 - 298 (1957).

\bibitem{AAJEfim} Alahmadi, A., Alsulami, H., Jain, S.K., Zelmanov,
E.: Leavitt path algebras of finite Gelfand-Kirillov dimension. J.
Algebra Appl. {\bf11}, No. 6, 1250225-1--1250225-6 (2012)

\bibitem{AAJEfim2013} Alahmadi, A., Alsulami, H., Jain, S.K., Zelmanov,
E.: Structure of Leavitt path algebras of polynomial growth.
www.pnas.org/cgi/doi/10.1073/pnas.1311216110

\bibitem{Artamonov} Artamonov, V.A.: Clones of multilinear operations
and multioperator for algebras. Uspekhi  Mat. Nauk. {\bf 24} (145),
47-59 (1969)

\bibitem{Artin26} Artin, E.: Theory der Z\"opf. Abh. Math. Sem. Hamburg Univ. {\bf 4},   47-72 (1926)

\bibitem{Artin47}Artin, E.:  Theory of braids. Ann. Math. {\bf48}, 101-126 (1947)


\bibitem{Ates1} Ates, F., Karpuz,  E., Kocapinar,  C.,  Cevik,  A.S.:
Gr\"{o}bner--Shirshov  bases of some monoids, Discrete Mathematics,
{\bf311}(12), 1064-1071  (2011)


\bibitem{Bakhturin11} Bahturin, Yu.A., Olshanskii, A.: Filtrations and Distortion in
Infinite-Dimensional Algebras.
 J. Algebra  {\bf327}, 251-291 (2011)

\bibitem{Bakhturin92}
Bahturin, Yu.A.,  Mikhalev,  A.A., Petrogradskij, V.M., Zajtsev,
M.V.:  Infinite dimensional Lie superalgebras.  [B] De Gruyter
Expositions in Mathematics. 7. Berlin etc.: W. de Gruyter. x, 250 p.
(1992)




\bibitem{Belyaev} Belyaev, V.Ya.: Subrings of finitely presented associative
rings. Algebra  Logika  {\bf17},  627-638 (1978)


\bibitem{Be78}Bergman, G.M.: The diamond lemma for ring theory.  Adv.
Math. {\bf29},  178-218 (1978)

\bibitem{Berstel-J-P}
 Berstel,  J.D.,  Perrin, D.: The origins of combinatorics on words.
  Eur. J. Comb. {\bf28},    996-1022  (2007)

\bibitem{Birman-Ko-Lee}Birman, J., Ko, K.H., Lee, S.J.: A~new approach to the word and
conjugacy problems for the braid groups.  Adv. Math. {\bf139},
 322-353 (1998)


\bibitem{Bjoner}Bjoner, A., Brenti, F.: Combinatorics of Coxeter
Groups. Graduate Texts in Math. 231, Springer 2005


\bibitem{bokut63} Bokut, L.A.:
A~base of free polynilpotent Lie algebras. Algebra   Logika {\bf2},
   13-19 (1963)



\bibitem{bokut66} Bokut, L.A.: On one property of the Boone group.  Algebra  Logika
{\bf5},  5-23 (1966)

\bibitem{bokut67} Bokut, L.A.: On the Novikov groups.  Algebra
Logika {\bf6},   25-38  (1967)

\bibitem{Bok1968} Bokut, L.A.:  Degrees of unsolvability of the conjugacy problem for finitely presented groups, Algebra Logika {\bf5, 6}, 4-70, 4-52 (1968)

\bibitem{bokut69i-iii} Bokut, L.A.: Groups of fractions for the multiplicative
semigroups of certain rings I-III. Sibirsk. Mat. Zh. {\bf10},
246-286, 744-799, 800-819 (1969)

\bibitem{bokut69iv} Bokut, L.A.: On the Malcev problem. Sibirsk. Mat. Zh.
{\bf10}, 965-1005 (1969)

\bibitem{Bo72} Bokut, L.A.: Insolvability of the word problem for Lie algebras, and
subalgebras of finitely presented Lie algebras.  Izvestija AN~USSR
(mathem.)   {\bf36},  1173-1219  (1972)

\bibitem{Bo76} Bokut, L.A.: Imbeddings into simple associative
algebras. Algebra   Logika  {\bf15},  117-142 (1976)




\bibitem{Bo08} Bokut, L.A.:  Gr\"obner--Shirshov bases for braid groups in Artin-Garside
generators. J. Symbolic Computation  {\bf43},  397-405 (2008)

\bibitem{Bo09} Bokut, L.A.: Gr\"obner--Shirshov bases for the braid group in the
Birman--Ko--Lee generators. J. Algebra {\bf321},  361-379  (2009)


\bibitem{Bo-Ch-Shum} Bokut, L.A., Chainikov, V.V., Shum, K.P.:
Markov and Artin normal form theorem for braid groups. Commun.
Algebra  {\bf35}, 2105-2115  (2007)


\bibitem{BokutChainikov08} Bokut, L.A.,   Chainikov, V.V.: Gr\"obner--Shirshov bases of Adjan extension
of the Novikov group. Discrete Mathematics {\bf308}, 4916-4930
(2008).




\bibitem{BC07} Bokut, L.A.,  Chen, Y.Q.:  Gr\"{o}bner--Shirshov bases for Lie algebras: after
A.I. Shirshov.  Southeast Asian Bull. Math. {\bf31}, 1057-1076
(2007)

\bibitem{survey08} Bokut, L.A.,  Chen, Y.Q.:   Gr\"{o}bner--Shirshov
bases: some new results, Advance in Algebra and Combinatorics.
Proceedings of the Second International Congress in Algebra and
Combinatorics, Eds. K. P. Shum, E. Zelmanov, Jiping Zhang, Li
Shangzhi, World Scientific, 2008, 35-56.


\bibitem{ChLiJing-plactic} Bokut, L.A., Chen, Y.Q.,   Chen, W.P.,  Li, J.:
New approaches to plactic monoid via Gr\"{o}bner--Shirshov bases. arxiv.org/abs/1106.4753

\bibitem{BCC08} Bokut, L.A.,  Chen, Y.Q.,  Chen, Y.S.: Composition-Diamond lemma
for tensor product of free algebras. J. Algebra {\bf323}, 2520-2537
(2010)


\bibitem{BoChChen-Liecomm} Bokut, L.A.,  Chen, Y.Q.,  Chen, Y.S.:  Gr\"{o}bner--Shirshov
bases for Lie algebras over a~commutative algebra. J. Algebra
{\bf337}, 82-102 (2011)




\bibitem{BCD08} Bokut, L.A.,  Chen, Y.Q.,   Deng, X.M.:
Gr\"{o}bner--Shirshov bases for Rota--Baxter algebras.  Siberian Math.
J. {\bf51},    978-988  (2010)


\bibitem{BCH13} Bokut, L.A.,  Chen, Y.Q.,  Huang, J.P.: Gr\"{o}bner--Shirshov bases for
L-algebras. Internat. J. Algebra Comput. {\bf23}, 547-571 (2013)


\bibitem{BoChLi09} Bokut, L.A.,  Chen, Y.Q.,  Li, Y.: Anti-commutative Gr\"obner--Shirshov
basis of a~free Lie algebra. Science in China Series A: Mathematics
 {\bf52},  244-253  (2009)

\bibitem{BoChLi-GSB-rightsym}  Bokut, L.A.,  Chen, Y.Q.,  Li, Y.: Gr\"{o}bner--Shirshov bases for
Vinberg--Koszul--Gerstenhaber right-symmetric algebras. J. Math. Sci.
{\bf166},  603-612 (2010)



\bibitem{BoChLi-CD-category}  Bokut, L.A.,  Chen, Y.Q.,  Li, Y.: Gr\"{o}bner--Shirshov bases for
categories. Nankai Series in Pure, Applied Mathematics and
Theoretical Physical, Operads and Universal Algebra  {\bf9}, 1-23
(2012)


\bibitem{BoChLi10} Bokut, L.A.,  Chen, Y.Q.,   Li, Y.: Lyndon--Shirshov words and anti-commutative
algebras. J. Algebra {\bf378}, 173-183 (2013)




\bibitem{BCL08} Bokut, L.A.,  Chen, Y.Q.,  Liu, C.H.:
Gr\"{o}bner--Shirshov bases for dialgebras. Internat. J. Algebra
Comput.  {\bf20},  391-415  (2010)


\bibitem{BokutChenMo-embed} Bokut, L.A.,  Chen, Y.Q.,   Mo, Q.H.: Gr\"{o}bner--Shirshov bases and
embeddings of algebras. Internat. J. Algebra Comput.  {\bf20},
 875-900 (2010)

\bibitem{BCM13} Bokut, L.A.,  Chen, Y.Q.,  Mo, Q.H.: Gr\"{o}bner--Shirshov bases for
semirings.   J. Algebra {\bf385}, 47-63 (2013)



\bibitem{BoChQiu-CD-Omega} Bokut, L.A.,  Chen, Y.Q.,   Qiu, J.J.: Gr\"{o}bner--Shirshov bases
for associative algebras with multiple operations and free
Rota--Baxter algebras. J. Pure Applied Algebra {\bf214}, 89-100
(2010)



\bibitem{BCS}  Bokut, L.A.,  Chen, Y.Q.,   Shum, K.P.: Some new results on Gr\"{o}bner--Shirshov
bases.
Proceedings of International Conference on Algebra 2010, Advances in
Algebraic Structures, 2012, pp.53-102.

\bibitem{BoChZhang-n-conf} Bokut, L.A.,  Chen, Y.Q.,  Zhang, G.L.:
Composition-Diamond lemma for associative n-conformal algebras.
arXiv:0903.0892

\bibitem{BoChZhao-inverse-sg} Bokut, L.A.,  Chen, Y.Q.,   Zhao, X.G.: Gr\"{o}bner--Shirshov beses
for free inverse semigroups. Internat. J. Algebra Comput. {\bf19},
 129-143 (2009)



\bibitem{BFK04}Bokut, L.A., Fong, Y.,  Ke, W.-F.: Composition Diamond lemma for
 associative conformal algebras.  J.  Algebra {\bf272}, 739-774
 (2004)


\bibitem{BFKK00}Bokut, L.A.,  Fong, Y.,  Ke, W.-F., Kolesnikov, P.S.: Gr\"obner and
Gr\"obner--Shirshov bases in algebra and conformal algebras.
Fundamental and Applied Mathematics  {\bf6}, 669-706   (2000)

\bibitem{BFKS02}Bokut, L.A., Fong, Y., Ke, W.-F., Shiao, L.-S.:
Gr\"{o}bner--Shirshov bases for the braid semigroup.  Shum, K.P.
(ed.) et al.  Advances in algebra. Proceedings of the ICM satellite
conference in algebra and related topics, Hong Kong, China, August
14-17  (2002)


\bibitem{BokutKLM}Bokut, L.A., Kang, S.-J., Lee, K.-H.,  Malcolmson, P.:
Gr\"obner--Shirshov bases for Lie superalgebras and their universal
enveloping algebras, J. Algebra {\bf217}, 461-495 (1999)

\bibitem{BokutKlein96}Bokut, L.A., Klein, A.A.: Serre relations and Gr\"obner--Shirshov
bases for simple Lie algebras I, II. Internat. J. Algebra Comput.
{\bf6},  389-400, 401-412  (1996)


\bibitem{BokutKlein98}Bokut, L.A., Klein, A.A.: Gr\"obner--Shirshov bases for exceptional
Lie algebras  I. J. Pure Applied Algebra {\bf133}, 51-57  (1998)

\bibitem{BokutKlein99}Bokut, L.A., Klein, A.A.: Gr\"obner--Shirshov bases for exceptional
Lie algebras E6, E7, E8. Algebra and Combinatorics (Hong Kong),
37-46, Springer, Singapore  1999




\bibitem{BK03}Bokut, L.A.,  Kolesnikov, P.S.: Gr\"obner--Shirshov bases: from their
incipiency to the present.  J. Math. Sci.  {\bf116}, 2894-2916
(2003)


\bibitem{BK05}Bokut, L.A., Kolesnikov, P.S.: Gr\"obner--Shirshov bases, conformal
algebras and pseudo-algebras.  J. Math. Sci.  {\bf131}, 5962-6003
(2005)



\bibitem{BokutKukin}Bokut, L.A., Kukin, G.P.:  Algorithmic and Combinatorial
Algebra.  Mathematics  and its Applications,  Kluwer Academic
Publishers Group, Dordrecht  (1994)



\bibitem{BoMa}Bokut, L.A., Malcolmson, P.:  Gr\"{o}bner--Shirshov bases
for quantum enveloping algebras.  Isr. J.   Math. {\bf96}, 97-113
(1996)

\bibitem{BoMa97}Bokut, L.A.,  Malcolmson, P.: Gr\"obner--Shirshov bases for Lie and associative
algebras.
Collection of Abstracts, ICAC,97, Hong Kong, 1997, 139-142.

\bibitem{BoMa99}Bokut, L.A.,  Malcolmson, P.: Gr\"obner--Shirshov bases for relations
of a~Lie algebra and its enveloping algebra.   Shum, Kar-Ping (ed.)
et al., Algebras and combinatorics. Papers from the international
congress, ICAC'97, Hong Kong, August 1997. Singapore: Springer.
47-54 (1999)

\bibitem{bs}Bokut, L.A., Shiao, L.-S.: Gr\"{o}bner--Shirshov bases for Coxeter
groups. Commun. Algebra   {\bf29},  4305-4319   (2001)

\bibitem{BokutShum}Bokut, L.A., Shum, K.P.: Relative Gr\"{o}bner--Shirshov bases for algebras and
groups. St. Petersbg. Math. J. {\bf19},  867-881 (2008)



\bibitem{Boon} Boone, W.W.: The word problem. Ann. Math.
{\bf70},   207-265  (1959)

\bibitem{Borcherds86}Borcherds, R.E.: Vertex algebras, Kac--Moody algebras, and the
monster.  Proc. Natl. Acad. Sci. USA {\bf84}, 3068-3071 (1986)


\bibitem{Borcherds88}Borcherds, R.E.: Generalized Kac--Moody algebras.   J.
Algebra {\bf115}, No.2, 501-512 (1988)

\bibitem{Borcherds90}Borcherds, R.E.:
The monster Lie algebra. Adv. Math. {\bf83}, No.1, 30-47 (1990)


\bibitem{BrieskornSaito}Brieskorn, E., Saito, K.: Artin-Gruppen und Coxeter-Gruppen.
 Invent. Math. {\bf17}, 245-271 (1972)

\bibitem{bu65} Buchberger, B.: An~algorithm for finding a~basis for the
residue class ring of a~zero-dimensional polynomial ideal. Ph.D.
thesis, University of Innsbruck, Austria  (1965)


\bibitem{bu70} Buchberger, B.: An~algorithmical criteria for the
solvability of algebraic systems of equations.   Aequationes Math.
{\bf4}, 374-383  (1970)


\bibitem{Bu87}Buchberger, B.:  History and basic feature of the critical-pair/completion procedure. J. Symbolic
Computation  {\bf3},   3-38 (1987)

\bibitem{BuCL} Buchberger, B., Collins, G.E., Loos, R., Albrecht, R.: Computer
algebra, symbolic and algebraic computation, Computing Supplementum,
Vol.4, New York: Springer-Verlag  (1982)

\bibitem{Portugal}Cain, A.J., Gray, R., Malheiro, A.: Finite  Gr\"{o}bner--Shirshov bases for Plactic algebras and
biautomatic structures for Plactic monoids,  arXiv:1205.4885v2

\bibitem{Cartier} Cartier, P.: Remarques sur le th$\acute{\mbox e}$or$\grave{\mbox e}$me de
Birkhoff-Witt, Annali della Scuola Norm. Sup. di Pisa s$\acute{\mbox
e}$rie III vol XII(1958), 1-4.


\bibitem{Jc01} Cassaigne, J., Espie, M., Krob, D., Novelli, J.C., Hivert, F.: The Chinese
Monoid. Internat. J. Algebra Comput.  {\bf11}, 301-334  (2001)



\bibitem{CFL} Chen, K.-T., Fox, R., Lyndon, R.: Free differential calculus IV: The quotient group of the
lower central series.  Ann. Math. {\bf68}, 81-95 (1958)

\bibitem{Chen08a}Chen, Y.Q.:  Gr\"{o}bner--Shirshov basis for Schreier extensions of
groups. Commun. Algebra {\bf36}, 1609-1625  (2008)

\bibitem{Chen09a}Chen, Y.Q.:  Gr\"{o}bner--Shirshov basis for extensions of
algebras. Algebra Colloq. {\bf16}  283-292 (2009)

\bibitem{CC12}Chen, Y.S., Chen, Y.Q.:  Gr\"{o}bner--Shirshov bases for matabelian Lie
algebras.   J. Algebra {\bf358},  143-161 (2012)



\bibitem{ChChLi-GSB-diff}Chen, Y.Q., Chen, Y.S.,  Li, Y.: Composition-Diamond lemma
for differential algebras. The Arabian Journal for Science and
Engineering  {\bf34},  135-145 (2009)





\bibitem{ChChLuo}Chen, Y.Q.,  Chen, W.S.,  Luo, R.I.: Word problem for Novikov's
and Boone's group via Gr\"{o}bner--Shirshov bases. Southeast Asian
Bull. Math.  {\bf32},  863-877   (2008)

\bibitem{CCZ}Chen, Y.Q.,   Chen, Y.S.,  Zhong, C.Y.: Composition-Diamond
lemma for modules.  Czechoslovak Math. J. {\bf60},  59-76 (2010)

\bibitem{CL}Chen,Y.Q., Li,Y.:
Some remarks for the Akivis algebras and the Pre-Lie algebras.
  Czechoslovak Math. J. {\bf61}(136), 707-720 (2011)


\bibitem{ChenLiTang}Chen, Y.Q.,   Li, Y.,   Tang, Q.Y.: Gr\"{o}bner--Shirshov bases for
some Lie algebras. Siberian Math.
J. to appear. arXiv:1305.4546

\bibitem{CLZ} Chen,Y.Q., Li, J., Zeng, M.J.:
Composition-Diamond Lemma for Non-associative Algebras over a
Commutative Algebra. Southeast Asian Bull. Math. {\bf34}, 629-638
(2010).


\bibitem{ChMo-GSB-Braid-Artin}Chen, Y.Q.,   Mo, Q.H.: Artin-Markov normal form for braid
group. Southeast Asian Bull. Math. {\bf33},   403-419 (2009)



\bibitem{ChenMo-AMS} Chen, Y.Q.,   Mo, Q.H.: Embedding dendriform algebra into its
universal enveloping Rota--Baxter algebra.  Proc. Am. Math. Soc.
{\bf139},   4207-4216 (2011)



\bibitem{ChMo-Partialcomm}Chen, Y.Q.,  Mo, Q.H.: Gr\"{o}bner--Shirshov bases for free partially commutative
Lie algebras. Commun. Algebra {\bf 41}, 3753-3761 (2013)

\bibitem{Chen-Qiu}Chen, Y.Q., Qiu, J.J.: Gr\"{o}bner--Shirshov basis for the Chinese
monoid. J. Algebra Appl. {\bf7}, 623-628 (2008)


\bibitem{ChShaoShum}Chen, Y.Q.,  Shao, H.S., Shum, K.P.: On Rosso-Yamane theorem
on PBW basis of
$U_q(A_N)$.
CUBO a~Mathematical Journal {\bf10},
171-194 (2008)


\bibitem{ChenZhong08}Chen, Y.Q.,  Zhong, C.Y.: Gr\"{o}bner--Shirshov basis for HNN extensions of groups
and for the alternative group.  Commun. Algebra {\bf36}, 94-103
(2008)

\bibitem{ChenZhong11}Chen, Y.Q., Zhong, C.Y.:  Gr\"{o}bner--Shirshov basis for some one-relator
groups.  Algebra Colloq.  {\bf19}, 99-116  (2011)

\bibitem{ChZhong-Braid}Chen, Y.Q.,  Zhong, C.Y.:  Gr\"{o}bner--Shirshov bases for braid
groups in Adjan--Thurston generators.  Algebra Colloq.   {\bf20},
 309-318 (2013)




\bibitem{Chibrikov04}Chibrikov, E.S.: On free conformal Lie algebras.
  Vestn. Novosib. Gos. Univ.  Ser. Mat.
Mekh. Inform. {\bf4}, No. 1, 65-83 (2004)


\bibitem{Chibrikov06a}Chibrikov, E.S.: A~right normed basis for free Lie algebras and
Lyndon--Shirshov words.  J. Algebra {\bf302},   593-612 (2006)

\bibitem{Chibrikov06b}Chibrikov, E.S.: The right-normed basis for a~free
Lie superalgebra and Lyndon--Shirshov words.  Algebra Logika {\bf45},
No. 4, 458-483 (2006)

\bibitem{Chibrikov11}Chibrikov, E.S.:  On free Sabinin algebras.   Commun.
Algebra {\bf39},  4014-4035 (2011)

\bibitem{Chibrikov12}Chibrikov, E.S.: On some embedding of Lie algebras.   J.
Algebra Appl. {\bf11}, No. 1,  12 p. (2012)






\bibitem{Cohn63} Cohn, P.M.: A~remark on the Birkhoff-Witt theorem.
 Journal London Math. Soc.  {\bf38},  197-203 (1963)



\bibitem{Cohnbook65} Cohn, P.M.:   Universal algebra.   Harper's
Series in Modern Mathematics. New York-Evanston-London: Harper and
Row, Publishers xv, 333 p. (1965). Second edition: Reidel,
Dordrecht-Boston (1981)


\bibitem{Collins}Collins, D.J.: Representation of Turing reducibility by word and
conjugacy problems in finitely presented groups.  Acta Math.
{\bf128}, 73-90 (1972)


\bibitem{DK10} Dotsenko, V.,  Khoroshkin, A.: Gr\"{o}bner bases for
operads. Duke Mathematical Journal  {\bf153}, 363-396 (2010)



\bibitem{eisenbud} Eisenbud, D., Peeva, I., Sturmfels, B.:
Non-commutative Gr\"{o}bner bases for commutative algebras. Proc.
Am. Math. Soc.  {\bf126},  687-691  (1998)


\bibitem{[Et]} Etingof, P., Henriques, A.,  Kamnitzer, J., Rains,
E.M.: The cohomology ring of the real locus of the moduli space of
stable curves of genus 0 with marked points. Ann. Math. {\bf171},
731-777 (2010)

\bibitem{FFG} Farkas, D.R., Feustel, C., Green, E.I.: Synergy in the
theories of Gr\"obner bases and path algebras. Canad. J. Math. {\bf45}, 727-739
(1993)


\bibitem{Guoli2013} Gao, X., Guo, L.,  Zheng, S.H.: Constrction of free
commutative integro-differential algebras by the method of
Gr\"obner--Shirshov bases. J. Algebra Applications, to appear.

\bibitem{Garside1969} Garside, A.F.:  The braid group and other groups. Q. J. Math.
Oxford {\bf20}, 235-254  (1969)

\bibitem{Gelfand} Gelfand, S.I.,  Manin, Y.I.: Homological Algebra.
Springer-Verlag  (1999)


\bibitem{GerdtKornyak95}Gerdt, V.P., Kornyak, V.V.: Lie algebras and superalgebras defined by
a~finite number of relations: computer analysis. J. Nonlinear Math.
Phys. {\bf2},  No.3-4, 367-373 (1995)

\bibitem{GerdtKornyak96}Gerdt, V.P., Robuk, V.N., Sever'yanov, V.M.: The construction of
finitely represented Lie algebras.  Comput. Math. Math. Phys.
{\bf36},  No.11, 1493-1505 (1996)

\bibitem{GerdtKornyak97}Gerdt, V.P., Kornyak, V.V.: Program for constructing a~complete
system of relations, basis elements, and commutator table for
finitely presented Lie algebras and superalgebras.   Program.
Comput. Softw. {\bf23},  No.3, 164-172 (1997)

\bibitem{Golod}Golod, E.S.: Standard bases and homology. in: Algebra: Some
current trends. Lect. Notes Math. 1352, 88-95 (1988)


\bibitem{Green} Green, D.J.:  Gr\"{o}bner Bases and the Computation
of Group Cohomology, Springer-Verlag Berlin Heideberg 2003


\bibitem{GreenJA} Green, J.A.:  Hall algebras, hereditary algebras and guantum algebras. Invent. Math.
{\bf120}, 361-377 (1985)

\bibitem{Guoli2012} Guo, L., Sit, W., Zhang, R.: Differential type operators and
Gr\"{o}bner--Shirshov bases. J. Symbolic Comput. {\bf52}, 97-123
(2013)


\bibitem{GUmirbaev99} Gupta, C.K., Umirbaev, U.U.: The occurrence problem for free
metanilpotent Lie algebras.  Commun. Algebra {\bf27},  5857-5876
(1999)


\bibitem{MH} Hall, M.: A~basis for free Lie rings and higher
commutators in free groups.  Proc. Am. Math. Soc.   {\bf3}, 575-581
(1950)


\bibitem{PH} Hall, P.: A~contribution to the theory of groups of prime
power order.  Proc. London Math. Soc. Ser. {\bf36}, 29-95 (1933)



\bibitem{Higman} Higman, G.: Subgroups of finitely presented groups.
Proc. Royal Soc. London (Series A)  {\bf262}, 455-475 (1961)

\bibitem{Jones} Jones, V.F.R.: Hecke algebra representations of braid
groups and link polynimials. Ann. Math. {\bf128}, 335-388  (1987)

\bibitem{Kac} Kac, G.: Infinite dimensional Lie algebras. Cambridge
University Press, Cambridge (1990)

\bibitem{Kalorkoti82}Kalorkoti, K.: Decision problems in group theory.  Proc.
Lond. Math. Soc., III. Ser. 44, 312-332 (1982)


\bibitem{Kalorkoti06}Kalorkoti, K.: Turing degrees and the word and conjugacy problems for
finitely presented groups. Southeast Asian Bull. Math. {\bf30},
855-887 (2006)


\bibitem{Kalorkoti09}Kalorkoti, K.: A~finitely presented group with almost solvable
conjugacy problem.  Asian-Eur. J. Math. {\bf2},  611-635 (2009)


\bibitem{Kalorkoti11}Kalorkoti, K.: Sufficiency conditions for Bokut' normal forms.
 Commun. Algebra {\bf39},   2862-2873 (2011)

\bibitem{Kandri-Rody}Kandri-Rody, A., Weispfenning, V.: Non-commutative Gr\"{o}bner bases in
algebras of solvable type A. J. Symbolic Comput. {\bf9},
1-26 (1990)


\bibitem{kl1} Kang, S.-J., Lee,  K.-H.:  Gr\"{o}bner--Shirshov
bases for representation theory.  J. Korean Math. Soc.  {\bf37},
55-72 (2000)


\bibitem{KL} Kang, S.-J., Lee,  K.-H.: Gr\"obner--Shirshov bases for irreducible
$sl_{n+1}$-modules.  J. Algebra {\bf232},  1-20  (2000)



\bibitem{KangLL02}Kang, S.-J.; Lee, I.-S., Lee, K.-H., Oh, H.: Hecke
algebras, Specht modules and Gr\"{o}bner--Shirshov bases.  J. Algebra
{\bf252},  258-292 (2002)

\bibitem{KangLL04}Kang, S.-J.; Lee, I.-S., Lee, K.-H., Oh, H.:  Representations of Ariki-Koike algebras and
Gr\"{o}bner--Shirshov bases.   Proc. Lond. Math. Soc., III. Ser.
{\bf89},   54-70 (2004)

\bibitem{kl2} Kang, S.-J., Lee, K.-H.: Linear algebraic approach to Gr\"{o}bner--Shirshov basis
theory.  J. Algebra  {\bf313}, 988-1004 (2007)

\bibitem{Karpuz} Karpuz, E.G.: Complete Rewriting System for the Chinese Monoid.
Applied Mathematical Sciences {\bf4},  1081-1087 (2010)



\bibitem{Ates3} Karpuz, E.G., Cevik, A.S.: Gr\"{o}bner--Shirshov  bases for
extended modular, extended Hecke, and Picard groups, Marhematical
Notes, Volume: 92,  636-642 (2012)

\bibitem{Ates4} Karpuz, E.G., Ates, F., Cevik, A.S.: Gr\"{o}bner--Shirshov
bases of some Weyl groups, Rocky Mont. J. Math., to appear.

\bibitem{Kharchenko99} Kharchenko, V.K.: A~quantum analog of the Poincar¨¦-Birkhoff-Witt
theorem.   Algebra Logika {\bf38}, No. 4, 476-507 (1999)

\bibitem{Kharchenko02} Kharchenko, V.K.: A~combinatorial approach to the quantification of
Lie algebras. Pac. J. Math. {\bf203},   191-233 (2002)



\bibitem{Kharchenko05} Kharchenko, V.K.: Braided version of Shirshov-Witt theorem.
J. Algebra {\bf294},   196-225 (2005)

\bibitem{Kharchenko08}  Kharchenko, V.K.: PBW-bases of coideal subalgebras and a~freeness
theorem.  Trans. Am. Math. Soc. {\bf360}, 5121-5143 (2008)

\bibitem{Kharchenko10} Kharchenko, V.K.: Triangular decomposition of right coideal
subalgebras.  J. Algebra {\bf324},  3048-3089 (2010)


\bibitem{Sapir94}  Kharlampovich, O.G., Sapir, M.V.: Algorithmic problems in
varieties. Inter. J. Algebra Compt.  {\bf5},  379-602 (1995)

\bibitem{Knuth} Knuth D.E.:
Permutations, matrices, and generalized Young tableaux.
Pacific J.\ Math.\ {\bf 34}, 709-727 (1970).

\bibitem{Knuth-Bendix} Knuth, D.E., Bendix, P.B.: Simple word problems in universal
algebras. In: Computational problems in abstract algebra. Pergamon
Press, Oxford
 263-297 (1970)

\bibitem{Ates2}Kocapinar,  C.,  Karpuz, E., Ates, F.,  Cevik, A.S.:
Gr\"{o}bner--Shirshov  bases of generalized Bruck-Reilly
*-extension, Algebra Colloquium, {\bf19},  813-820 (2012)


\bibitem{Kolchin1973} Kolchin, E.R.: Differential algebras and algebraic groups. Academic
Press, New York  (1973)



\bibitem{KMLU08}Kozybaev, D.,   Makar-Limanov, L., Umirbaev, U.:
 The Freiheitssatz and autoumorphisms of free
right-symmetric algebras.  Asian-European Journal of Mathematics
{\bf1}, 243-254 (2008)

\bibitem{Kukin} Kukin, G.P.: On the word problem for Lie algebras.
Sibirsk. Mat. Zh. {\bf18}, 1194-1197  (1977)

\bibitem{Kurosh} Kurosh, A.G.: Nonassociative free algebras and free
products of algebras. Mat. Sb.  {\bf20}, 239-262 (1947)


\bibitem{Kurosh69} Kurosh, A.G.: Multioperator ringpond algebras.
Uspekhi Mat. Nauk    {\bf24} (145), 3-15 (1969)

\bibitem{Scala-Leva}La Scala, R., Levandovskyy, V.: Letterplace ideals and
non-commutative Gr\"{o}bner bases.   J. Symb. Comput. {\bf44},
1374-1393 (2009)

\bibitem{Lambek}Lambek, J., Deductive system and categories II: standard constructions and closed categories, in:
Lecture Notes in Math. {\bf86}, Springer 1969


\bibitem{M.L} Lascoux, A., Leclerc, B., Thibon, J.Y.: The plactic monoid. in:
Algebraic Combinatorics on Words. Cambridge University Press,
Cambridge (2002)

\bibitem{Latyshev98}Latyshev, V.N.: General version of standard bases in linear
structures.   Bahturin, Yu. (ed.), Algebra.   Berlin: Walter de
Gruyter. 215-226 (2000)

\bibitem{Latyshev00}Latyshev, V.N.: An~improved version of standard
bases. Krob, Daniel (ed.) et al., Formal power series and algebraic
combinatorics.  Berlin: Springer. 496-505 (2000)


\bibitem{Lazard60} Lazard, M.: Groupes, anneaux de Lie et probl$\grave{e}$me
de Burnside. Istituto Matematico dell' Universit$\grave{a}$
di Roma
(1960)

\bibitem{Lee}Lee, D.V.: Gr\"{o}bner--Shirshov bases and normal forms for the Coxeter
groups
$E_6$
and
$E_7$.
 Shum, K.P. (ed.) et al., Advances in
algebra and combinatorics.  Hackensack, NJ: World Scientific.
243-255 (2008)



\bibitem{Lothaire-etc}Lothaire, M.:    Combinatorics on words.   Addison-Wesley Publishing Company,
xix, 238 p. (1983). Second edition: Cambridge University Press,
Cambridge (1977)


\bibitem{Lothaire02}  Lothaire, M.: Algebraic Combinatorics on Words. Cambridge University Press,
Cambridge (2002)

\bibitem{Lusztig}  Lusztig, G.: Canonical bases arising from quantized
enveloping algebras. J. Amer. Math. Soc. {\bf3}, 447-498(1990)


\bibitem{Lusztig2003}  Lusztig, G.:  Hecke algebras with unequal parameters. CRM monograph series vol.18, American Mathematical Society (2003)

\bibitem{Lyndon}  Lyndon, R.C.: On Burnside's problem I.  Trans. Amer. Math. Soc. {\bf77}, 202-215 (1954)

\bibitem{Sa}Maclane, S.: Homology. Springer-Verlag  (1963)



\bibitem{MakarLimanov94}Makar-Limanov, L.: A~version of the Poincar$\acute{e}$-Birkhoff-Witt theorem.
 Bull. Lond. Math. Soc. {\bf26}, No.3, 273-276 (1994)



\bibitem{MakarLimanovU11}Makar-Limanov, L., Umirbaev, U.U.: The Freiheitssatz for
Poisson algebras.  J. Algebra {\bf328}, No. 1, 495-503 (2011)




\bibitem{Markov1945} Markov, A.A.:  An~introduction to the algebraical theory of
braids.  Proceedings  of the Steklov Mat. Ins. RAS, vol. 16, 1945.

\bibitem{Markov}Markov,  A.A.: Impossibility of some algorithms in the
theory of some associative system. Dokl. Akad. Nauk SSSR {\bf55},
 587-590 (1947)

\bibitem{May}May, P.: The geometry of Iterated Loop Space.   Lecture Notes in Math. {\bf271}, Springer
1972

\bibitem{Michel75} Michel, J.: Bases des alg$\grave{e}$bres de Lie et
s$\acute{e}$rie de Hausdorff,  Semin. P. Dubreil, 27e annee 1973/74,
Algebre, Fasc. 1, Expose 6, 9 p. (1975)



\bibitem{Michel76}Michel, J.:  Calculs dans les alg$\grave{e}$bres
de Lie libres: la s$\acute{e}$rie de Hausdorff et le
probl$\grave{e}$me de Burnside,   Ast$\acute{e}$risque 38/39,
139-148 (1976)

\bibitem{Mikhalev89}Mikhalev, A.A.: A~composition lemma and the equality problem for
color Lie superalgebras.   Mosc. Univ. Math. Bull. {\bf44}, No.5,
87-90 (1989)

\bibitem{Mikhalev92}Mikhalev, A.A.: The composition lemma for color Lie superalgebras and
for Lie
$p$-superalgebras.   Algebra, Proc. Int. Conf. Memory A.I.
Mal'cev, Eds. L.A. Bokut, Yu.L. Ershov, A.I. Kostrikin,
Novosibirsk/USSR 1989, Contemp. Math. 131, Pt. 2, 91-104 (1992)

\bibitem{Mikhalev96} Mikhalev, A.A.: Shirshov composition techniques in
Lie superalgebras (noncommutative Gr\"obner bases). J. Math. Sci.,
New York {\bf80}, No.5, 2153-2160 (1996)


\bibitem{MSUmirbaev04} Mikhalev, A.A., Shpilrain, V., Umirbaev, U.U.: On
isomorphism of Lie algebras with one defining relation.  Int. J.
Algebra Comput. {\bf14}, No. 3, 389-393 (2004)


\bibitem{MikhalevZ95}Mikhalev, A.A., Zolotykh, A.A.: Combinatorial aspects of Lie
superalgebras.  Boca Raton, FL: CRC Press. viii, 260 p. (1995)





\bibitem{MZ}Mikhalev, A.A., Zolotykh, A.A.: Standard Gr\"obner--Shirshov bases of free algebras over
rings, I. Free associative algebras. Internat. J. Algebra Comput.
{\bf8}, 689-726 (1998)

\bibitem{Mora} Mora, F.:
Gr\"obner bases for non-commutative polynomial rings. in: Algebraic
algorithms and error-correcting codes.  Lect. Notes Comput. Sci.
{\bf229}, 353-362 (1986)



\bibitem{Newman} Newman, M.H.A.: On theories with a~combinatorial definition of
`equivalence'.  Ann. Math. {\bf43},  223-243 (1942)


\bibitem{Novikov} Novikov, P.S.: On algorithmic undecidability of the
word problem in the theory of groups.  Trudy Mat. Inst. Steklov.
{\bf44}, 1-144 (1955)


\bibitem{Abdu3} Obul, A., Yunus, G.: Gr\"{o}bner--Shirshov basis of quantum group
of type
$E_6$.
J. Algebra {\bf346}, 248-265 (2011)



\bibitem{Odesskii} Odesskii, A.: Introduction to the theory of
elliptic algebras. data.imf.au.dk/conferences/FMOA05/

\bibitem{Schein}Poliakova, O., Schein, B.M.: A~new construction for
free inverse semigroups.  J. Algebra  {\bf288}, 20-58 (2005)

\bibitem{Polishchuk} Polishchuk, A., Positselski, L.: Quadratic
algebras, AMS  (2005)


\bibitem{poroshenkoB}Poroshenko, E.N.: Gr\"{o}bner--Shirshov bases
for Kac--Moody algebras
$A^{(1)}_n$
and
$B^{(1)}_n$.
Krob, Daniel
et al. Formal power series and algebraic combinatorics. Berlin
Springer 552-563 (2000)


\bibitem{poroshenkoC-D}Poroshenko, E.N.: Gr\"{o}bner--Shirshov bases
for Kac--Moody algebras of types
$C^{(1)}_n$
and
$D^{(1)}_n$.
Vestn.
Novosib. Gos. Univ. Ser. Mat. Mekh. Inform. {\bf2}, 58-70 (2002)

\bibitem{poroshenkoA}Poroshenko, E.N.: Gr\"{o}bner--Shirshov bases
for Kac--Moody algebras of type
$A^{(1)}_n$.
Commun. Algebra {\bf30},
2617-2637 (2002)



\bibitem{poroshenko}Poroshenko, E.N.: Bases for partially commutative Lie
algebras. Algebra  Logika {\bf50}, 405-417 (2011)


\bibitem{Post46} Post, E.: A~variant of a~recursively unsolvable
problem. Bull. Amer. Math. Soc. {\bf52}, 264-268  (1946)

\bibitem{Post47} Post, E.:  Recursive unsolvability of a~problem of Thue. J. Symbolic Logic {\bf1}, 1-11 (1947)

\bibitem{Priddy} Priddy, S.B.: Koszul resolutions. Trans. Amer. Math. Soc. {\bf152},
39-60 (1970)

\bibitem{ChQiu-Cd-diff}Qiu, J.J., Chen,  Y.Q: Composition-Diamond lemma for
$\lambda$-differential associative algebras with multiple operators.
J. Algebra Appl.   {\bf9}, 223-239 (2010)

\bibitem{Qiu} Qiu, J.J.: Gr\"{o}bner--Shirshov bases for commutative
algebras with multiple operators  and free commutative Rota--Baxter
algebras.  Asian-European Jour. Math.  to appear.

\bibitem{Rabin}Rabin,M.: Recursice unsolvability of group thepretic
problems. Ann. of Math. {\bf67}, (1) 172-194, (1958).

\bibitem{Razmyslov}Razmyslov, Yu.P.:  Identities of algebras and their representations.
Providence, RI: AMS xiii, 318p. (1994)



\bibitem{Abdu1}Ren, Y.H.,  Obul, A.: Gr\"{o}bner--Shirshov basis of quantum group
of type
$G_2$.
 Commun. Algebra {\bf39}, 1510-1518 (2011)



\bibitem{Reutenauer}Reutenauer,  C.: Free Lie algebras. Oxford University Press, New York  (1993)

\bibitem{Ringel90} Ringel, C.M.:  Hall algebras and quantum groups.
Invent. Math.  {\bf101}, 583-592 (1990)

\bibitem{Ringel96}Ringel, C.M.:  PBW-bases of quantum groups.
J. reine angew. Math.  {\bf170}, 51-88 (1996)

\bibitem{Ritt1950}Ritt, J.F. Differential algebras. AMS, New York
(1950)

\bibitem{ro99} Roitman, M.: On the free conformal and vertex
algebras.  J. Algebra {\bf217},  496-527 (1999).


\bibitem{Rosso89}Rosso,  M.: An~analogue of the Poincare-Birkhoff-Witt theorem
and the universal R-matrix of
$U_q(sl(N+1))$.
 Commun. Math. Phys.
{\bf124},  307-318  (1989)

\bibitem{Sch63}
 Sch\"utzenberger, M.P., Sherman, S.:  On a~formal product over the
conjugate classes in a~free group.  J. Math. Anal. Appl. {\bf7},
 482-488 (1963)


\bibitem{Sch65}Sch\"utzenberger, M.P.:  On a~factorization of free
monoids. Proc. Am. Math. Soc. {\bf16}, 21-24 (1965)


\bibitem{Se94}Segal, D.:  Free left-symmetric algebras and an~analogue
of the Poincar$\acute{e}$-Birkhoff-Witt Theorem.  J. Algebra
 {\bf164},  750-772 (1994)


\bibitem{S53a} Shirshov, A.I.: Some problems in the theory of non-associative
rings and algebras, Candidate of Science Thesis, Moscow State
University, 1953. see http://math.nsc.ru/LBRT/a1/ShirshovPhD.djvu


\bibitem{Sh53}Shirshov, A.I.:  Subalgebras of free Lie algebras. Uspekhi Mat.
Nauk   {\bf8}, (3)  173 (1953)

\bibitem{Sh53b}Shirshov, A.I.: On the representation of Lie rings in associative rings.  Uspekhi Mat.
Nauk N. S.  {\bf8}, (5)(57)  173-175 (1953)

\bibitem{Sh54}Shirshov, A.I.:  Subalgebras of free commutative and free anticommutative
algebras. Mat. Sb.   {\bf4}, (1)  82-88 (1954)

\bibitem{Sh58}Shirshov, A.I.:  On free Lie rings.  Mat. Sb. {\bf45}, (2) 113-122
(1958)

\bibitem{Sh58b}Shirshov, A.I.:  Some problems in the theory of rings that are
nearly associative. Uspekhi Mat. Nauk {\bf13},  no. 6 (84) 3-20
(1958)

\bibitem{Sh62c}Shirshov, A.I.:  On the bases of a~free Lie algebra.
Algebra   Logika  {\bf1}, (1) 14-19 (1962)





\bibitem{Sh62a}Shirshov, A.I.:  Some algorithmic problem for
$\varepsilon$-algebras. Sibirsk. Mat. Zh. {\bf3},  132-137 (1962)


\bibitem{Sh62b}Shirshov, A.I.:  Some algorithmic problem for Lie
algebras. Sibirsk. Mat. Zh. {\bf3}, (2) 292-296  (1962); English
translation in SIGSAM Bull.   {\bf33},   3-6 (1999)

\bibitem{Sh62d}Shirshov, A.I.:  On a~hypothesis in the theory of Lie
algebras. Sibirsk Mat. Zh. {\bf3}, (2) 297-301 (1962)

\bibitem{Shirshov-Selected} Selected works of
A.I. Shirshov. Eds. Bokut, L.A., Latyshev, V., Shestakov, I.,
Zelmanov, E., Trs. Bremner, M., Kochetov, M.,  Birkh\"auser, Basel,
Boston, Berlin (2009)


\bibitem{Tits} Tits, J.: Le probl$\grave{e}$me des mots dans les groupes de Coxeter. Symposia Math. {\bf1},
 175-185 (1968)


\bibitem{Turing} Turing, A.M.: The word problem in semi-groups with
cancellation.  Ann.  Math. {\bf52}, 191-505  (1950)


\bibitem{Ufnarovski95}Ufnarovski, V.A.:  Combinatorial and asymptotic methods in algebra, Algebra
 VI (57) (1995)  1-196. Springer, Berlin Heidelberg New York.


\bibitem{Ufnarovski98}Ufnarovski, V.A.: Introduction to noncommutative  Gr\"{o}bner bases theory.
 Buchberger, Bruno (ed.) et al., Gr\"{o}bner bases and
applications.  Cambridge: Cambridge University Press. Lond. Math.
Soc. Lect. Note Ser. 251, 259-280 (1998)


\bibitem{Umirbaev84}Umirbaev, U.U.: Equality problem for center-by-metabelian Lie
algebras.   Algebra Logic {\bf23}, 209-219 (1984)



\bibitem{Umirbaev93a}Umirbaev, U.U.:
The occurrence problem for Lie algebras.  Algebra Logic {\bf32},
No.3, 173-181 (1993)

\bibitem{Umirbaev93b}Umirbaev, U.U.: Algorithmic problems in associative algebras.
  Algebra Logic {\bf32}, No.4, 244-255
(1993)






\bibitem{Viennot78}Viennot, G.: Algebras de Lie libres et monoid libres. Lecture Notes in
Mathematics, 1978, No.691.


\bibitem{Witt56}  Witt, E.: Die Unterringe der freien Lienge Ringe.
Math. Z.  {\bf64},  195-216 (1956)

\bibitem{Wu1978}Wu, W.-T.:  On the decision problem and the mechanization of
 theorem proving in elementary geometry. Scientia Sinica {\bf21},
 157-179 (1978)


\bibitem{Ya} Yamane, I.: A~Poincare-Birkhoff-Witt theorem for
quantized universal enveloping algebras of type
$A_N$.
 Publ.  RIMS.
Kyoto Univ.  {\bf25},  503-520 (1989)


\bibitem{Abdu2}Yunus, G., Obul, A.:  Gr\"{o}bner--Shirshov basis of quantum group
of type
$D_4$.
 Chin. Ann. Math. {\bf32}, B(4) 581-592 (2011)

\bibitem{Abdu4}Yunus, G.,  Gao,  Z.Z., Obul, A.: Gr\"{o}bner--Shirshov bases of quantum
groups.  Algebra Colloq. to appear.

\bibitem{Xia-Ming}Zhang, X., Jiang M.: On Post's and Markov's examples of semigroups with unsolvable word problem.
Southeast Asian Bull. Math. {\bf37}, 465 -473 (2013)


\bibitem{Zelmanov92}Zelmanov, E.:  Nil rings and periodic groups.   KMS
Lecture Notes in Mathematics. Seoul: Korean Mathematical Society, x,
79 p. (1992)

\bibitem{Zhou-Winkler}Zhou, M., Winkler,  F.: Gr\"obner bases in
difference-differential modules and difference-differential
dimension polynomials. Science in China, Ser. A: Mathematics
{\bf51},  1732-1752  (2008)


\bibitem{Zhukov}Zhukov, A.I.: Complete systems of defining relations in noassociative
algebras. Mat. Sbornik  {\bf69}, (27) 267-280 (1950)


\bibitem{MikhalevZ97}Zolotykh, A.A., Mikhalev, A.A.: A~complex of algorithms for
computations in Lie superalgebras.  Program. Comput. Softw. {\bf23},
No.1, 8-16 (1997)

\bibitem{MikhalevZ98}Zolotykh, A.A., Mikhalev, A.A.: Algorithms for construction of
standard  Gr\"{o}bner--Shirshov bases of ideals of free algebras over
commutative rings.   Program. Comput. Softw. 24, No.6, 271-272
(1998)


\end{thebibliography}
\end{document}